\documentclass[11pt]{amsart}
\usepackage{latexsym}
\usepackage{amscd, amsfonts, eucal, mathrsfs, amsmath, amssymb, amsthm}
\usepackage[OT2,OT1]{fontenc}
\input xy
\xyoption{all}

\usepackage{vmargin}

\setpapersize{USletter}
\newcommand{\Caldararu}{C\u ald\u araru}

 \DeclareFontFamily{U}{wncy}{}
 \DeclareFontShape{U}{wncy}{m}{n}{<->wncyr10}{}
 \DeclareSymbolFont{mcy}{U}{wncy}{m}{n}
 \DeclareMathSymbol{\Sha}{\mathord}{mcy}{"58}

\newcommand{\field}[1]{\mathbf #1}
\newcommand{\mf}[1]{\mathfrak #1}
\newcommand{\mc}[1]{\mathcal #1}
\newcommand{\ms}[1]{\mathscr #1}
\newcommand{\widebar}[1]{\overline{#1}}

\newcommand{\R}{\field R}
\renewcommand{\L}{\field L}
\newcommand{\C}{\field C}
\newcommand{\F}{\field F}
\newcommand{\Z}{\field Z}
\newcommand{\Q}{\field Q}

\newcommand{\simto}{\stackrel{\sim}{\to}}

\renewcommand{\phi}{\varphi}

\renewcommand{\O}{\ms O}
\renewcommand{\hom}{\operatorname{Hom}}
\newcommand{\shom}{\ms H\!om}

\newcommand{\rshom}{\mathbf{R}\shom}
\newcommand{\saut}{\ms A\!ut}
\DeclareMathOperator{\sh}{Sh}
\DeclareMathOperator{\rhom}{\operatorname{{\bf R}Hom}}
\newcommand{\send}{\ms E\!nd}
\newcommand{\rsend}{\mathbf{R}\ms E\!nd}

\newcommand{\spec}{\operatorname{Spec}}

\renewcommand{\P}{\field P}
\newcommand{\A}{\field A}

\DeclareMathOperator{\Pic}{Pic}
\newcommand{\td}{\operatorname{Td}}
\DeclareMathOperator{\chern}{ch}
\DeclareMathOperator{\muhat}{\widehat{\mu}}
\DeclareMathOperator{\supp}{Supp}
\DeclareMathOperator{\pr}{pr}
\DeclareMathOperator{\quot}{Quot}
\DeclareMathOperator{\ass}{Ass}

\DeclareMathOperator{\ob}{\operatorname{Ob}}
\newcommand{\finitemod}{\operatorname{Mod}^{\operatorname{finite}}}

\newcommand{\m}{\boldsymbol{\mu}}

\newcommand{\G}{\field G} 

\newcommand{\etale}{\operatorname{\acute{e}t}}

\renewcommand{\H}{\operatorname{H}}
\newcommand{\bH}{\mathbf{H}}

\newcommand{\GL}{\operatorname{GL}}
\newcommand{\PGL}{\operatorname{PGL}}

\DeclareMathOperator{\ext}{\operatorname{Ext}}
\DeclareMathOperator{\sext}{\ms E\!xt}

\newcommand{\fppf}{\operatorname{fppf}}

\DeclareMathOperator{\dimhom}{hom}
\DeclareMathOperator{\dimext}{ext}

\newcommand{\ch}{\operatorname{char}}
\DeclareMathOperator{\per}{per}
\DeclareMathOperator{\ind}{ind}

\DeclareMathOperator{\D}{\operatorname{\bf D}}

\renewcommand{\[}{[\!\hspace{0.03em}[}
\renewcommand{\]}{]\!\hspace{0.03em}]}

\renewcommand{\)}{)\!\hspace{0.03em})}
\DeclareMathOperator*{\tensor}{\otimes}
\DeclareMathOperator*{\ltensor}{\stackrel{\field L}{\otimes}}
\newcommand{\im}{\operatorname{im}}

\DeclareMathOperator{\Tr}{\operatorname{Tr}}
\DeclareMathOperator{\rk}{\operatorname{rk}}

\newcommand{\inj}{\hookrightarrow}
\newcommand{\id}{\operatorname{id}}

\newcommand{\xto}{\xrightarrow}

\DeclareMathOperator{\End}{\operatorname{End}}
\DeclareMathOperator{\aut}{\operatorname{Aut}}
\DeclareMathOperator{\isom}{\operatorname{Isom}}
\DeclareMathOperator{\M}{\operatorname{M}}
\DeclareMathOperator{\Br}{\operatorname{Br}}

\DeclareMathOperator*{\ctensor}{\widehat{\tensor}}

\newcommand{\invlim}{\varprojlim}

\newcommand{\Tw}{\mathbf{Tw}}
\DeclareMathOperator{\mTw}{Tw}
\newcommand{\Sh}{\mathbf{Sh}}
\DeclareMathOperator{\mSh}{Sh}
\DeclareMathOperator{\B}{\operatorname{B\!}}

\DeclareMathOperator{\expdim}{exp dim}

\newtheorem{lem}{Lemma}[subsubsection]
\newtheorem{thm}[lem]{Theorem}
\newtheorem*{theorem}{Theorem}

\newtheorem{prop}[lem]{Proposition}

\newtheorem{cor}[lem]{Corollary}

\theoremstyle{definition}
\newtheorem{defn}[lem]{Definition}

\newtheorem{example}[lem]{Example}
\newtheorem{hyp}[lem]{Hypothesis}

\newtheorem{para}[lem]{\indent}
\newtheorem{situation}[lem]{Situation}

\theoremstyle{remark}
\newtheorem{remark}[lem]{Remark}
\newtheorem{remarks}[lem]{Remarks}
\newtheorem{notn}[lem]{Notation}
\newtheorem{notation}[lem]{Notation}
\newtheorem{ques}[lem]{Question}

\newenvironment{rem}{\begin{remark}}{\nolinebreak{\hfill $\blacklozenge$} \end{remark}}

\newenvironment{exa}{\begin{example}}{\nolinebreak{\hfill $\lozenge$} \end{example}}

\author{Max Lieblich}
\thanks{Work on this paper was partially supported by a National 
Science Foundation Graduate 
Fellowship, a Clay Mathematics Institute Liftoff Fellowship, and a 
National Science Foundation Postdoctoral Fellowship.}
\title{Moduli of twisted sheaves}
\address{Department of Mathematics, Princeton University, Princeton NJ 08544}
\email{{\tt lieblich@math.princeton.edu}}
\date{}
\begin{document}

\bibliographystyle{plain}

\begin{abstract} We study moduli of semistable twisted sheaves
on smooth proper morphisms of algebraic spaces.  In the 
case of a relative curve or surface, we prove results on the structure 
of these spaces.  For curves, they are essentially isomorphic to spaces 
of semistable vector bundles.  In the case of surfaces, we show (under 
a mild hypothesis on the twisting class) that the spaces are asympotically 
geometrically irreducible, normal, generically smooth, and l.c.i.\ over 
the base.  We also develop general tools necessary for these results: 
the theory of associated points and purity of sheaves on Artin 
stacks, twisted Bogomolov inequalities, semistability and boundedness 
results, and basic results on twisted Quot-schemes on a surface.
\end{abstract}

\maketitle

\tableofcontents

\section{Introduction}

We have only recently begun to understand the role played in algebraic and arithmetic geometry by twisted sheaves.  Originally studied by mathematical physicists,
most research to date has focused on their derived category.  In his thesis \cite {caldararu}, \Caldararu\ extended the classical construction of the Fourier-Mukai transform to study equivalences of derived categories of twisted sheaves on elliptic fibrations.  This was further extended in \cite {pantev-donagi} to the case of genus 1 fibrations without sections which appear as elements of the Tate-Shafarevich group of a fixed elliptic fibration.  These twisted Fourier-Mukai transforms arise in the presence of non-fine moduli spaces; instead of having a universal sheaf, one has a universal \emph {twisted\/} sheaf.

Implicit in these constructions is a theory of moduli for twisted sheaves.  Recent work (apparently roughly simultaneous with the work described here) has produced such a theory tailored to special cases (in particular, K3 surfaces) \cite {yoshioka3}, and this theory has been used to prove a general conjecture of \Caldararu\ about the relationship between Hodge isometries and twisted Fourier-Mukai transforms for K3 surfaces \cite {huybrechts-stellari}.  The construction of the moduli space in \cite {yoshioka3}, using a definition of stability proposed in \cite{caldararu2}, is easily seen to be equivalent to the classical construction of Simpson in the case of moduli of modules.

We develop in this paper a general theory of moduli of twisted sheaves on algebraic spaces, and then apply it to study twisted sheaves on curves and surfaces.  Applications to geometry and arithmetic are taken up in \cite {pgl-bundles} and \cite {period-index-paper}.  The crucial observation is that by thinking of twisted sheaves as sheaves on certain stacks $\ms X$, one can carry out an analysis very similar to the analysis of sheaves on varieties.  Using the rational Chow theory of $\ms X$ and Riemann-Roch theorems for representable morphisms of Deligne-Mumford stacks, one can define a notion of semistability for twisted sheaves.  With this in hand, it is possible to ``twist'' classical tools such as Harder-Narasimhan filtrations, elementary transforms, and determinants, and arrive at constructions and results familiar from the classical case.  This is most interesting in the case of surfaces, where we prove the following main structure theorem.  Let $X \to S$ be a projective relative surface with smooth connected geometric fibers, and let $\ms X \to X$ be a $\m _ n$-gerbe.  (A review of gerbes and twisted sheaves may be found in sections \ref {S:sheaves and gerbes on stacks} and \ref{S:twisted sheaves} below.)

\begin{theorem}
The stack of semistable $\ms X$-twisted sheaves is an Artin stack locally of finite presentation over $S$.  If $\ms X$ is optimal then the substack of twisted sheaves of fixed determinant and sufficiently large discriminant is generically smooth, normal, geometrically irreducible, and l.c.i.\ along the stable locus.
\end{theorem}
The condition of ``optimality'' is added to make the results characteristic-free; it is relatively clear from our methods of proof that this hypothesis can be eliminated in characteristic 0, and recent work of Langer in the untwisted case should help eliminate this hypothesis in arbitrary characteristic (although the methods involved are slightly different from our own).  Since we have not carried out either exercise, we will only give the proof here in the optimal case.

This theorem lives at the junction of several roads in algebraic geometry.  First, it is a ``case study'' in the moduli of sheaves on tame smooth Deligne-Mumford stacks.  In fact, many of the techniques of this paper -- the ``geometric Hilbert polynomial,'' boundedness results, etc. -- can be carried over to the general case (in preparation); the case of gerbes is then seen as the most tractable case of a more general theory.  Second, by a rigidification mechanism, one can use our methods to study compactified moduli spaces of $\PGL _ n$-bundles (or Azumaya algebras, for the arithmetically minded).  The resulting statements give a first-order algebraic approximation to results of Taubes on the stable topology of the space of connections on a fixed smooth bundle.  These matters will be taken up in detail in \cite {pgl-bundles}.

The fact that $\PGL _ n$-bundles are the same thing as Azumaya algebras relates the constructions of this paper to classical problems about the Brauer groups of function fields.  In particular, our techniques permit more efficient and conceptual proofs of the basic facts about Brauer groups of schemes, including the results of Gabber's thesis concerning the relation between the Brauer group and the cohomological Brauer group.  (In fact, de Jong's recent proof of Gabber's theorem on the cohomological Brauer group of a quasi-projective scheme uses techniques similar to those of this paper.) Moreover, the structure theory for moduli spaces of twisted sheaves gives new results about the period-index problem for surfaces over finite and local fields, generalizing well-known recent results of de Jong for surfaces over algebraically closed fields.  In particular, one can prove the following.

\begin{theorem}
Let $X$ be a smooth projective geometrically connected surface over a field $k$ and $\alpha \in \Br (X)$ a Brauer class of order prime to the characteristic of $k$.
\begin{enumerate}

\item If $k$ is algebraically closed, then $\per (\alpha) = \ind (\alpha)$.

\item If $k$ is finite, then $\ind (\alpha) \mid \per (\alpha) ^ 3$.  If in addition $\alpha$ is unramified, then $\ind (\alpha) = \per (\alpha)$.

\item If $k$ is local, $\alpha$ is unramified, and $X$ has smooth reduction, then $\ind (\alpha) \mid \per (\alpha) ^ 2$.  If in addition $\alpha$ is unramified on a smooth model of $X$ over the integers of $k$, then $\ind (\alpha) = \per (\alpha)$.
\end{enumerate}
\end{theorem}
\noindent The first statement is the well-known result of de Jong; the others are new.  These ideas are discussed in \cite {period-index-paper}.

The results described here stand as yet another example of how the use of stack-theoretic methods can clarify and extend classical results, and suggest new approaches and connections.  This is the overarching philosophy of this work, and we hope that the reader will take this, if nothing else, away from this paper.

As this is a young field, we provide relatively 
complete foundations for the abstract theory of twisted sheaves.  We devote all of section 
\ref{S:twisted sheaves big} to a development of the algebraic theory 
of twisted sheaves, including general nonsense on their deformation 
theory, as well as a theory of associated points and purity for 
sheaves on Artin stacks.  We end the section with a discussion of 
Riemann-Roch theorems for gerbes and the basic properties of Quot 
spaces.  In section \ref{S:moduli} we show that the stack of semistable 
twisted sheaves is an Artin stack by using Artin's representability 
theorem.  Finally, section \ref{S:curves and surfaces} is devoted to studying the 
resulting stacks of twisted sheaves on curves and surfaces.  In the 
last subsection, we prove twisted analogues of O'Grady's results on 
asymptotic properties of the moduli spaces, including the first theorem above.

\section*{Acknowledgements}

This paper represents a sizable chunk of my thesis \cite{mythesis}.  
I am indebted to my advisor, Aise Johan de 
Jong, for many fruitful discussions and crucial suggestions.  I would also like to thank 
Brian Conrad, Jacob Lurie, and Jason Starr for helpful 
conversations and comments, and Martin Olsson for thoroughly reading 
the thesis from which this paper is derived and offering numerous 
corrections and improvements.

It has come to my attention that Stuhler and Hoffmann have 
recently obtained some of the results of this paper independently 
\cite{stuhler-hoffmann}.

\section*{Notation}

Following standard conventions, we use $=$ for canonical 
isomorphisms.  

Every locally free sheaf is assumed to have finite rank 
everywhere.

As should be universal by now, ``algebraic stack'' will mean 
``algebraic stack in the sense of Artin.''  Deligne-Mumford stacks 
will be called ``DM stacks.''  All algebraic stacks are quasi-separated, as is any base scheme 
appearing in this paper.

In order to prevent psychological problems, when given a topos $X$, we 
will write $U$ in place of $X_{/U}$ to stand for the restriction of 
$X$ to the object $U\in X$.  For the sake of intuition, we will also interchangeably 
refer to ``sheaves on $X$'' and ``objects in $X$'' depending upon the 
context.

Following Huybrechts and Lehn \cite{h-l}, we use the notation 
``$\dimhom$'' for ``$\dim\hom$'' and 
``$\dimext$'' for ``$\dim\ext$.''  In general, we have tried to keep 
notations in common with their book when treating the twisted 
analogues of classical theorems so beautifully discussed there.

We will occasionally call a stack \emph{quasi-proper\/} if it is universally closed over Noetherian base spaces but is not necessarily separated.  This arises quite often in the theory of stacks, e.g.\ when dealing with strictly semistable objects -- GIT stacks for example can only reasonably be expected to be separated along the stable locus.

There is one pedantic grammatical convention we adopt which we hope will spread: a number 
with mathematical meaning is always written as a numeral, occasionally 
in contradiction to accepted rules of grammar (e.g., 
``rank 1'', ``characteristic 0'').  

\section{Twisted sheaves}\label{S:twisted sheaves big}
\subsection{Preliminaries: twisted sheaves on ringed topoi}\label{S:prelim}
In this section, we lay the foundations for the theory of twisted 
sheaves on algebraic spaces and stacks.  The reader will note that much of the first three 
sections is written in the language of ringed topoi.  We encourage 
those uncomfortable with this notion to substitute ``ringed site'' or 
even ``ringed space'' for ``ringed topos''; the exposition will remain 
more or less the same after this substitution (but the reader should 
note that sites larger than 
the Zariski site of a scheme are essential for the theory to actually be interesting).  
One reason to write in this 
degree of generality is to make the theory apply to algebraic stacks, 
where one can only really understand the theory of sheaves from the 
topos-theoretic point of view.

In order to link Giraud's ideas with subsequent developments in the 
theory of algebraic stacks \cite{l-mb}, we review foundations on the 
sites associated to a stack, sheaves on those sites, and classifying 
topoi associated to gerbes on (ringed) topoi.  We only consider stacks 
in groupoids here; the task of extending the results to 
stacks in arbitrary (small) categories is left to the (odd) reader.  (It 
will primarily consist of adding the word ``Cartesian'' in a few 
places.)

\subsubsection{Sheaves and gerbes on stacks}\label{S:sheaves and gerbes on stacks}
Let $X$ be a topos and $F:\ms S\to X$ a stack on $X$.  The topology 
on $X$ naturally induces a topology on $\ms S$.

\begin{defn} The \emph{site of $\ms S$\/}, denoted $\ms S^{s}$, has as underlying category
    \begin{enumerate}
	\item[] Objects: morphisms $f:S\to\ms S$ of fibered categories over 
	$X$, where $S$ ranges over all sheaves on (=objects of) $X$
	\item[] Morphisms: a morphism from $f:S\to\ms S$ to $g:S'\to\ms S$ is 
	a pair $(\phi,\psi)$ where $\phi:S\to S'$ is a morphism in $X$ and 
	$\psi:f\simto g\circ\phi$ is a 2-isomorphism.
    \end{enumerate}
A covering is given by a morphism $(\phi,\psi)$ with $\phi$ a covering.
\end{defn}

\begin{defn} The \emph{classifying topos\/} of $\ms S$, denoted 
$\widetilde{\ms S}$, is the topos of sheaves on the site of $\ms S$.
\end{defn}

There is a morphism of topoi $\pi:\widetilde{\ms S}\to X$: 
given a sheaf $\ms F$ on $X$, one gets a sheaf $\pi^{\ast}\ms F$ on 
$\widetilde{\ms S}$ by assigning to $f:S\to \ms S$ the object $\ms 
F(S)$.  
The obvious exactness properties show that this is the pullback of a 
morphism of topoi.  In particular, when $X$ is ringed, say by $\ms 
O$, $\widetilde{\ms S}$ is naturally ringed by $\pi^{\ast}\ms O$.  

\begin{rem}\label{R:we agree with Giraud} The reader can easily check that 
the formation of the classifying 
topos is functorial, and that our description agrees with Giraud's 
original definition:
$\widetilde{\ms S}=\operatorname{Cart}(\ms 
S,\operatorname{Fl(X)})$ \cite[\S5.1]{giraud}.
\end{rem}

\begin{para} Given a stack $\ms S\to X$, there is an associated 
stack $\ms I(\ms S)\to\ms S$ called the \emph{inertia stack\/}.  To give the stack $\ms I(\ms S)$ it is enough to describe its sheaf of sections.
\end{para}

\begin{defn} The assignment $(f:S\to\ms S)\mapsto\saut(f)$ is a sheaf. 
The corresponding stack is denoted $\ms I(\ms S)\to\ms 
S$ and called the \emph{inertia stack\/} of $\ms S$.
\end{defn}

\begin{lem} For any stack $\ms S$, $\ms I(\ms S)=\ms S\times_{\ms S\times\ms S}\ms S$.
\end{lem}

Thus, when $X$ is the topos of sheaves on the big \'etale topology on affine 
schemes over a fixed base $B$ and $\ms S$ is an algebraic stack, then 
one easily sees that $\ms I(\ms S)$ is also an algebraic stack and 
the morphism $\ms I(\ms S)\to\ms S$ is representable, quasi-compact, 
and separated.  We leave the following standard lemmas to the 
reader.

\begin{lem}\label{L:inertia is functorial} Given a 1-morphism $f:\ms 
S\to\ms S'$ of stacks, there is an induced map $\ms I(\ms S)\to 
f^{\ast}\ms I(\ms S)$ in $\widetilde{\ms S}$.
\end{lem}

\begin{lem}\label{L:universal inertia action} Given any sheaf $F$ on $\widetilde{\ms S}$, there is a 
natural right group action $\mu:F\times\ms I(\ms S)\to F$.
\end{lem}

\begin{para}\label{Par:gerbes} We will be concerned throughout this 
paper with gerbes.
\end{para}
\begin{defn}\label{D:gerbe defn} The stack $\ms S$ is a \emph{gerbe\/} on $X$ if 
    \begin{enumerate}
	\item For any $U\in X$ there exists a covering $U'\to U$ such 
	that $\ms S_{U'}\neq\emptyset$.
    
	\item For any $U\in X$ and any $s,s'\in\ms S_{U}$, there 
	exists a covering $U'\to U$ such that $s|_{U}$ is isomorphic 
	to $s'|_{U'}$.
    \end{enumerate}
    
\end{defn}
In looser language, $\ms S$ has local sections everywhere and any two 
sections are locally isomorphic.  There is a ``moduli-theoretic'' 
interpretation of this definition.

\begin{defn}\label{D:associated sheaf} The \emph{sheaf associated to $\ms S$\/}, denoted 
$\sh(\ms S)$, is the sheafification of 
the presheaf whose sections over $U\in X$ are isomorphism classes of 
objects in the fiber category $\ms S_{U}$.
\end{defn}

\begin{lem}\label{L:gerbe if sheaf is base} The stack $\ms S$ is a gerbe on $X$ if and only if the 
natural map $\sh(\ms S)\to e_{X}$ is an isomorphism in $X$.
\end{lem}
Here $e_{X}$ denotes the final object of the topos $X$.  This will 
often be written as $X$ by abuse of notation.
\begin{proof} Suppose $\ms S$ is a gerbe.  By functoriality of the natural 
map, it is enough to demonstrate the claim when $\ms S$ has a global 
section $\sigma$ over $X$.  But then every local section is locally isomorphic 
to $\sigma$, hence $\sh(\ms S)$ is a singleton and the natural map is 
an isomorphism.

Suppose conversely that $\sh(\ms S)\to e_{X}$ is an isomorphism.  In 
particular, $$\sh(\ms S)(U)=\{\emptyset\}$$ for any $U\in X$.  By the 
definition of $\sh(\ms S)$ and of sheafification, this says precisely that 
conditions 1 and 2 of 
\ref{D:gerbe defn} is satisfied.  
\end{proof}

\begin{lem}\label{L:descend sheaf on gerbe} If $\ms S\to X$ is a gerbe and $F$ is a 
sheaf on $\ms S$ such 
that the inertia action $F\times\ms I(\ms S)\to F$ is trivial, then 
$F$ is naturally the pullback of a unique sheaf on $X$ up to isomorphism.
\end{lem}
\begin{proof} We claim that $\pi^{\ast}\pi_{\ast}F\to F$ is an 
isomorphism.  To verify this, it suffices to work locally on $X$, so 
we may assume that $\ms S$ has a section.  One can then check using 
the hypothesis on the action that the pullback of $F$ along this 
section equals the pushforward of $F$ along the structure morphism.  
The result follows.
\end{proof}

\begin{rem} This holds more generally when $X$ is the coarse moduli space of a 
Deligne-Mumford stack $\ms S$ (with the action of inertia being 
studied in the big \'etale topology), but the proof is slightly more 
difficult: it follows without too much difficulty from the \'etale 
local structure of the stack as a finite group quotient 
stack \cite[\S6]{l-mb}, \cite{keel-mori}.
\end{rem}

\begin{lem} If $\pi:\ms S\to X$ is a gerbe and $\ms I(\ms S)$ is an abelian 
sheaf on $\ms S$, then there is an abelian sheaf $A$ on $X$ and an 
isomorphism $\widetilde{\pi}^{\ast}A\cong\ms I(\ms S)$ as objects 
of $\widetilde{\ms S}$.
\end{lem}
\begin{proof} This is an application of \ref{L:descend sheaf on gerbe}.
\end{proof}

\subsubsection{Twisted sheaves}\label{S:twisted sheaves}

Let $(X,{\ms O})$ be a ringed topos, $A$ a sheaf of commutative groups on 
$X$, and $\chi:A\to\G_{m}$ a character.

\begin{defn} An \emph{$A$-gerbe\/} on $X$ is a gerbe $\ms S\to X$ 
along with an isomorphism $A_{\ms S}\simto\ms I(\ms S)$ in 
$\widetilde{\ms S}$.
\end{defn}
When $A$ is non-commutative, this definition is not correct.  One must 
instead choose any isomorphism as in the definition in the category of 
\emph{liens\/} on $\ms S$ rather than the category of sheaves.  The basic 
reason for this may be seen by thinking about the gerbe $\B{G}$ for a 
non-commutative group $G$.  In general, the automorphism group of a 
left $G$-torsor is an inner form of $G$, not $G$ itself.  (The 
stack of liens on a topos is the universal stack receiving a 1-morphism from 
the stack of groups on that topos such that two inner forms naturally 
map to isomorphic objects.)  This is of course described in great 
detail in \cite{giraud}.

Given a cohomology class 
$\alpha\in\H^{2}(X,A)$, there is a corresponding equivalence class of 
$A$-gerbes on $X$.  We will fix such a class $\alpha$ and an 
$A$-gerbe $\ms X\to X$.  The goal of this section is 
to single out a subcategory of sheaves on $\ms X$ which will play a 
fundamental role in what follows.

Given an ${\ms O}_{\ms X}$-module $\ms F$,
the module action $m:\G_{m}\times\ms F\to\ms F$ yields an 
associated right action $m':\ms F\times\G_{m}\to\ms F$ with 
$m'(s,\phi)=m(\phi^{-1},s)$.  This will always be called the 
\emph{associated right action\/}.

\begin{defn} A \emph{$d$-fold $\chi$-twisted sheaf\/} on $\ms X$ is an ${\ms O}_{\ms 
X}$-module $\ms F$ such that the natural action $\mu:\ms F\times A\to\ms 
F$ given by the $A$-gerbe structure makes the diagram
$$\xymatrix{\ms F\times A\ar[r]\ar[d]_{\chi^{d}} & \ms F\ar[d]^{\id}\\
\ms F\times\G_{m}\ar[r]^{m'} & \ms F}$$
commute, where $\chi^{d}(s)=\chi(s)^{d}$.  A $1$-fold twisted sheaf 
will be called simply a twisted sheaf.
\end{defn}

We will see below that if $\ms X^{(d)}$ is a gerbe representing 
$d\cdot\alpha\in\H^{2}(X,A)$, then twisted sheaves on $\ms X^{(d)}$ 
are equivalent to $d$-fold twisted sheaves on $\ms X$.  Note that $\ms X$-twisted sheaves naturally form a fibered subcategory of the classifying topos $\widetilde{\ms X}$, viewed as a fibered category over $X$ via the natural map $\widetilde{\ms X}\to X$ of topoi.

\begin{prop}\label{P:twisted sheaves form a stack} The fibered category of $\ms X$-twisted sheaves is a naturally a stack $X$.
\end{prop}
\begin{proof} This follows from the fact that the condition that a sheaf on $\ms X$ be an $\ms X$-twisted sheaf is local on the site of $\ms X$ along with the fact that any morphism $Y\to X$ of topoi defines a natural stack on $X$ by restriction.
\end{proof}

\begin{notn}\label{N:curly T} Given a ($\G_m$- or $\m_n$-) gerbe $\ms X$ on a morphism of ringed topoi $X\to S$, the stack of $S$-flat quasi-coherent $\ms X$-twisted sheaves locally of finite presentation will be denoted $\ms T_{\ms X/S}$.
\end{notn}

An attentive reader may have noticed that the morphism $\chi$ yields a 
``change of structure group'' for the gerbe $\ms X$.

\begin{prop}\label{L:comparison of gerbes} Let $f:A\to B$ be a morphism of abelian sheaves and 
$\alpha\in\H^{2}(X,A)$ a cohomology class with direct image
$f_{\ast}(\alpha)=\beta\in\H^{2}(X,B)$.  Given gerbes $\ms X_{\alpha}$ and $\ms 
X_{\beta}$ representing $\alpha$ and $\beta$, there is a 1-morphism 
$F:\ms X_{\alpha}\to\ms X_{\beta}$ over $X$ such that for any 
section $\sigma:S\to\ms X_{\alpha}$, the induced morphism 
$A_{S}=\aut(\sigma)\to\aut(F(\sigma))=B_{S}$ is $f_{S}$.  Given a character $\chi:B\to\G_m$, this induces an identification of the stack of $\chi$-twisted sheaves with the stack of $\chi\circ f$-twisted sheaves.
\end{prop}
\begin{proof} The existence of the morphism of stacks $\ms 
X_{\alpha}\to\ms X_{\beta}$ is part of 
Giraud's theory of non-abelian cohomology.  One can also see this 
explicitly by re-expressing the gerbe as a stack of twisted torsors and using the natural contraction of torsors along a group homomorphism to define the morphism of gerbes.  
This is done in \cite{mythesis}.  To see that the stacks of twisted sheaves are identified, it is enough to prove that for the morphism $\chi:B\to\G_m$, pullback along $\chi$ yields a 1-isomorphism from the stack of $\id$-twisted sheaves to the stack of $\chi$-twisted sheaves.  Since these are both stacks, it suffices to prove this locally on $X$, so we may assume that in fact both gerbes admit global sections, say $\sigma_{\alpha}$ and $\sigma_{\beta}$.  It is easy to see that there is an invertible sheaf $\ms L$ on $\ms X$ such that for any $\chi$-twisted sheaf $\ms F$ on $\ms X_{\beta}$, there is a natural isomorphism $\sigma_{\alpha}^{\ast}F^{\ast}\ms F\simto\ms L\tensor\sigma_{\beta}^{\ast}\ms F$.  It is therefore enough to prove that for any $A$-gerbe $\ms X$ with a section $\sigma$ and any character $\chi:A\to\G_m$, the pullback functor $\sigma^{\ast}$ defines an equivalence of the stack of $\chi$-twisted sheaves with the stack of modules on $\ms X$.  But any such $\ms X$ is isomorphic to $\B{A}$, so the stack of all modules on $\ms X$ is equivalent to the stack of $\ms O_X$-modules with a right $A$-action, and the stack of $\chi$-twisted sheaves is equivalent to the stack of $A$-equivariant $\ms O_X$-modules such that $A$ acts via $\chi$.  It is clear that this last category is equivalent to the category of modules by simply forgetting the $A$-action, which is precisely the pullback functor.
\end{proof}

Under this identification, the stack of $\chi$-twisted sheaves is 
identified with the stack of $\iota$-twisted sheaves, where 
$\iota:\im\chi\inj\G_{m}$ is the natural inclusion of the image.  
Thus, we have gained very little but canonicity by our formalism.  
However, one might in the future try something similar when $\ms S$ 
is not a gerbe and the 
inertia stack $\ms I(\ms S)$ is \emph{not\/} constant, in which case 
this setup is the correct one.  Such a study is related to moduli of 
ramified Azumaya algebras and the ramified period-index problem.  
These issues will be explored in future work.

\subsubsection{Comparison with the formulation of 
\Caldararu}\label{S:cald comp}
We explain in this section how our definition of twisted sheaves 
squares with that used by \Caldararu\ in \cite{caldararu}.  The reader will note 
that his formulation seems more ``user-friendly.''  We hope to make 
clear below, especially in our discussion of deformations and 
obstructions, why the more abstract approach is essential.

Throughout, we retain the notation of the previous section: $(X,{\ms O})$ 
is a ringed topos and $\chi:A\to\G_{m}$ is a character of a sheaf of 
commutative groups.  Let 
$\alpha\in\H^{2}(X,A)$ be a fixed cohomology class.  By a theorem of 
Verdier \cite[Expos\'e V.7]{sga4-2}, there is a hypercovering $U_{\bullet}\to X$ and a cocycle 
$a\in\Gamma(U_{2},A)$  
which represents $\alpha$ in cohomology.  We fix such a 
representative in this section.  We also fix a choice 
of $A$-gerbe $\ms X$ representing the cohomology class $\alpha$.

\begin{defn}[\Caldararu] A $\chi$-twisted sheaf on $X$ is a pair 
$(F,g)$, where $F$ is an ${\ms O}_{U_{0}}$-module and $g:(\pr_1^{U_1})^{\ast}F\simto(\pr_0^{U_1})^{\ast}F$ 
is a gluing datum on $U_{1}$ such that $\delta g\in\aut((\pr_0^{U_2})^{\ast}F)$ 
equals the cocycle $\chi(a)$.
\end{defn}
We will (temporarily) call such an object a \Caldararu-$\chi$-twisted 
sheaf.

\begin{exa} Suppose $A=\G_{m}$, $\chi=\id$, and $X$ is a complex analytic space.  
We may take the hypercovering $U_{\bullet}$ to be the \v Cech 
hypercovering generated by an open covering of $X$, i.e., we may 
replace $U_{\bullet}$ by an open covering $\{U_{i}\}$ of $X$.  Then 
a $\chi$-twisted sheaf on $X$ is given by 
\begin{enumerate}
    \item a sheaf of modules $\ms F_{i}$ on each $U_{i}$

    \item for each $i$ and $j$ an isomorphism of modules $g_{ij}:\ms 
    F_{j}|_{U_{ij}}\simto\ms F_{i}|_{U_{ij}}$

\end{enumerate}
subject to the requirement that on $U_{ijk}$, 
$g_{ik}^{-1}g_{ij}g_{jk}:\ms F_{k}|_{U_{ijk}}\simto\ms 
F_{k}|_{U_{ijk}}$ is equal to multiplication by the scalar 
$a\in\G_m(U_{ijk})$ giving the 2-cocycle. 
\end{exa}

\begin{prop} There is a natural equivalence of fibered categories between 
$\chi$-twisted sheaves and \Caldararu-$\chi$-twisted sheaves.
\end{prop}
Surprisingly, the proof is not obvious.  It is written in gory detail 
in section 2.1.3 of \cite{mythesis}.  It can at times be useful to know that 
these two categories are equivalent, as certain statements are completely obvious in 
one and completely mysterious in the other.  The canonical example of 
this is furnished by quasi-coherent and coherent twisted sheaves on a 
Noetherian gerbe.  One sees easily using the stack-theoretic language 
that any quasi-coherent twisted sheaf is the colimit of its coherent 
subsheaves, whereas from the \Caldararu\ point of view this is far 
from obvious.  We will summarize the important properties of 
quasi-coherent twisted sheaves below; the proofs are exercises and have been omitted or briefly sketched.

It is obvious that the \Caldararu-$\chi$-twisted 
sheaves are equivalent to twisted sheaves for the cocycle in $\G_{m}$ 
induced by $\chi$.  (This is an example of a property which is easier to detect using the \Caldararu\ formalism -- witness the proof of \ref{L:comparison of gerbes}.)  Let $\ms X_{\alpha}\to\ms X_{\chi(\alpha)}$ be as in 
\ref{L:comparison of gerbes}, where $\chi(\alpha)\in\H^{2}(X,\G_{m})$.  

\begin{defn} If $\chi:A\to\G_{m}$ is the natural inclusion of a 
subsheaf (e.g., $\m_{n}$ or $\G_{m}$), a $\chi$-twisted sheaf on $\ms 
X$ will be called an \emph{$\ms X$-twisted sheaf\/}.
\end{defn}

\begin{cor}\label{C:character means nothing} The map $\ms X_{\alpha}\to\ms X_{\chi(\alpha)}$ induces 
by pullback a 1-isomorphism of the stack of $\ms 
X_{\chi(\alpha)}$-twisted sheaves with the stack of $\chi$-twisted 
sheaves. 
\end{cor}

\subsection{The case of a scheme}\label{S:the case of a scheme}

\subsubsection{Quasi-coherent twisted sheaves}

Let $X$ be a scheme, $A$ a group scheme which is faithfully flat and locally of finite 
presentation over $X$, 
$\alpha\in\H^{2}(X_{\fppf},A)$ a flat cohomology class, and $\chi:A\to\G_{m}$ an 
algebraic character.  Fix a gerbe $\ms X$ representing $\alpha$ in the 
big fppf topology on $X$.  When $A$ is smooth, a theorem of 
Grothendieck \cite[Appendix]{grothbrauer2} says that the restriction of $\ms X$ to the (big or small) 
\'etale topos of $X$ is an $A$-gerbe (and this defines an 
isomorphism $\H^{2}(X_{\fppf},A)\simto\H^{2}(X_{\etale},A)$).  (In fact, 
Grothendieck's theorem holds for the cohomology in all 
degrees.)  We recall for the reader some of the basic facts about 
quasi-coherent twisted sheaves.

\begin{lem} The gerbe $\ms X$ is an algebraic stack locally of finite 
presentation over $X$.  If $X$ is quasi-separated and $A$ is finitely presented then $\ms X$ is 
finitely presented.  The scheme $X$ is (locally) Noetherian if and 
only if $\ms X$ is (locally) Noetherian.
\end{lem}

\begin{remark} In our study of twisted sheaves on surfaces (when we 
actually want to say something!), we will take $A=\m_{n}$ with $n$ prime 
to the characteristics of $X$.  In this case, the reader will 
immediately verify that any $A$-gerbe is in fact a DM stack.
\end{remark}

When $X$ is Noetherian we can define quasi-coherent and coherent 
twisted sheaves.

\begin{lem} The obvious morphism $F:\widetilde{\ms 
X}_{\fppf}\to\widetilde{\ms X}_{\etale}$ induces by pullback an 
equivalence of the stacks of quasi-coherent sheaves.  These stacks are naturally equivalent to the stack of twisted sheaves in the lisse-\'etale topos of $\ms X$.  
\end{lem}

The notion of quasi-coherence is also independent of the group chosen.

\begin{lem}\label{L:quasi-coherent is unambiguous} Under the equivalence of \ref{C:character means 
nothing}, quasi-coherent sheaves are taken to quasi-coherent sheaves. 
\end{lem}
\begin{proof} This follows from the fact that a sheaf on an algebraic stack $\ms X$ is 
quasi-coherent if and only if it pulls back to a quasi-coherent 
big \'etale sheaf on any scheme mapping to $\ms X$.  Thus, if $T\to\ms X_{\alpha}$, the compatibility of pullbacks shows that given any $\ms F$ on $\ms X_{\chi(\alpha)}$, we have that $\ms F|_{\ms X_{\alpha}}|_T$ is naturally isomorphic to the pullback along the induced morphism $T\to\ms X_{\chi(\alpha)}$, whence it is quasi-coherent by the assumption on $\ms F$.
\end{proof}

\begin{prop}\label{C:colimit} Suppose $X$ is Noetherian and $A$ is 
group scheme faithfully flat of finite presentation over $X$.  A quasi-coherent $\chi$-twisted 
sheaf is the colimit of its coherent $\chi$-twisted subsheaves.
\end{prop}
This proposition turns out to be quite useful in re-proving the basic 
facts about the Brauer group of a scheme.  We refer the interested 
reader to \cite{period-index-paper}.

In fact, when $A$ is diagonalizable we can split up the category of quasi-coherent 
sheaves into pieces indexed by characters.  Suppose $D$ is a diagonalizable affine group scheme (i.e., the 
Cartier dual of $D$ is a constant finitely generated abelian group).  
Write $C$ for the dual group of $D$, which is the group of 
homomorphisms $D\to\G_{m}$.  Let $\ms X$ be a $D$-gerbe and $\ms F$ a 
quasi-coherent sheaf on $\ms X$.  Given $\chi\in C$, there is a 
$\chi$-eigensheaf $\ms F_{\chi}\subset\ms F$.  The following proposition follows easily from the representation theory of diagonalizable group schemes (and reduction to the case of a trivial gerbe).

\begin{prop}\label{P:q-coh splits}  Suppose $\ms F$ is a 
quasi-coherent sheaf on $\ms X$.  The natural maps induce an 
isomorphism
$$\bigoplus_{\chi\in C(X)}\ms F_{\chi}\simto\ms F.$$
The eigensheaves $\ms F_{\chi}$ are quasi-coherent.
\end{prop}

Let $Y\to X$ be a quasi-compact morphism of schemes and $\ms X$ a $D$-gerbe on $X$.  
Define $\ms Y:=Y\times_{X}\ms X$; this is naturally a $D$-gerbe on 
$Y$.  Denote the morphism $\ms Y\to\ms X$ by $\pi$.

\begin{lem}\label{L:q-coh split is natural} If $\ms F$ is 
a quasi-coherent sheaf on $\ms Y$, then the natural map $\pi_{\ast}\ms 
(\ms F_{\chi})\to\pi_{\ast}\ms F$ identifies $\pi_{\ast}(\ms F_{\chi})$ 
with $(\pi_{\ast}\ms F)_{\chi}$.
\end{lem}

\subsubsection{Gabber's theorem and Morita equivalence}\label{S:gabber's theorems}

Using twisted sheaves, de Jong \cite{dejong-gabber} has recently proven the following 
result of Gabber (vastly generalizing a result of his thesis 
\cite{gabber}).

\begin{thm}\label{T:gabber de jong} If $X$ is a quasi-compact separated 
scheme admitting an ample invertible sheaf then $\Br(X)=\Br'(X)$.
\end{thm}

For the reader unfamiliar with the Brauer group, this result is 
equivalent to the following corollary.

\begin{cor}\label{C:gabber de jong} Given $X$ as in \ref{T:gabber de 
jong} and a $\m_{n}$-gerbe $\ms X\to X$, there is a locally free $\ms 
X$-twisted sheaf $\ms V$ of constant non-zero rank.
\end{cor}

In general, the question of when there exists a locally free $\ms 
X$-twisted sheaf is delicate and interesting question.  It is 
equivalent to $\ms X$ being a quotient stack; the study of this 
question for gerbes is closely related the question of when an arbitrary  
(tame) Deligne-Mumford stack is a quotient stack.  This has been 
studied by Vistoli, Kresch, Hasset, Edidin and others (see 
\cite{vistoli-kresch}, \cite{vistoli-kresch-etc}, \cite{totaro} and 
the references therein).

As a consequence of \ref{C:gabber de jong}, we can rewrite the 
theory of twisted sheaves on such a scheme $X$ in terms of modules 
over the algebra $\send(\ms V)$.

\begin{prop}\label{P:mini-morita} Given the notation of \ref{C:gabber de jong}, the 
functor $\ms W\mapsto\shom(\ms V,\ms W)$ establishes an equivalence of 
fibered categories between $\ms X$-twisted sheaves and right 
$\send(V)$-modules.
\end{prop}
\begin{proof}[Sketch of proof] See Theorem 1.3.7 of \cite{caldararu}.  This is a special case of ``fibered Morita 
equivalence'' which is studied in gory generality (in an arbitrary 
topos with sufficiently many points) in \cite{mythesis}.
\end{proof}

\subsubsection{Deformations and obstructions}\label{S:deformations of twisted sheaves}
Since twisted sheaves are modules in a topos, we can try to apply the deformation theory of 
Illusie to study deformations and obstructions of twisted sheaves.  The condition that 
a deformation preserve the character of the inertial action and that an 
obstruction take this into account makes the situation slightly more 
complicated than in Illusie's bare theory.  We present an alternative 
approach 
to the deformation theory of twisted sheaves, parallel to the 
approach of Grothendieck sketched in \cite{illusie}.  Since extensions of quasi-coherent twisted sheaves in the 
category of all sheaves on a gerbe are well-behaved (and stay 
twisted!), it is easy to see that Illusie's theory (with its attendant 
functorialities) applies perfectly.  Thus, the reader who trusts that 
one can develop the deformation theory of twisted sheaves ``from 
scratch'' in the more complicated cases can skip ahead to the study 
of Artin's conditions for this deformation theory 
\ref{Para:obstruction}.  We note that the non-Illusian approach 
described here may have wider applicability.  An example of this is 
furnished by recent work \cite{moduli-complexes} on constructing an algebraic 
stack of complexes; one can approach the deformation theory of 
objects in the derived category using a derived version of 
Grothendieck's approach which is not clearly related to Illusie's 
approach.

\begin{para} In this 
section, we work in the Abelian category of twisted sheaves.  Thus, 
all $\ext$ groups are computed in this category and not in the larger 
category of all sheaves on $\ms X$.  When we specialize to 
quasi-coherent twisted sheaves, this will no longer matter, as both 
$\ext$ spaces are naturally isomorphic.  For the moment, fix a topos 
$X$ and a $\G_{m}$-gerbe $\ms X$ on $X$.
\end{para}
\begin{lem}\label{L:enough of everything} The category of $\ms X$-twisted sheaves contains enough 
injectives and enough flat objects.
\end{lem}
\begin{proof} Let $U\in X$ be an object over which $\ms X$ splits and 
let $\ms U=\ms X\times_{X}U$, with natural map $f:\ms U\to\ms X$.  
Then $\ms U\cong\B{\G_{m,U}}$ as $U$-stacks.  The usual abstract 
nonsense shows that there are enough twisted injectives and flat objects on $\ms 
U$.  Taking $f_{\ast}$ of injectives and $f_{!}$ of flat objects 
yields the desired result.  The details are left to the reader (or see 
\cite{caldararu}).
\end{proof}

In fact, recent work \cite{locoresc} yields the existence of 
$K$-injective and $K$-flat resolutions in the homotopy category of 
complexes of twisted sheaves.  This justifies the use of all of the 
usual derived functors for twisted sheaves: $\rshom$, $\rhom$, 
$\ltensor$, $\L f^{\ast}$, etc., as well as allowing proofs of the 
assertions below.  We refer the reader unfamiliar with these 
ideas to \cite{spaltenstein} and \cite{locoresc}.  

We will use the ``cher \`a Cartan'' isomorphism to produce a na\"ive 
deformation and obstruction theory for twisted sheaves (without making 
use of the whole topos of sheaves on the gerbe).  First, we recall 
(without proof) what the cher \`a Cartan isomorphism is, in our 
context.

\begin{prop}\label{C:general cher} Let $B\to B_{0}$ be a morphism of 
rings in $X$ and $\ms X\to X$ a $\G_{m}$-gerbe.  Given a complex of $\ms 
X$-twisted $B$-modules $M$ and a complex of $\ms X$-twisted 
$B_{0}$-modules $J$, there is a natural 
isomorphism in the derived category
$$\rhom_{B}(M,J)\simto\rhom_{B_{0}}(M\ltensor_{B}B_{0},J).$$
\end{prop}

This is the derived 
adjointness of $\ltensor_{B}B_{0}$ and (derived) restriction of scalars 
to $B$.  (Note that this also applies when the gerbe $\ms X$ is 
trivial, so in particular in any topos.)  We will use this below in \ref{P:loc and constr} to deduce 
localization and constructibility properties for the deformation 
and obstruction theory of twisted sheaves.

\begin{cor}[cher \`a Cartan]\label{C:cher} Let $B\to B_{0}$ be a surjection of rings in $X$ and $\ms 
X\to X$ a $\G_{m}$-gerbe.  Given two $\ms X$-twisted $B_{0}$-modules 
$M_{0}$ and $J$, there is a natural isomorphism
$$\rhom_{B}(M_{0},J)\to\rhom_{B_{0}}(M_{0}\ltensor_{B}B_{0},J)$$
in the derived category of $\ms X$-twisted sheaves.
\end{cor}

Let $B\to B_{0}$ is a square-zero extension of rings in 
$X$ with kernel $I$.  Suppose $M_{0}$ and $J$ are twisted 
$B_{0}$-modules.  We wish to know when there exists an $I$-flat 
extension of $M_{0}$ by $J$.  Given any such extension, there is a 
naturally resulting morphism $I\tensor_{B_{0}}M_{0}\to J$, which is an 
isomorphism if and only if the extension is $I$-flat.  Fix a morphism 
$u:I\tensor_{B_{0}}M_{0}\to J$.  As in \S IV.3.1 of \cite{illusie}, we have the 
following proposition.

\begin{prop}\label{P:cher deformation theory}  There is an exact sequence
    $$0\to\ext^{1}_{B_{0}}(M_{0},J)\to\ext^{1}_{B}(M_{0},J)\to\hom_{B_{0}}(I\tensor_{B_{0}}M_{0},J)\xto{\partial}\ext^{2}_{B_{0}}(M_{0},J)$$
with the property that there exists an extension with associated 
morphism $u$ if and only if $\partial(u)=0$.  The space of all such 
extensions is a torsor under $\ext^{1}_{B_{0}}(M_{0},J)$.
\end{prop}
\begin{proof} The exact sequence is the sequence of low degree terms arising 
from \ref{C:cher} and the composition of functors spectral sequence.  
That the maps agree with the interpretation given is checked carefully 
in \cite[p.\ 252ff]{illusie}.  Note that Illusie's proof works in the derived 
category of twisted sheaves; it is not necessary to work in the 
category of all modules in the topos. (The Ext groups are different, but the 
functorialities are the same.)
\end{proof}

\begin{para}\label{Para:obstruction} We now have enough information to describe an obstruction theory for 
the problem of twisted sheaves on a scheme.  In this section, $f:X\to S$ 
will be a proper morphism of Noetherian excellent algebraic spaces and
$\ms X\to X$ will be a 
fixed $\G_{m}$-gerbe.  (An algebraic space is excellent if 
\emph{every\/} \'etale chart is excellent.  In 18.7.7 of \cite{ega4-4}, the reader 
will find an example of a non-excellent scheme with an excellent 
finite \'etale cover.  We thank Brian Conrad for pointing out this example.)  
We will develop the deformation-theoretic tools 
necessary to apply Artin's Existence Theorem.  Let $A_{0}$ be a 
reduced Noetherian ring.  We recall some terminology from Artin's 
paper \cite{artin}.
\end{para}
\begin{defn} A \emph{deformation situation\/} is a commutative 
diagram of Noetherian rings $A'\to A\to A_{0}$ such that 
\begin{enumerate}
    \item $A\to A_{0}$ and $A'\to A$ are infinitesimal extensions 
    (i.e., they have nilpotent kernels)

    \item $\ker(A'\to A)=M$ is a finite $A_{0}$-module.
\end{enumerate}
\end{defn}
In the classical study of versal deformations, one often takes $A_{0}$ 
to be a field and $A$, $A'$ to be local Artinian rings with residue 
field $k$.

Let $F$ be a stack on $S$.
\begin{defn} An \emph{obstruction theory\/} for $F$ consists of two 
parts.
\begin{enumerate}
    \item[(i)] For each infinitesimal extension $A\to A_{0}$ and 
    element $a\in F(A)$, a functor 
    $\ob_{a}:\operatorname{Mod}_{A_{0}}^{\operatorname{finite}}\to\operatorname{Mod}_{A_{0}}^{\operatorname{finite}}$

    \item[(ii)] For each deformation situation and $a\in F(A)$, there 
    is an element $o_{a}(A')\in\ob_{a}(M)$ which vanishes if and only 
    if there is an element $F(A')$ whose reduction to $A$ is 
    isomorphic to $a$.
\end{enumerate}
These data are subject to two further constraints:
\begin{enumerate}
    \item[\bf F:] Given a diagram 
    $$\xymatrix{B\ar[dr]\ar[dd]^{f} &  \\ & A_{0} \\
		A\ar[ur] & }$$
    of infinitesimal extensions, one has $\ob_{a}=\ob_{f(a)}$ as functors 
    $\finitemod_{A_{0}}\to\finitemod_{A_{0}}$.
    \item[\bf L:] For any diagram of deformation situations
    $$\xymatrix{A'\ar[r]\ar[dd]^{g} & A\ar[dd]\ar[dr] & \\ & & A_{0}\\
		B'\ar[r] & B\ar[ur] &}$$
    giving rise to an $A_{0}$-linear map of kernels $M_{A}\to M_{B}$, we get 
    for any $a\in F(A)$ an $A_{0}$-linear map 
    $\ob_{a}(M_{A})\to\ob_{a}(M_{B})$ taking $o_{a}(A')$ to $o_{a}(B')$.            
\end{enumerate}
We will call {\bf F} \emph{functoriality\/} and {\bf L} 
\emph{linearity\/} of the obstruction theory.
\end{defn}
Let $F$ be the stack which assigns to any Noetherian affine scheme $\spec A\to 
S$ the groupoid of $A$-flat families of coherent $\ms X$-twisted 
sheaves $\ms F$ on $X\tensor_{S}A$.
\begin{lem}\label{L:finite coho} If $\ms F$ and $\ms G$ are $A$-flat coherent $\ms 
X\tensor A$-twisted 
sheaves then $\ext^{i}(\ms F,\ms G)$ is a finite $A$-module.
\end{lem}
\begin{proof} This follows from the local to global spectral sequence 
for Ext and the finiteness of coherent cohomology for a proper 
morphism.  (Coherence of the sheaf Exts is a local computation in the 
big \'etale topology of $\ms X$, hence follows from the 
corresponding fact for locally Noetherian schemes.)
\end{proof}

\begin{prop} The following give an obstruction theory for $F$. 
    \begin{enumerate}
	\item Given an infinitesimal extension $A\to A_{0}$ and $\ms F\in 
	F(A)$, $$\ob_{\ms F}(M)=\ext^{2}_{X\tensor A_{0}}(\ms 
	F_{0},M\tensor_{A_{0}}\ms F_{0})=\ext^{2}_{X\tensor A}(\ms 
	F,M\tensor_{A}\ms F).$$
    
	\item Given a deformation situation $A'\to A\to A_{0}$ with 
	kernel $M$, $o_{\ms F}(A')=\partial(\id:M\tensor_{A}\ms F\to 
	M\tensor_{A}\ms F)$ in \ref{P:cher deformation theory}
    \end{enumerate}
\end{prop}
The equality of (1) above is a simple consequence of the cher \`a 
Cartan isomorphism and the fact that $\ms F\ltensor_{A}A_{0}=\ms 
F_{0}$ (by flatness).
\begin{proof} Using \ref{P:cher deformation theory}, it suffices to 
check {\bf F} and {\bf L} and prove that $\ext^{2}(\ms 
F_{0},M\tensor_{A_{0}}\ms F_{0})$ is a finite $A_{0}$-module.  {\bf F} follows from the 
description of the obstruction group in terms of $A_{0}$ and {\bf L} 
follows from the naturality of \ref{P:cher deformation theory}.  The 
finiteness is \ref{L:finite coho}.  
\end{proof}

We can use this formalism to prove several ``localization and constructibility'' 
results about the deformation theory of coherent twisted sheaves.  
These are the ``conditions (4.1)'' of Artin's famous \cite{artin}.  

\begin{prop}\label{P:loc and constr} Let 
$A_{0}$ be a reduced Noetherian ring, $A_{0}\to B_{0}$ a flat ring 
extension, $f:X\to\spec A_{0}$ a proper 
morphism, $\ms X\to X$ a $\G_{m}$-gerbe, $\ms F$ an $A_{0}$-flat 
family of coherent $\ms X$-twisted sheaves, and $M$ an 
$A_{0}$-module.  For any $i\geq 0$ the following hold.
    \begin{enumerate}
	\item $\ext^{i}_{X}(\ms F,M\tensor\ms 
	F)\tensor_{A_{0}}B_{0}\cong\ext^{i}_{X_{B_{0}}}(\ms 
	F_{B_{0}},M_{B_{0}}\tensor\ms F_{B_{0}})$.
    
	\item If $\mf m\subset A_{0}$ is a maximal ideal then 
	$$\ext^{i}_{X}(\ms F,M\tensor\ms 
	F)\tensor\widehat{A_{0}}\cong\invlim\ext^{i}_{X}(\ms F,M/\mf 
	m^{n}M\tensor\ms F),$$ the completion of $A_0$ being taken with respect 
	to $\mf m$.
    
	\item There is a dense open set of points (of finite type) 
	$p\in\spec A_{0}$ such that $$\ext^{i}(\ms F,M\tensor\ms 
	F)\tensor\kappa(p)\cong\ext^{i}_{X_{\kappa(p)}}(\ms 
	F_{\kappa(p)},M_{\kappa(p)}\tensor\ms F_{\kappa(p)}).$$
    \end{enumerate}
\end{prop}
\begin{proof}  The proof of 1 is immediate.  To prove 2, we work in 
\Caldararu\ form.  This makes it clear that one can easily understand 
formal twisted sheaves on the formal completion of a scheme along a 
closed subscheme.  We wish to prove that if $X\to\spec A$ is a proper 
scheme over a complete Noetherian local ring and $\ms F$ and $\ms G$ 
are coherent twisted sheaves on $X$ then $$\ext^{i}_{X}(\ms F,\ms 
G)=\invlim\ext^{i}_{X\tensor A/\mf m^{n}}(\ms F\tensor A/\mf m^{n},\ms 
G\tensor_{A}A/\mf m^{n}).$$  This works just as in \cite[4.5]{ega3-1}: one 
shows that the completion of the sheaf $\sext^{i}(\ms F,\ms G)$ along 
the closed fiber is naturally isomorphic to the sheaf 
$\sext^{i}(\widehat{\ms F},\widehat{\ms G})$ of extensions over the 
formal scheme.  The rest comes by taking the local-to-global Ext 
spectral sequence and using the finiteness of coherent cohomology to 
make an Artin-Rees argument.  The interested reader should consult 
\S III.4.5 of \cite{ega3-1} for further details.

The proof of 3 is slightly subtle.  Suppose first that $A_{0}$ is a 
finite type $\Z$-algebra.  Localizing, we may suppose that $A_{0}$ is 
a regular Noetherian ring of finite Krull dimension $d$ and that $X$ 
and $M$ are $A_{0}$-flat.  Thus, any 
$A_{0}$-module has homological dimension at most $d$.  Let $\ms 
C=\rshom(\ms F,M\ltensor\ms F)\in\D(X)$.  It is easy to see (using the 
bound on the homological dimension) that one has $$\ms C\ltensor 
\kappa(p)=\rshom(\ms F,M\ltensor\ms F\ltensor\kappa(p)).$$  
Furthermore, the bound on homological dimension implies that 
$\tau_{\leq i+d}\ms C\to\ms C$ remains a quasi-isomorphism (of 
bounded-below complexes) in 
degrees $\leq i$ upon any base change.  Thus, we see that we are 
concerned with the base change properties of $\R f_{\ast}\tau_{\leq 
i+d}\ms C$.  Localizing $A_{0}$, we may suppose that the (coherent) 
cohomology sheaves of $\tau_{\leq i+d}\ms C$ are all $A_{0}$-flat.  
By the standard cohomology and base change argument, we then see that the 
formation of $\R f_{\ast}\tau_{\leq i+d}\ms C$ is compatible with base 
change.  Passing to the open subscheme of $\spec A_{0}$ over which the 
cohomology sheaves of $\R f_{\ast}\tau_{\leq i+d}\ms C$ are flat 
yields the open subset we seek.

In the case of arbitrary reduced $A_{0}$, note that since $X$ is of 
finite type (hence of finite presentation over $A_{0}$, everything 
being Noetherian) and $\ms F$ is coherent, we can descend $X$,  
$\ms F$, and $M$ to a finite type $\Z$-subalgebra $B\subset A_{0}$.  
Localizing, we may take $B$ to be regular of finite Krull dimension.  
Applying the previous paragraph, upon shrinking $\spec B$ enough, we 
find a perfect complex $\ms P$ whose cohomology universally computes 
the Ext spaces in question.  Thus, the question reduces to the case 
already treated.
\end{proof}

\subsubsection{Determinants and equideterminantal deformations}\label{S:equi}

Given a perfect object of the derived category of twisted sheaves on 
a topos, one can use the construction of Mumford and Knudsen 
\cite{mumford-knudsen} to define its determinant.  Given a ringed topos 
$X$ and a $\G_{m}$-gerbe $\ms X\to X$, let $\D^{\tau}(\ms X)$ denote 
the derived category of $\ms X$-twisted sheaves.  There is a natural 
map $\D^{\tau}(\ms X)\to\D(\ms X)$, but it is not clear what 
properties this map has.  Of course, if $X$ is an algebraic space and 
once considers only quasi-coherent cohomologies, then the natural 
functor is an equivalence onto a direct summand triangulated category.)

\begin{defn} Let $X$ be a topos, $\ms X\to X$ a $\G_{m}$-gerbe, and 
$\ms F$ an $\ms X$-twisted sheaf on $X$ such that $\ms F$ is perfect 
as an object of $\D^{\tau}(\ms X)$.  The \emph{determinant\/} of $\ms F$, 
denoted $\det\ms F$, is the Knudsen-Mumford determinant of the 
complex $\ms F\in\D^{\tau}(\ms X)$.
\end{defn}
It is clear from the construction that $\det\ms F$ can be computed in 
either $\D(\ms X)$ or $\D^{\tau}(\ms X)$.  More generally, if $X$ is 
an algebraic space then it is clear that the restriction of the functor $\D^{\tau}(\ms 
X)\to\D(\ms X)$ to the sub-triangulated category $\D^{\tau}(\ms 
X)_{\text{parf}}$ of perfect complexes 
induces an equivalence with a triangulated direct summand of the triangulated 
category $\D(\ms X)_{\text{parf}}$.  This will not be of any use to us.

Our goal in this section is to study deformations of a twisted sheaf 
which fix its determinant.  Let $I\to A\to A_{0}$ be a small 
extension of Noetherian rings over $S$ and $X\to S$ a proper algebraic 
space of finite presentation.  We assume 
that $n$ is invertible in $A_{0}$ in what follows.  Fix a 
$\m_{n}$-gerbe $\ms X$ on $X_{A}$.  (By standard results in \'etale 
cohomology, $\ms X$ is in fact determined by the $\m_{n}$-gerbe 
structure on $\ms X_{A_{0}}$; this fact is relevant to the study 
of deformations of Azumaya algebras and their generalizations.  The 
reader is referred to \cite{pgl-bundles} for details.)

\begin{lem}\label{L:cher a cartan revisited} Let $(X,{\ms O})$ be a ringed topos and $A$, $B$, and $C$ 
complexes of ${\ms O}$-modules.  There is a natural isomorphism
$$\rshom(A\ltensor B,C)\simto\rshom(A,\rshom(B,C))$$
and a natural isomorphism
$$\rhom(A\ltensor B,C)\simto\rhom(A,\rshom(B,C)).$$
\end{lem}
\begin{proof}[Sketch of proof]  For further details see e.g.\ \cite{spaltenstein}.  
This is a close relative of \ref{C:general cher} and 
can proven similarly using the techniques of Neeman and 
Spaltenstein: replace $A$ and $B$ by $K$-flat resolutions $F_{A}$, 
$F_{B}$, and $C$ by 
a $K$-injective resolution $I_{C}$.  Then $\shom(F_{B},I_{C})$ is weakly 
$K$-injective, hence $\rshom(A,\rshom(B,C))$ is computed by 
$$\shom(F_{A},\shom(F_{B},F_{C})).$$  Using the hom-tensor adjunction on 
modules, this is naturally isomorphic to $$\shom(F_{A}\tensor 
F_{B},I_{C}),$$ which computes $\rshom(A\ltensor B,C)$ as usual.  The 
last formula follows upon taking derived global 
sections of the sheafified version.
\end{proof}

If $\ms F$ is \emph{perfect\/}, then there is a natural isomorphism $\ms F\simto\ms 
F^{\vee\vee}$.  Applying \ref{L:cher a cartan revisited} 
we see that to the identity in $\End(\ms F)$ corresponds some morphism 
$\hom_{\D}(\ms F\ltensor\ms F^{\vee},{\ms O})$.  This gives rise to a 
morphism $\rshom(\ms F,\ms F)\to{\ms O}$, called the \emph{trace 
morphism\/}, which we will denote $\Tr$.  In what follows, we will let $\ms A=\rshom(\ms F,\ms F)$; this is an example of a \emph{generalized Azumaya algebra\/}, the collection of which may be used to compactify moduli of Azumaya algebras.  This is discussed in \cite{pgl-bundles}.

\begin{defn} The homotopy fiber of $\Tr:\ms A\to{\ms O}$ in $\D(X)$ is the 
\emph{traceless part\/} of $\ms A$ and denoted $s\ms A$.
\end{defn}

\begin{lem}\label{L:dualize me} Under the natural isomorphisms $\ms A^{\vee}\simto\ms A$ 
and ${\ms O}^{\vee}\simto{\ms O}$, the trace is dual to the unit ${\ms O}\to\ms 
A$.
\end{lem}
\begin{proof} By functoriality, we can localize and assume that $\ms 
F$ is a strict perfect complex, where this is just a computation.
\end{proof}

\begin{lem}\label{L:trace splitting} If $\ms F$ is a perfect complex of 
${\ms O}$-modules, the composition $${\ms O}\to\rshom(\ms F,\ms 
F)\xto{\operatorname{Tr}}{\ms O}$$ is 
equal to multiplication by the rank of $\ms F$.
\end{lem}
\begin{proof} If $\ms F$ is a strict perfect complex, i.e., there is a 
quasi-isomorphism $\ms V\simto\ms F$ with $\ms V$ a finite complex of 
locally free modules, this comes down 
to checking that the adjunction is induced by the obvious maps.  As 
every perfect complex is locally quasi-isomorphic to such a complex, 
this will prove the general case by functoriality.
\end{proof}

\begin{defn} The \emph{reduced trace\/} of $\ms A$ is the map 
$\tau=\frac{1}{\rk\ms F}\Tr:\rsend(\ms F)\to{\ms O}$.
\end{defn}

\begin{prop}\label{P:split me} Let $f:A\to B$ be a map in the derived category $\D(\ms C)$ 
of an abelian category.  If $f$ has a section $g:B\to A$ then there is 
an isomorphism $\operatorname{holim}(g)\cong\operatorname{hocolim}(f)$.
\end{prop}
\begin{proof} In other words, the homotopy fiber of $g$ is isomorphic 
to the homotopy cofiber (``mapping cone'') of $f$.  This is a straightforward exercise which works in any 
triangulated category.
\end{proof}

\begin{cor}\label{C:duality of p and s} The third vertex $p\ms A$ of the unit ${\ms O}\to\ms A$ is isomorphic 
to the traceless part $s\ms A$.
\end{cor}
\begin{proof} This is an application of \ref{P:split me} to 
\ref{L:trace splitting} and \ref{L:dualize me}.
\end{proof}

The main result of this section is that the traceless part of $\rhom(\ms 
F,\ms F)$ governs the equideterminantal deformation theory of $\ms F$ 
(as long as $\det\ms F$ is unobstructed).  
(In fact, one can also see that the traceless part governs the 
deformation theory of the derived algebra $\rshom(\ms F,\ms F)$ in a 
precise manner.  This is discussed in \cite{pgl-bundles} and in 
\cite{mythesis} in great detail.)  
In the general case (when the rank is not invertible on the base), a more 
subtle analysis is called for.  It is not especially difficult, and 
may be found in \cite[\S8.4]{artin-dejong}, but we will not make use of it here.

For the sake of concreteness, we assume that the determinant is 
trivialized.  In general, one need only know that the determinant is 
unobstructed to make the following proposition valid as stated.

\begin{prop}\label{P:equideterminantal deformations} Let $\ms F$ be an $A_{0}$-flat 
$\ms X_{A_{0}}$-twisted coherent sheaf with torsion free fibers of rank $n$ and trivial 
determinant ${\ms O}_{X_{A_{0}}}\simto\det\ms F$.  Let $\ms A=\rshom(\ms 
F,\ms F)$.
\begin{enumerate}
    \item The obstruction to deforming $\ms F$ while preserving the 
determinant lies in the hypercohomology $\bH^{2}(I\ltensor s\ms 
A)=\ext^{2}(\ms F,I\tensor\ms F)_{0}$.

    \item The isomorphism classes of equideterminantal deformations of $\ms F$ are a 
    principal homogeneous space under the hypercohomology 
    $$\bH^{1}(I\ltensor s\ms A)=\ext^{1}(\ms F,I\tensor\ms F)_{0}.$$
    
    \item The 
    determinant-preserving infinitesimal automorphisms of a deformation 
    are equal to $\bH^{0}(I\ltensor s\ms A)=\hom(\ms F,I\tensor\ms 
    F)_{0}$.
\end{enumerate} 
\end{prop}
\begin{proof} According to Illusie's standard deformation theory of sheaves 
in topoi (which applies verbatim as $\ms F_{0}$ is coherent), we have only to show that 
the trace of the obstruction of $\ms F$ is the obstruction of 
$\det\ms F$.  As we will only use this in the case where $X/A_{0}$ is 
smooth and projective, we may assume that $\ms X$ has enough locally free twisted sheaves (perhaps this should be called ``twisted locally factorial'') and that every coherent twisted sheaf admits a finite locally free resolution.  (When $X$ is projective, the existence of locally free twisted sheaves follows from \ref{C:gabber 
de jong}.  Note that it is not known if being regular and separated ensures the existence of sufficiently many locally free twisted sheaves.)  The argument one can use to prove this is 
practically identical to the argument of Artamkin \cite{artamkin} and 
proceeds by induction on the homological dimension of $\ms F$.  If 
$\ms F$ is locally free, the statement is quite easy.  The inductive 
step works as follows: choose a surjection $0\to\ms K\to \ms V\to\ms F\to 0$ 
with $\ms V$ a locally free twisted sheaf whose deformation is 
unobstructed.  Then the obstruction to deforming $\det\ms F$ is the 
same as the obstruction to deforming $\det\ms K$.  Furthermore, $\ms 
K$ has smaller homological dimension, hence the obstruction of 
$\det\ms K$ is the trace of the obstruction of $\ms K$.  A simple 
argument shows that the trace of the obstruction of $\ms F$ equals the 
trace of the obstruction of $\ms K$.

The second statement works in a similar way and uses \ref{C:duality of 
p and s}.  The last statement is left to the reader.
\end{proof}

\subsubsection{Optimality}\label{S:optimality}

The following notion will appear from time to time throughout this 
paper, so we honor it with its own subsubsection.

\begin{defn} Given a $\m_{n}$-gerbe $\ms X\to X$, the \emph{index\/} 
of $\ms X$, denoted $\ind(\ms X)$ is the minimal rank of a locally free $\ms 
X$-twisted sheaf over the generic scheme of $X$.  The \emph{period\/} of $\ms X$, denoted $\per(\ms 
X)$, is the order of the 
image of $[\ms X]$ in $\Br(X)$.  
\end{defn}

\begin{defn}\label{D:optimal}
A $\m_n$-gerbe $\ms X$ is \emph{optimal\/} if the period is $n$.
\end{defn}

When $X$ is regular, one immediately seems upon taking determinants that  
$\per(\ms X)|\ind(\ms X)$: any $\ms X$-twisted sheaf of rank $m$ over the generic scheme of $X$ extends to a coherent $\ms X$-twisted sheaf of rank $m$, and forming $\det\ms F$ yields an $m$-fold $\ms X$-twisted \emph{invertible\/} sheaf, which yields $m[\ms X]=0\in\H^2(X,\G_m)$.  (For details the reader is referred to \cite{period-index-paper}.)  
Using Galois cohomology at generic points, one also sees that $\ind(\ms X)|\per(\ms X)^{m}$ for 
some $m$ (see \cite{lang-tate}).  When $X$ is a surface (over an algebraically closed field), a theorem of 
de Jong (which is re-proven in \cite{period-index-paper} using the 
theory developed here and the theorem of Graber-Harris-Starr-de Jong 
\cite{dejong-starr}) shows that $\per(\ms 
X)=\ind(\ms X)$.  Thus, \emph{on a surface\/}, the index of a 
$\m_{n}$-gerbe divides $n$.  It is easy to show that the rank of any locally 
free $\ms X$-twisted sheaf is divisible by $\ind(\ms X)$.

Studying moduli of twisted sheaves of rank $n$ on an optimal 
$\m_{n}$-gerbe $\ms X\to X$ is a non-commutative analogue of the Picard scheme, in 
the sense that such sheaves are essentially rank $1$ right modules 
over an Azumaya algebra on $X$.  Thus, the stability condition which 
we define in section \ref{S:semistability} becomes vacuous, making certain 
proofs technically easier.  While we have not yet written out the 
general proofs, we strongly feel that this is a non-essential 
distinction for the theorems of \ref{S:asymptotic properties} to 
hold.  This deficit will be addressed in future work.

\subsubsection{Purity of sheaves on Artin stacks}\label{S:purity}
In this section, we will study \emph{purity\/} of 
twisted sheaves as a precursor to \ref{S:semistability}, where we 
will study various stability conditions on twisted sheaves.  The ultimate goal is to produce a tractable algebraic stack parametrizing a well-behaved collection of twisted sheaves.  We develop most of this section in much greater generality for coherent sheaves on algebraic stacks.

\begin{para} \emph{Support of twisted sheaves\/}.  Twisted sheaves 
may be viewed both as objects on $X$ and as 
objects on a $\m_{n}$-gerbe $\ms X$ over $X$.  This 
leads to two natural definitions of support for a twisted sheaf, 
which coincide.  (In the sequel, when the gerbe 
is understood we will often 
refer to ``twisted sheaves on $X$'' for the sake of notational 
simplicity.)
\end{para}

\begin{defn} Given a $\m_{n}$-gerbe $\ms X\to X$ and an $\ms 
X$-twisted coherent sheaf $\ms F$, the \emph{support\/} of $\ms F$ is the 
closed substack of $\ms X$ defined by the kernel of the map $\ms 
O_{\ms X}\to \send_{\ms X}(\ms F)$, which is a quasi-coherent sheaf of 
ideals.  The \emph{schematic support\/} 
of $\ms F$ is the scheme-theoretic image in $X$ of the support of $\ms 
F$.
\end{defn}
Since $\send_{\ms X}(\ms F)$ is the pullback to $\ms X$ of a coherent 
${\ms O}_{X}$-algebra, it is immediate that the support of $\ms F$ is the 
preimage of the schematic support. 
In particular, a twisted sheaf $\ms F$ with 
schematic support $Y\subset X$ is naturally a $\ms 
X\times_{X}Y$-twisted sheaf with full schematic support (on $Y$).  
Thus, considering the support of a sheaf does not nullify its ``twistedness.''

\begin{para} \emph{Associated points on Artin stacks\/}.  
    We can now define a torsion 
filtration on a twisted sheaf.  To do this properly, we will briefly develop the 
theory of associated points and torsion subsheaves on an arbitrary 
Noetherian algebraic stack.  (When trying to generalize these results 
to the non-Noetherian case, certain equivalences will fail, making the 
theory developed here only one possibility.)  Throughout, we 
systematically work with the underlying topological space $|\ms X|$ of 
a Noetherian algebraic stack.  The support of a 
sheaf will be taken to mean simply the underlying set of points of 
$|\ms X|$, or the reduced closed substack structure on that set 
when it is closed (e.g., if $\ms F$ is coherent).  We will not require 
(as is typical) that $\supp(\ms F)$ is the closure of the set of 
points where $\ms F$ is supported.
\end{para}
Let $\ms F$ be a quasi-coherent sheaf on $\ms X$.
\begin{defn}\label{D:ass pt} A point $p\in|\ms X|$ is an \emph{associated point\/} 
of $\ms F$ if there is a quasi-coherent subsheaf $\ms G$ such that $p\in\supp(\ms 
G)\subset\widebar{\{p\}}$.  The set of associated points of $\ms F$ 
will be written $\ass(\ms F)$.
\end{defn}
If $\ms F$ is \emph{coherent\/}, this is the same as requiring that 
$\supp(\ms G)=\widebar{\{p\}}$.  In general, this is not the case, as 
supports need not be closed for quasi-coherent sheaves.

\begin{rem} When $\ms X$ is a Noetherian scheme, this is the same as the usual 
notion (essentially because one can extend quasi-coherent subsheaves off 
of generic points).  More generally, if $\ms X$ is a Noetherian DM stack, one can 
say that a geometric point $p\to\ms X$ is associated to $\ms F$ if $p$ 
is an associated point for the stalk of $\ms F$ at $p$ (as a module 
over ${\ms O}_{p,X}^{\text{sh}}$).  By an argument similar to \ref{P:pull 
back ass} below, a point of $|\ms X|$ is associated iff some (and hence any) geometric 
point lying over it is associated, so this also yields the same notion as 
\ref{D:ass pt}.
\end{rem}

\begin{prop}\label{P:pull back ass} Let $f:X\to\ms X$ be a flat surjection, with $X$ a
Noetherian scheme. If $\ms F$ is a quasi-coherent sheaf on $\ms X$, 
then $\ass(\ms F)=f(\ass(\ms F|_{X}))$.
\end{prop}
\begin{proof} Write $\ms F'=\ms F|_{X'}$.  Given a point $p\in\ass(\ms F)$, it is easy to see that 
a generic point of $f^{-1}(\widebar{\{p\}})$ will be in $\ass(\ms 
F')$.  Conversely, let $q\in\ass(\ms F')$ and let 
$Y=f^{-1}\widebar{\{f(q)\}}$ as a reduced closed subscheme of $X$.  
Let $\ms G\subset\ms F'$ be the maximal quasi-coherent subsheaf supported on $Y$.  
(It is not true that $\supp\ms G=Y$, but we at 
least know that $q\in\supp\ms G$.)  We claim that $\ms G$ descends to 
a subsheaf of $\ms F$ with support containing $f(q)$ and contained in 
$\widebar{\{f(q)\}}$.  To see this, it is enough to show that the two 
pullbacks of $\ms G$ to $X\times_{\ms X}X$ are equal as subsheaves.  
In fact, by colimit considerations, we may assume that $\ms F$ is 
coherent.  We are reduced to the following situation: given a flat surjection 
$g:Z\to W$ of Noetherian algebraic stacks with $W$ an affine scheme, a closed subspace $Y\subset W$, and a coherent sheaf $\ms 
F$ on $W$, let $\ms G_{Y}\subset\ms F$ denote the maximal subsheaf $\ms G$ with 
$\supp(\ms G)=Y$.  It suffices to show that $g^{\ast}(\ms G_{Y})$ is the maximal subsheaf of $\ms F$ with support on $f^-1(Y)$.  To 
prove this, let $\ms I$ be the ideal cutting out the reduced 
structure on $Y$.  By flatness, $\ms J=g^{*}\ms I$ is a sheaf of ideals 
cutting out a substack of $Z$ supported on $g^{-1}(Y)$.  To say 
that $\ms G_{Y}$ is maximal is the same as saying that the sheaf $\shom({\ms O}/\ms 
I^{n},\ms F/\ms G_{Y})$ vanishes for all $n>0$.  By flat pullback, we conclude that 
$\shom_{{\ms O}_{Z}}({\ms O}_{Z}/\ms J^{n},\ms F_{Z}/g^{\ast}\ms G_{Y})=0$, 
whence $g^{\ast}\ms G_{Y}$ is maximal.
\end{proof}

\begin{cor}\label{C:ass flat} If $f:\ms X'\to\ms X$ is a flat 
surjection of Noetherian algebraic stacks and $\ms F$ is a 
quasi-coherent ${\ms O}_{\ms X}$-module, then $\ass(\ms F)=f(\ass(\ms 
F|_{\ms X'}))$.
\end{cor}
\begin{proof} Choosing a smooth cover of $\ms X'$ reduces this to 
\ref{P:pull back ass}.
\end{proof}

\begin{cor}\label{C:ass finite} If $\ms F$ is a coherent sheaf on $\ms X$ then $\ass(\ms F)$ is finite.
\end{cor}
\begin{proof} The stack $\ms X$ has a smooth cover by a Noetherian 
scheme $X'$.  By \ref{P:pull back ass}, we are reduced to the case of 
a scheme, where this is a classical result \cite[6.5]{matsumura}.
\end{proof}
Points of a stack are subject to the relations of specialization and 
generization in the usual way. This gives $\operatorname{Ass}(\ms F)$ 
the structure of partially ordered set.  By \ref{C:ass finite}  
there are well-defined \emph{minimal\/} elements of 
$\operatorname{Ass}(\ms F)$.

It is easy to check that $\widebar{\operatorname{Ass}(\ms 
F)}=\widebar{\operatorname{Supp}(\ms F)}$ and that the minimal points of 
$\operatorname{Ass}(\ms F)$ coincide with the minimal points of 
$\operatorname{Supp}(\ms F)$.

\begin{lem}\label{L:ass lemma} Suppose $\ms F,\ms G,\ms H$ are three coherent sheaves 
on $\ms X$ fitting into an exact sequence $0\to \ms F\to\ms G\to\ms H\to 0$.
\begin{enumerate}
    \item $\operatorname{Ass}(\ms F)\subset\operatorname{Ass}(\ms G)\subset\operatorname{Ass}(\ms 
    F)\cup\operatorname{Ass}(\ms H)$.  If the sequence is split exact, 
    the second inclusion is a bijection.

    \item  The minimal points of $\operatorname{Ass}(\ms H)$ are 
    contained in $\operatorname{Ass}(\ms G)$

    \item If $\ms G\neq 0$, then $\operatorname{Ass}(\ms 
    G)\neq\emptyset$.
\end{enumerate}
\end{lem}
\begin{proof} This is precisely analogous to the classical proof 
\cite[\S6]{matsumura}.
\end{proof}

\begin{defn} A \emph{torsion subsheaf\/} of $\ms F$ is a subsheaf 
$\ms F'\subset\ms F$ with the property that none of the minimal 
points of $\operatorname{Ass}(\ms F)$ are contained in 
$\operatorname{Ass}(\ms F')$.
\end{defn}
Note that any minimal point of $\operatorname{Ass}(\ms F)$ which is 
also associated to a subsheaf $\ms F'$ will be minimal in 
$\operatorname{Ass}(\ms F')$.
\begin{lem}\label{L:the torsion subsheaf} The sum of any two torsion subsheaves of $\ms F$ is a 
torsion subsheaf.  There is a unique maximal coherent torsion subsheaf of 
$\ms F$.
\end{lem}
\begin{proof} Suppose $\ms F'$ and $\ms F''$ are torsion subsheaves 
of $\ms F$.  By \ref{L:ass lemma}, the minimal points of $\operatorname{Ass}(\ms F'+\ms 
F'')$ are contained in $\operatorname{Ass}(\ms 
F')\cup\operatorname{Ass}(\ms F'')$.  This proves the first 
statement.  The second follows by taking the sum of all torsion 
subsheaves of $\ms F$ (which is allowable because they form a set).
\end{proof}
The maximal torsion subsheaf of $\ms F$ will be called \emph{the\/} 
torsion subsheaf of $\ms F$ and denoted $T(\ms F)$.
\begin{lem}\label{L:T(F) eats the ass} Any non-minimal point of $\operatorname{Ass}(\ms F)$ is 
contained in $\operatorname{Ass}(T(\ms F))$.
\end{lem}
\begin{proof} Immediate from the definition!
\end{proof}
\begin{remark} When $\ms X$ is a gerbe bound by a diagonalizable group 
scheme, the decomposition \ref{P:q-coh splits} respects torsion 
subsheaves, so we see that we have also developed a good theory of 
torsion subsheaves for twisted sheaves.
\end{remark}

\begin{defn} A coherent sheaf $\ms F$ is \emph{pure\/} if $T(\ms F)=0$.
\end{defn}

\begin{remark} By \ref{L:T(F) eats the ass}, we see that $\ms F$ is 
pure if and only if $\operatorname{Ass}(\ms F)$ consists solely of minimal points, 
i.e., the partial ordering on $\operatorname{Ass}(\ms F)$ is trivial.
\end{remark}

\begin{lem}\label{L:purification no torsion} Given any coherent sheaf $\ms F$ on $\ms X$, the sheaf 
$\ms F/T(\ms F)$ is pure.
\end{lem}
\begin{proof} There is an exact sequence $0\to T(\ms F)\to\ms G\to 
T(\ms F/T(\ms F))\to 0$.  It follows from \ref{L:ass lemma} and the 
definition of the torsion subsheaf that $\ms 
G\subset T(\ms F)$.  Thus, $T(\ms F/T(\ms F))=0$.
\end{proof}

\begin{lem}\label{L:pure smooth ass} If $X\to\ms X$ is a smooth cover, then $\ms F$ is pure if 
and only if $\ms F|_{X}$ is pure.
\end{lem}
\begin{proof} As in \ref{P:pull back ass}, it suffices to show that 
if $Z\to W$ is a smooth map of schemes then the pullback of a torsion 
free sheaf is torsion free.  As this is a local property on the source 
and target and is obviously true for arbitrary 
quasi-finite flat morphisms (hence for \'etale morphisms), we see that 
it suffices to prove that the pullback of a torsion free sheaf on $W$ 
to $\A^{n}_{W}$ is torsion free.  Again by \ref{P:pull back ass}, we 
see that any torsion subsheaf on $\A^{n}_{W}$ must have all 
associated points lying over minimal (generic) points of $W$.  Thus, 
we may assume $W$ is the spectrum of an Artinian local ring $R$ and we 
wish to show that the pullback of any finite $R$-module to $\A^{n}_{R}$ 
cannot have torsion.  Taking a composition series, we may assume 
that $R$ is a field.  The result follows from the fact that 
$\A^{1}_{K}$ is Cohen-Macaulay.
\end{proof}

Let $\pi:\ms X'\to\ms X$ be a flat surjection 
of Noetherian algebraic stacks representable by an open immersion into 
an integral ring extension and $\ms F$ a coherent sheaf on $\ms X$.
\begin{prop} $\ms F$ is pure if and only if $\ms F|_{\ms X'}$ is pure.
\end{prop}
\begin{proof} By the going-up lemma \cite[9.4]{matsumura} and 
flatness (which implies the going-down lemma \cite[9.5]{matsumura}), 
a morphism such as $\pi$ has the property that $\pi(p)$ is minimal if 
and only if $p$ is minimal.  The result now follows from \ref{P:pull back ass}.
\end{proof}
\begin{cor} If $\ms X$ is over a field $k$, then the 
purity of a coherent sheaf is invariant under finite extensions of 
$k$.  If $\ms X$ is finite type over $k$ then purity is geometric.
\end{cor}
The finite type hypothesis in the second statement serves only to 
ensure that $\ms X\tensor K$ is Noetherian for any extension $K\supset 
k$, so that our theory applies.

\begin{remark} When $X$ is an integral universally catenary scheme of 
finite Krull dimension (for example, a 
projective variety) and $\ms F$ is an $\ms X$-twisted sheaf with 
support of dimension $d$, we can filter $T(\ms F)$ by the dimension of 
support: let $T_{e}(\ms F)$ be the maximal subsheaf of $\ms F$ whose 
support is of dimension at most $e$.  Then $T(\ms F)=T_{d-1}(\ms 
F)\supset T_{d-2}(\ms F)\supset\cdots\supset T_{0}(\ms F)$.  This 
filtration can be useful when considering various notions of 
semistability, as in \S 1.6 of \cite{h-l}; it will not come up in the sequel.
\end{remark}

\begin{para} We will now show that the property of being pure is open 
in flat families of coherent sheaves on a proper algebraic stack.  
\end{para}
\begin{prop} Let $\pi:\ms X\to S$ be a proper morphism of finite 
presentation from an algebraic stack to an algebraic space.  Suppose $\ms F$ is an $S$-flat family of coherent sheaves.  
The locus of points $s\in S$ such that $\ms F_{s}$ is pure is open. 
\end{prop}
\begin{proof} We may reduce to the case where $S$ is affine, Noetherian, and 
even excellent (in fact, affine of finite type over 
$\Z$).  Indeed, we may present the stack $\ms X$ as a groupoid 
$X_{1}\to X_{0}\times X_{0}$ of 
finite presentation between two schemes of finite presentation 
over $S$.  Thus, we may descend $\ms X$ to a Noetherian base (using the results of \S 8 of \cite{ega4-3}).  Having 
done this, note that a coherent sheaf on $\ms X$ is given by a 
coherent sheaf on $X_{0}$ with an action of the groupoid, i.e., an 
isomorphism of the pullbacks to $X_{1}$ which is compatible with the 
groupoid structure.  By Grothendieck's theory of limits, we can 
descend these data to a finite level.

Consider the set $\Xi$ of points $x\in|\ms X|$ with the property ``$x$ is 
contained in the support of $T(\ms F_{\pi(x)})$.''  It suffices to show that  
$\pi(\Xi)$ is constructible and that $\pi(\Xi)$ is closed under 
specialization when $\ms X$ is proper over $S$ (and $S$ is Noetherian).

The second statement is immediate: it suffices to check this when $S$ 
is the spectrum of a discrete valuation ring.  By flatness, the 
minimal points of $\operatorname{Ass}(\ms G)$ all lie in the generic 
fiber for any coherent subsheaf $\ms G\subset\ms F$.  Thus, the 
torsion subsheaf $T(\ms F)$ of the total family $\ms F$ is non-zero if 
and only if the torsion subsheaf of the generic fiber is non-zero.  On 
the other hand, since $S$ is Dedekind, $T(\ms F)$ is also a flat 
family, and in particular has constant fiber dimension, like $\ms F$.  
Finally, the cokernel $\ms F/T(\ms F)$ is pure, hence $S$-flat, by 
\ref{L:purification no torsion}.  These facts combine to yield the 
statement about specialization.  (Properness is not necessary for 
this, as long as we assume that the specialization on the base is 
contained in the image of $\pi$.)

The first statement (that $\pi(\Xi)$) is constructible) is more 
subtle.  It is easily reduced to showing that if $S$ is an integral and 
Noetherian affine scheme and the generic fiber of $\Xi$ is non-empty, then $\Xi$ is 
non-empty over an open subscheme of $S$.  By \ref{P:pull back ass}, we 
may assume that $\ms X$ is in fact a scheme and that $\pi$ is 
surjective.  The argument for schemes is classical, and is left to 
the reader.  (It is also written out in full detail as Proposition 
4.1.2.21 of \cite{mythesis}.)
\end{proof}
 
As a consequence of the Proposition, when $\ms X$ is a $\G_{m}$-gerbe, there is an open substack of 
$\ms T_{\ms X/S}$ (see \ref{N:curly T}) representing families of pure twisted sheaves.  Note 
that since the support of a flat family over a dvr is itself flat 
over the dvr, the dimension of the 
fibers of a flat family of coherent $\ms X$-twisted sheaves over a locally Noetherian base scheme is 
locally constant.
\begin{cor} Let $X\to S$ be a proper flat morphism of finite presentation 
with geometrically integral fibers.  There is an open substack of $\ms 
T_{\ms X/S}$ consisting of 
families of torsion free sheaves, i.e., pure sheaves of maximal 
dimension.
\end{cor}
\begin{defn} If $X\to S$ is a proper flat morphism of finite 
presentation and $\ms X\to X$ is a $\m_{n}$-gerbe, the open substack 
parametrizing families with torsion free fibers is denoted $\Tw_{\ms 
X/S}$.
\end{defn}

\begin{para} Suppose $X$ is a smooth projective variety over a field 
$k$ and $\ms X\to X$ is a $\m_{n}$-gerbe with $n\in k^{\times}$.  
Let $\Tw_{\ms X/k}(n)$ denote the stack parametrizing torsion free 
twisted sheaves of rank $n$.  Since $\ms X$ is smooth, any $S$-flat 
family of twisted sheaves $\ms F$ on $\ms X\times S$ has finite 
homological dimension everywhere.  In other words, $\ms F$ is 
\emph{perfect\/} as an object of the derived category.  Using the 
constructions of section \ref{S:equi}, we can thus define the 
determinant of $\ms F$,  which will be the pullback to $\ms 
X\times S$ of an invertible sheaf 
on $X\times S$ (as $\ms F$ has rank $n$).  This yields a morphism of algebraic stacks
$$\det:\Tw_{\ms X/k}(n)\to\ms Pic_{X/k}.$$
Given an invertible sheaf $\ms L$ on $X$, one can form the fiber of 
$\det$ over the resulting $k$-point of $\ms Pic_{X/k}$.
\end{para}
\begin{defn} $\Tw_{\ms X/k}(n,\ms L):=\Tw_{\ms X/k}\times_{\ms 
Pic_{X/k}}k$
\end{defn}
Chasing through the definition of the natural 1-fiber product of 
stacks shows that the objects of $\Tw_{\ms X/k}(n,\ms L)$ are pairs 
$(\ms F,\phi)$ consisting of a torsion free twisted sheaf $\ms F$ of 
rank $n$ and a chosen isomorphism $\det\ms F\simto\ms L$.  
The deformation theory for $\Tw(n,\ms L)$ was 
developed above in section \ref{S:equi}; as we show in 
\cite{pgl-bundles}, it is this deformation theory which governs the 
stack of Azumaya algebras ($\PGL_{n}$-bundles) and its 
compactification by ``generalized Azumaya algebras.''

\subsubsection{Riemann-Roch, Hilbert polynomials, and Quot spaces}\label{S:hilbert polynomials and quot spaces}
In this section, we develop a notion of Hilbert polynomial for twisted 
sheaves which we will ultimately use to define semistability.  In 
later sections, we will show that on a surface (and more generally on 
a variety which carries a twisted sheaf with sufficiently many vanishing Chern 
classes), our notion agrees with Simpson's notion \cite[\S3]{simpson} and thus yields a 
GIT quotient corepresenting the stack of semistable twisted sheaves.  For higher dimensional ambient varieties, 
it will still be possible to show that the stack of stable twisted 
sheaves is a gerbe over an algebraic space, but dealing with 
properly semistable points is difficult in the absence of a GIT description of the moduli problem.

\begin{para} In order to define our semistability condition, and for 
future reference, we briefly recall the basic facts about rational 
Chow rings of DM stacks over a field.  Vistoli \cite{vistoli} and 
Gillet \cite{gillet} have defined Chow 
theories which only work rationally but which are formally identical 
to the usual Chow theory: in Vistoli's approach, one takes the Chow groups to be generated by 
integral closed substacks modulo rational equivalence (suitably 
defined).  There is a refined theory due to Edidin and 
Graham \cite{edidin-graham} which 
applies to quotient stacks to yield an integral Chow theory which 
agrees with Vistoli's theory when tensored with $\Q$.  A further 
refinement of the integral theory for algebraic stacks stratified by 
quotient stacks was developed by Kresch in his 
thesis \cite{kresch-cycle}.  We will denote 
the rational Chow groups (which are the same in all of these theories) by $A_{\Q}$, and we will write $A^{n}_{\Q}$ for 
the group generated by cycles of codimension $n$.  When the 
underlying stack is smooth, the graded group $\oplus A^{n}_{\Q}$ has a 
commutative ring structure.  As usual, there is a theory of Chern classes 
and a splitting principle.  The theory admits proper pushforwards, 
flat pullbacks, and Gysin maps \cite{vistoli}.  It is useful to note 
that since the splitting principle uses only the construction of 
projective bundles, one need never leave the category of smooth tame 
DM stacks with (quasi-)projective coarse moduli spaces, if one so 
desires.

Given a proper DM stack $\ms X$ with moduli space $X$, 
one can show that the proper pushforward $A_{\Q}(\ms X)\to A_{\Q}(X)$ is an 
isomorphism which respects the ring structure when both are smooth (see \cite{vistoli}).  
In particular, when $\ms X$ is of dimension $n$, there is a rational 
degree function $\deg:A^{n}(\ms X)_{\Q}\to\Q$.  Given 
any element $\alpha$ of the graded group $A^{\ast}_{\Q}(\ms X)$, we will 
let $\alpha_{n}$ denote the part in degree $n$.  Given a class 
$\beta\in A^{\ast}(\ms X)_{\Q}$, we will let $\deg\beta$ denote the 
degree of $\beta_{n}$.


Let $\ms X$ be a smooth proper DM stack of dimension $n$ over a field $k$ with projective 
moduli space $X$.  Recall that $K^{0}(\ms X)$ is the Grothendieck group of vector 
bundles on $\ms X$, while $K_{0}(\ms X)$ is the Grothendieck group of 
coherent sheaves.  When every coherent sheaf on $\ms X$ admits a 
finite resolution by locally free sheaves, it is easy to see that 
$K^{0}\cong K_{0}$.  In general, $K^{0}$ is a ring and $K_{0}$ is a 
$K^{0}$-module (via tensor product).  
One of the basic problems for arbitrary DM 
stacks is the fact that $K^{0}$ and $K_{0}$ are 
not isomorphic even on smooth DM stacks.  For smooth quotient stacks, they 
are naturally the same, which makes it easier to prove theorems.  (This 
is yet another place where \ref{C:gabber de jong} and its corollaries 
have a large impact.)  

Let $\alpha\in K^{0}(\ms X)$.  We will write $\td_{\ms X}$ for the 
Todd class of the tangent sheaf $T_{\ms X/k}$ of $\ms X$.\end{para}
\begin{defn} The \emph{geometric Euler characteristic\/} of $\alpha$ 
is
$$\chi^{g}(\alpha):=[\ms I(\ms X):\ms X]\deg(\chern(\alpha)\cdot\td_{\ms X}).$$
When $X$ is projective with chosen polarization ${\ms O}(1)$, the 
\emph{geometric Hilbert polynomial\/} of $\alpha$ is the function
$$n\mapsto P^{g}_{\alpha}(n)=\chi^{g}(\alpha(n)),$$
where $\alpha(n):=\alpha\tensor{\ms O}(n)$.
\end{defn}
To verify that $P^{g}_{\alpha}$ is a polynomial, it suffices to prove 
it when $\alpha=[\ms E]$, $\ms E$ a locally free sheaf on $\ms X$.  
This then follows by a simple splitting principle calculation left to 
the reader.

\begin{remark} The geometric Euler characteristic and Hilbert function are clearly 
additive functions on the category of (perfect) coherent 
sheaves.  When $\ms F$ is the pullback to $\ms X$ of a coherent sheaf 
on $X$, they agree with the usual Euler characteristic and Hilbert 
function by the Grothendieck-Hirzebruch-Riemann-Roch theorem.  
However, for sheaves which are not pullbacks, they do not agree with 
the usual cohomologically defined functions.  For a trivial example, 
consider the case of a gerbe $\ms X$ over an algebraic curve $X$.  In this case, there is an invertible sheaf $\ms 
L$ on $\ms X$ which is $\ms X$-twisted whose $n$th tensor power $\ms 
L^{\tensor n}$ is the pullback of an invertible sheaf $\ms M$ on $X$.  
The geometric Euler characteristic of $\ms L$ is easily seen to be 
$$\chi^{g}(\ms L)=\deg c_{1}(\ms L)+\chi({\ms O}_{X})=\frac{1}{n}\deg 
c_{1}(\ms M)+\chi({\ms O}_{X}).$$
Thus, one can easily produce gerbes $\ms X$ and 
$\ms X$-twisted sheaves with non-zero $\chi^{g}$.  On the other hand, 
if we use coherent cohomology to compute the cohomological Euler 
characteristic, we find $\chi(\ms L)=0$ when $\ms L$ has non-trivial 
stabilizer action.  (For an even more trivial example, let $\ms X$ be 
a gerbe over a point!)  There are ways of rectifying this difference, due to To\"en \cite{toen}, by instead working with the Chow theory of the inertia stack paying more careful attention to the representations of inertia on fibers of vector bundles.  It is interesting to note that in many cases the ``correct'' Riemann-Roch formula yields $0$, whereas the seemingly blunt instrument wielded here produces non-zero answers, thus somehow capturing geometric information about sheaves which is not cohomological and which is not visible on non-stacky varieties.
\end{remark}

Recent results of Vistoli-Kresch \cite{vistoli-kresch}, Edidin-Hassett-Kresch-Vistoli 
\cite{vistoli-kresch-etc}, and Gabber/de Jong \cite{dejong-gabber} (stated by Gabber and proven by Gabber 
and independently by de Jong) show that any separated smooth generically tame DM stack over a field with quasi-projective 
moduli space is a quotient stack, 
and that such a stack has the ``resolution property'': any coherent 
sheaf is a quotient of a locally free sheaf.  In these cases, the 
natural map $K^{0}\to K_{0}$ is thus an isomorphism.  We will denote 
it simply by $K(\ms X)$.

\begin{prop}\label{P:gen rr} Let $f:\ms X\to\ms Y$ be a projective 
l.c.i.\ morphism of DM stacks of finite type over a field 
such that $\ms Y$ admits a finite flat cover $\pi:Y\to\ms Y$.  
Given any class $\alpha\in K^{0}(\ms X)$, there is a natural equality
$$\chern(f_{\ast}\alpha)=f_{\ast}(\chern(\alpha)\cdot\td_{f})$$
in $A(\ms Y)_{\Q}$.
\end{prop}
\begin{proof} Form the Cartesian diagram
    $$\xymatrix{X\ar[r]^{\tilde\pi}\ar[d]^{\tilde f} & \ms X\ar[d]^{f}\\
    Y\ar[r]^{\pi} & \ms Y.}$$
The morphism $\tilde f$ is a projective l.c.i.\ morphism of schemes; 
hence the formula holds for $\tilde f$.  Furthermore, we have that 
$\pi_{\ast}\pi^{\ast}$ and $\tilde\pi_{\ast}\tilde\pi^{\ast}$ are both 
identified with multiplication by $\deg\pi$ (which follows immediately from the definitions of the pushforward and pullback functors just as in the classical situation); thus, to show an equality 
in rational Chow groups, it suffices to show equality after pulling 
back by $\pi$.  Furthermore, the 
formation of $\chern(\alpha)$ commutes with flat pullback and 
$\tilde\pi^{\ast}\td_{f}=\td_{\tilde f}$.  Finally, we know that 
flat pullback commutes with proper pushforward (3.9 of \cite{vistoli}).  
Combining these statements yields 
$$\pi^{\ast}f_{\ast}(\chern(\alpha)\cdot\td_{f})=\tilde 
f_{\ast}(\chern(\tilde\pi^{\ast}\alpha)\cdot\td_{\tilde f})=\chern(\tilde 
f_{\ast}\tilde\pi^{\ast}\alpha)=\pi^{\ast}\chern(f_{\ast}\alpha).$$
\end{proof}

\begin{cor}\label{P:rr} Let $f:\ms X\to \ms Y$ be a projective morphism of 
smooth pseudo-projective DM stacks.  Then for all $\alpha\in K(\ms 
X)$,
$$\chern(f_{\ast}\alpha)\cdot\td_{\ms Y}=f_{\ast}(\chern(\alpha)\cdot\td_{\ms 
X})$$
in $A(\ms Y)_{\Q}$.
\end{cor}
\begin{proof} Any such morphism must be l.c.i., so we can apply 
\ref{P:gen rr}, once we note that any smooth generically tame DM 
stack over a field admits a finite flat cover by a smooth scheme 
(apply the main result of \cite{dejong-gabber} to Theorem 2.2 of 
\cite{vistoli-kresch} and use this as input into Theorem 2.1 of 
\cite{vistoli-kresch}).
\end{proof}

\begin{cor}\label{C:hrr closed immersion} Let $\iota:\ms X\inj\ms Y$ be 
a closed immersion of smooth pseudo-projective DM stacks and $\ms F$ a coherent sheaf on $\ms X$.  
Then $\chi^{g}(\ms X,\ms F)=\chi^{g}(\ms Y,\iota_{\ast}\ms F)$.
\end{cor}

\begin{rem} When $\ms X$ is a $\m_{n}$-gerbe over a smooth projective 
variety $X$ and there is a locally free $\ms X$-twisted sheaf with 
sufficiently many vanishing Chern classes (e.g., $X$ is a surface), 
then the formula in the last sentence of the proof of \ref{L:semistab 
same as simpson} below gives a much more concrete proof of \ref{C:hrr 
closed immersion}.
\end{rem}

We fix a smooth pseudo-projective DM stack $\ms X$ with moduli space 
$\pi:\ms X\to X$ in what follows.  For the moment, we assume that the base is a field; we will see in a moment that the geometric Hilbert function is constant in flat families and invariant under extension of base field.  We also assume that the moduli space $X$ is quasi-projective over $k$, with fixed very ample $\O(1)$.  We prove our results under the following hypothesis.

\begin{hyp}\label{H:vanishing} For any sufficiently large integer $n>0$, a general section of $\O(n)$ has smooth vanishing locus on $\ms X$.
\end{hyp}
This is clearly satisfied when $\ms X\to X$ is a $\m_n$-gerbe (or, more generally, a $G$-gerbe with $G$ a smooth group scheme).  In fact, \ref{H:vanishing} holds whenever $\ms X$ is smooth tame DM stack (a stacky Bertini theorem), a fact which will appear in a paper currently in preparation \cite{orbifoldmoduli}.

\begin{lem} The geometric Hilbert function is geometric: if $k\subset K$ is an 
extension of fields and $X$ is a smooth geometrically connected projective variety over $k$, then 
for any coherent $\ms X$-twisted sheaf $\ms F$, $P^{g}_{\ms F}=P^{g}_{\ms 
F\tensor K}$ as functions on $\Z$.
\end{lem}
\begin{proof} This follows from the fact that Chern classes of 
arbitrary (perfect) coherent sheaves pull back 
under Tor-independent maps.
\end{proof}

\begin{notn} Following the conventions of Huybrechts and Lehn 
\cite[\S1.2]{h-l}, we write 
    $$P^{g}_{\ms F}(m)=\sum_{i=0}^{\dim\ms F}\alpha_{i}(\ms 
    F)\frac{m^{i}}{i!}.$$
With this definition the coefficients $\alpha_{i}$ need not be 
integers (contrary to \cite[p.\ 10]{h-l}).  
\end{notn}

\begin{defn} Given a coherent sheaf $\ms F$ of 
dimension $d$ on $\ms X$, the 
\emph{geometric rank\/} of $\ms F$ is defined to be $$\rk\ms F:=\alpha_{d}(\ms 
F)/\alpha_{d}({\ms O}_{\ms X}).$$  The \emph{geometric degree\/} of $\ms F$ is defined to be 
$$\deg\ms F=\alpha_{d-1}(\ms F)-\rk(\ms F)\cdot\alpha_{d-1}({\ms O}_{\ms X}).$$
\end{defn}
\begin{lem} Suppose $k$ is infinite.  Given a coherent sheaf $\ms F$ on $\ms X$, for any integer $n>0$ there is a 
global section $\sigma$ of ${\ms O}(n)$ such that $\sigma:\ms F(-n)\to\ms F$ 
is injective.
\end{lem}
\begin{proof} The set $\operatorname{Ass}\ms F$ is finite and 
determined by its image in $X$.  Since ${\ms O}(n)$ is very ample, 
there is a 
section missing these finitely many points.  It is easy to see that 
any associated point of the kernel of $\sigma$ must then be contained 
in the zero locus of $\sigma$, contradicting the choice of $\sigma$ 
and \ref{L:ass lemma}. 
\end{proof}

\begin{lem} For any coherent twisted sheaf, $\deg P^{g}_{\ms F}=\dim\ms 
F$.  In particular, if $\ms F$ is torsion free then geometric rank of $\ms F$ is non-zero.
\end{lem}
\begin{proof}  
Using \ref{H:vanishing}, this is clear by induction and the previous lemma, once we have 
verified it when $\dim X=0$.  In this case, the geometric Euler 
characteristic is just the dimension of the fiber of $\ms F$ over any 
geometric point of $\ms X$ divided by the degree of the inertia stack $\ms I(\ms X)\to\ms X$.  Indeed, the pushforward map $\ms X\to\spec(k)$ sends $1$ to $1/\deg(\ms I(\ms X)/\ms X)$.  (The denominator is just the cardinality of the stabilizer of a geometric point of $\ms X$.)  
\end{proof}

\begin{remark}\label{R:vanish if P does} In particular, the geometric Hilbert function of $\ms F$ vanishes if 
and only if $\ms F=0$.  Furthermore, one sees that the geometric rank 
of $\ms F$ is precisely the rank of $\ms F$ as an ${\ms O}$-module.  
Unfortunately, one cannot show this by arguing that $\ms F$ and 
${\ms O}^{\rk\ms F}$ agree on a dense open substack, as this is false.  
Instead, one must appeal directly to the Hirzebruch-Riemann-Roch 
formula (and the computation \cite[3.2.2]{fulton} of Chern classes of a twist). 
We leave the details to the 
reader.  The geometric degree of $\ms F$ is related to the degree of 
$\det(\ms F)$ just as in the case of ordinary sheaves: 
$\alpha_{d-1}=\deg\det(\ms F)$ (so the geometric degree is arrived at 
by a linear transformation familiar from \cite[1.6.8ff]{h-l}).  This will aid us in 
comparing various notions of semistability and slope-semistability to 
their classical counterparts (as in Simpson's theory for 
semistability of modules for sheaves of algebras \cite[\S3]{simpson}).
\end{remark}

\begin{para} For the rest of this section, we will consider only the 
case where $\ms X\to X$ is a $\m_{n}$-gerbe with $n\in{\ms O}(X)^{\times}$ 
and $X$ is a smooth projective scheme over a Noetherian affine base 
$S$.  Generalizations of these results and those of the following section to the case of a smooth tame DM stack will be considered in an upcoming paper \cite{orbifoldmoduli}.

We start with a refinement of \ref{C:gabber de jong} better suited to 
the eventual study of stability.

\begin{prop}\label{P:de jong gabber makes me quiver} Given a $\m_{n}$-gerbe on a smooth 
projective morphism $\ms X\to X\to S$ with Noetherian affine base $S$, there is a locally free $\ms 
X$-twisted sheaf $\ms V$ of constant non-zero rank and trivial determinant.
\end{prop}
\begin{proof} The existence of $\ms V$ is a non-trivial result which 
holds on any (separated) scheme with an ample invertible sheaf.  We 
refer the reader to the work of de Jong \cite{dejong-gabber} for the (upcoming) 
details.  To make the determinant trivial, first consider $\ms W:=\ms 
V^{\oplus n}$.  Since $n|\rk\ms W$, we have $\ms L=\det\ms W\in\Pic(X)$.  
Now $(\ms L^{\vee}\oplus{\ms O}^{\oplus \rk\ms W-1})\tensor\ms W$ is a 
locally free twisted sheaf of trivial determinant.
\end{proof}
\end{para}

Recall that $\ms X$ \emph{has the resolution property\/} if every 
coherent sheaf on $\ms X$ is the quotient of a locally free sheaf 
\cite{totaro}.  The present virtue of \ref{P:de jong gabber makes me quiver} 
lies in the following corollary.
\begin{cor} For any affine $T\to S$, the stack $\ms X_{T}$ has the 
resolution property.
\end{cor}
\begin{proof} By \ref{P:q-coh splits}, the category of coherent sheaves on $\ms X$ breaks up 
according to the degree of twisting.  It suffices to show that 
coherent $\ms X$-twisted sheaves have the resolution property.  
Applying the fibered Morita equivalence $\shom(\ms V,\ \cdot\ )$ 
reduces us to showing that $\send(\ms V)$-modules on $X$ have the 
resolution property.  This follows from the fact that coherent sheaves 
on a projective morphism have the resolution property.
\end{proof}

\begin{prop}\label{P:hilb poly constant} If $\ms F$ on $\ms X$ is $S$-flat, then $P^{g}_{\ms 
F_{s}}$ is constant for all geometric points $s\to S$.
\end{prop}
\begin{proof} Let $\ms G^{\bullet}\to\ms F$ be a locally free 
resolution of $\ms F$.  As $\ms X\to S$ is smooth, $\ms G^{\bullet}$ 
may be taken to be a finite resolution.  If $\ms F$ is flat, then for 
any $s\to S$, the complex $\ms G^{\bullet}_{s}$ is a resolution of $\ms 
F_{s}$.  Thus, to prove that $P^{g}$ is constant for $\ms F$, it 
suffices by additivity to prove it when $\ms F$ is assumed locally 
free.  In this case, we may globally apply the splitting principle 
(noting that the base change which filters the sheaf produces another 
proper smooth $S$-flat family).  Thus, it is enough to show that given 
invertible sheaves $L_{1},\ldots,L_{d}$ on $\ms X$ (with $d=\dim\ms 
X/S$), the intersection 
product $c_{1}(L_{1})\cdot\cdots\cdot c_{1}(L_{d})$ is constant in 
fibers.  As $A(\ms X)_{\Q}=A(X)_{\Q}$, it suffices (by raising each 
$L_{i}$ to the $n$th tensor power and using multi-linearity) to prove this for 
invertible sheaves on $X$.  This is now a standard calculation using 
the fact that Euler characteristics are constant in a flat family.  
(In other words, we return to Kleiman's definition of intersection product 
using Snapper's lemma \cite[\S1]{badescu}, \cite{kleiman}, where the intersection number appears as 
a coefficient in a polynomial Euler characteristic.)
\end{proof}
Thus, given $P$, the substack $\Tw_{\ms X/S}(P)\subset\Tw_{\ms X/S}$ 
consisting of twisted sheaves with fixed geometric Hilbert 
polynomial $P$ is open (in fact, a union of connected components).  
Since we will verify shortly that $\Tw_{\ms X/S}$ is an algebraic stack, it will immediately follow that $\Tw_{\ms X/S}(P)$ is an algebraic stack.

\begin{para} Let $P$ be a fixed polynomial and $\ms E$ a fixed coherent $\ms 
X$-twisted sheaf.  We will briefly study the space of quotients of 
$\ms E$ with a fixed geometric Hilbert polynomial.
\end{para}

\begin{defn} Let $\underline{\quot}^{P}_{\ms X/S}(\ms E)$ denote the 
functor on affine $S$-schemes which assigns to $T\to S$ the set of 
subsheaves $\ms G\subset\ms E_{T}$ such that $\ms E_{T}/\ms G$ is 
$T$-flat with geometric Hilbert polynomial $P$ in every fiber over $S$.
\end{defn}

\begin{prop}\label{P:quot proper} The functor $\underline{\quot}^{P}_{\ms X/S}(\ms E)$ is 
represented by locally projective scheme $\quot^{P}_{\ms X/S}(\ms 
E)$ over $S$.
\end{prop}
\begin{proof} It follows by an easy application of Artin's 
representability theorem 
that $\underline{\quot}$ is representable by an algebraic space which 
satisfies the valuative criterion of properness.  This is checked in 
great detail in \cite{olsson-starr} in a slightly different context (which 
is sufficiently close to ours to be a complete proof in our case as 
well).  The only fact that remains to prove is that the functor 
$\underline{\quot}$ is bounded in the sense of \cite[1.7.5]{h-l}.  In other 
words, we need to show that there is a quasi-compact scheme surjecting 
onto the functor.

Let $\ms V$ be a locally free $\ms X$-twisted sheaf.  Given any $\ms 
G\subset\ms E_{T}$ as above, note that the inclusion $\shom(\ms V_{T},\ms 
G)\subset\shom(\ms V_{T},\ms E_{T})$ has $T$-flat 
cokernel $\ms C$.  
If we knew that the (classical) Hilbert polynomial of the fibers of $\ms C$ were 
always the same (as $\ms G$ varies), we would be done.  Unfortunately, this is highly 
unlikely.  However, since we do know the geometric rank and geometric 
degree of $(\ms E_{T}/\ms G)_{s}$, we know the rank and degree of 
$\ms C_{s}$.  Furthermore, we know that the $\ms C_{s}$ are quotients 
of a fixed sheaf $\ms U:=\shom(\ms V,\ms E)$.  By a result of 
Grothendieck \cite[1.7.9]{h-l}, 
we know that the set of quotients of $\ms U_{s}$ with slope bounded above 
is bounded.  A consideration of the proof of Huybrechts 
and Lehn [\textit{ibid\/}.] shows that the set of Hilbert polynomials 
appearing in such quotients is finite and independent of $s$.  Thus, 
as the geometric Hilbert polynomial is locally constant, we see 
that $\quot^{P}_{\ms X}(\ms E)$ is a union of finitely many connected 
components of finitely many schemes of quotients of 
$\send(\ms V)$-modules which are themselves projective over $S$.  This 
completes the proof.
\end{proof}
\begin{rem} The proof actually works (with slight modification) 
for any coherent $\ms E$ on $\ms X$ (independent of twisting).
\end{rem}

We end this section with a lemma which will be useful in section 
\ref{S:asymptotic properties} and which shows some of the 
similarities between classical sheaves and twisted sheaves on surfaces.  Let $X$ be a smooth projective surface.  First, a provisional definition.

\begin{defn} Given a coherent sheaf $\ms F$ on $\ms X$ of dimension $0$, the \emph{length\/} of $\ms F$ is 
$$\ell(\ms F)=\chi^g(\ms F).$$ 
\end{defn}
It readily follows from the definition of $\chi^g$ that if $\ms F$ is supported on a residual gerbe with rank $r$, then the length of $\ms F$ is $r$.  (The salient observation is that the pushforward $A^2(\ms X)\to A^2(X)$ is naturally identified with multiplication by $1/n$.  On the other hand, if one desires that $\ell_Y(\ms F_Y)=\deg(Y/\ms X)\ell(\ms F)$ for a finite flat map $Y\to\ms X$, then one should omit the factor of $[\ms I(\ms X):\ms X]$ in the formula.  For the purposes of this paper, this is immaterial, so we have normalized everything to yield integer-valued functions.  In future work \cite{orbifoldmoduli} this distinction will be important, and we will correspondingly alter the definition.)

Given an integer $\ell>0$, we will also write $\quot(\ms E,\ell)$ for $\quot^{\ell}_{\ms X/S}(\ms E)$ (to be consistent with existing notations in the untwisted category).

\begin{lem}\label{L:irred quot} Suppose $X$ is a smooth surface.  
If $\ms E$ is a locally free $\ms X$-twisted sheaf of 
rank $r$, then $\quot(\ms E,\ell)$ is irreducible of dimension 
$\ell(r+1)$.
\end{lem}
\begin{proof} We use the highly non-trivial fact that this is true 
when $\ms X$ is trivial \cite[6.A.1]{h-l}.  First, note that there is an 
open subspace of the $\quot$ corresponding to length $\ell$ quotients 
which are just $\ell$ distinct ``twisted lines'' (invertible sheaves supported on residual gerbes).  This open subspace is 
isomorphic to an \'etale $(\P^{r-1})^{\ell}$-bundle over 
$\operatorname{Sym}^{\ell}(X)\setminus\Delta$, where $\Delta$ is the 
multidiagonal, hence is irreducible (and has the right dimension).  It is thus enough to show that 
the entire $\quot$ is the closure of this open, which is the same as 
showing that any quotient may be deformed into a quotient with reduced 
support.  Let $\ms E\to\ms Q$ be any quotient of length $\ell$.  Write 
the support (with its natural scheme structure) of $\ms Q$ as $Z$ (which will be the preimage of a closed 
subscheme of $X$).  The quotient map is the same as a quotient $\ms 
E_{Z}\to\ms Q$.  Since $Z$ is a scheme of finite length over an 
algebraically closed field, we have $\Br(Z)=0$.  Let $\ms L$ be a 
twisted invertible sheaf on $Z$; any two invertible twisted sheaves are 
in fact mutually isomorphic.  Twisting down by $\ms L$, we see 
that $\ms E_{Z}\to\ms Q$ is the same thing as a surjection ${\ms O}_{Z}^{r}\to 
Q$.  (In other words, $\ms E\tensor\ms L^{\vee}\cong{\ms O}^{r}$.)  By the irreducibility of $\quot({\ms O}^{r},\ell)$, we know that 
there is a complete discrete valuation ring $R$ containing $k$ and a 
flat family of quotients ${\ms O}^{r}_{X\tensor R}\to\widetilde Q$ on $X\tensor 
R$ whose special fiber is ${\ms O}^{r}\to Q$.  The support $S$ of 
$\widetilde Q$ will be finite over $R$, and hence will be strictly 
Henselian.  Thus, $\Br(S)=0$, and we may choose an invertible twisted 
sheaf $\widetilde{\ms L}$ on $S$ (for the pullback of $\ms X$ to 
$X\tensor R$).  Since $S$ is semilocal, it follows that $\widetilde{\ms 
L}^{r}\cong(\ms E\tensor 
R)_{S}$.  Thus, twisting the quotient $\widetilde Q$ by 
$\widetilde{\ms L}$, we find an effective deformation of $\ms Q$ into 
a quotient with reduced support.
\end{proof}

\subsection{Algebraic moduli}\label{S:moduli}
In this section, we will show that the stack of 
twisted sheaves is algebraic (in the sense of Artin).  In 
the process, we will develop a theory of semistable twisted sheaves 
and study the relation to Geometric Invariant Theory.

We prove that the stack of twisted sheaves on a proper 
morphism $X\to S$ of finite presentation over an excellent (quasi-separated) 
base is algebraic.  
This sets the stage for a study of 
stability of twisted sheaves (in its Mumford-Takemoto and Gieseker 
forms) when $X\to S$ is projective and its use in producing GIT quotient stacks 
and corepresenting projective schemes for stacks of semistable twisted sheaves in 
\ref{S:hilbert polynomials and quot spaces} and \ref{S:semistability}.  The work 
on Gieseker stability will require the definition of a suitable Hilbert 
polynomial.  We define and study this polynomial and state 
a Riemann-Roch theorem \ref{C:hrr closed immersion} which will be useful at 
various points throughout this work.

In the special case where $S$ is affine and $X\to S$ is projective of 
relative dimension 2, 
we can use \ref{C:gabber de jong} to ``drastically simplify'' the 
situation by reducing it to work of Simpson on stability of modules 
for an algebra.  Indeed, once there exists a locally free twisted 
sheaf $\ms V$, the category of twisted sheaves becomes equivalent (by 
\ref{P:mini-morita}) 
with the category of modules for the Azumaya algebra $\send(\ms V)$ 
on $X$.  Simpson has 
considered moduli of modules \cite{simpson} quite generally; being careful, we can choose $\ms 
V$ so that the stability condition considered by Simpson agrees with 
the stability condition defined here.  (In fact, in arbitrary 
dimension such a Morita equivalence will always preserve 
slope-stability.)  We will use this 
technique to transport Simpson's GIT approach to the twisted setting 
on a surface and to prove some boundedness theorems in arbitrary 
dimension by appeal to classical results after a Morita 
equivalence.

However, we wish to emphasize that this approach is fundamentally 
incorrect.  While it is useful to have a Morita equivalence handy for 
transporting classical theorems, it is always better to 
work intrinsically when possible.  Working directly on the gerbe is also a step toward 
producing a satisfactory theory of sheaves and bundles on (at least 
DM) stacks.  In \cite{orbifoldmoduli} (in preparation) we will present a different method for studying semistable sheaves on arbitrary polarized smooth tame DM stacks, including boundedness questions.

\subsubsection{Abstract existence}
Let $X\to S$ be an algebraic 
space which is proper of finite presentation over a locally 
Noetherian scheme, 
and let $\ms X\to X$ be a fixed $\m_{n}$-gerbe, where $n$ is prime 
to $\ch(X)$.  Consider the 
$S$-groupoid $\ms T_{\ms X/S}$ which assigns to an affine scheme $\spec R\to S$ over 
$S$ the category whose objects are $R$-flat families of coherent $\ms 
X$-twisted sheaves.  (We reserve the notation $\Tw$ for twisted 
sheaves without embedded points; \ref{S:purity} above shows that 
$\Tw\subset\ms T$ is an open substack.)  Our goal in this section is to apply Artin's 
Theorem \cite{artin} to prove the following.

\begin{prop}\label{P:abstract existence} $\ms T_{\ms X/S}$ is an algebraic stack locally of finite 
presentation over $S$.
\end{prop}

\begin{lem}\label{C:general abstract existence} The result of 
\ref{P:abstract existence} is true if and only if it is true when $S$ 
is excellent and Noetherian.
\end{lem}
\begin{proof} Since 1) $X$ is of finite presentation, 2) being algebraic is 
local on $S$ and stable under base change, 3) the formation of \'etale cohomology is compatible 
with affine limits \cite[I.4]{freitag-kiehl}, and 4) the formation of the stack $\ms T$ 
is compatible with base change, we may replace $S$ with a finite type 
$\Z$-algebra.
\end{proof}

Most of the components necessary to apply Artin's Theorem are described 
in the deformation theory of \ref{S:deformations of twisted 
sheaves}.

\begin{lem}\label{L:existence theorem} Let $R$ be a complete local Noetherian 
ring, and suppose $S=\spec R$ above.  Given a compatible system of 
twisted sheaves $\ms F_{i}$ on $\ms X\tensor R/\mf m_{R}^{i+1}$, there 
is a twisted sheaf $\ms F$ on $\ms X$ whose reduction modulo $\mf 
m_{R}^{i+1}$ is compatibly isomorphic to $\ms F_{i}$.
\end{lem}
\begin{proof} This follows directly from the result of 
Olsson and Starr for sheaves on DM stacks (Proposition 2.1 of \cite{olsson-starr}), 
generalizing earlier work of Abramovich and Vistoli (appendix to 
\cite{abramovich-vistoli}).  (One could also use Olsson's general version of the existence theorem for Artin stacks, proved in \cite{olsson-chow} as a consequence of Chow's lemma for such stacks.)  If $X\to S$ is a projective 
morphism of schemes, then by \ref{C:gabber 
de jong} and Morita equivalence, the category of coherent twisted sheaves is 
equivalent to the category of $\send_{\ms X}(\ms V)$-modules, where 
$\ms V$ is a faithful locally free twisted sheaf.  But 
then we are reduced to the classical form of 
Grothendieck's Existence Theorem for modules over a coherent algebra 
\cite[\S5]{ega3-1}.
\end{proof}

\begin{lem}[Schlessinger]\label{L:schlessinger} Suppose $A_{1}\to A_{0}\leftarrow A_{2}$ is a diagram of commutative 
rings such that $A_{2}\to A_{0}$ is a surjection with nilpotent kernel $J$.  
Suppose give a diagram of flat modules $M_{1}\to M_{0}\leftarrow M_{2}$ 
over the diagram of rings inducing isomorphisms $M_{i}\tensor A_{0}\cong 
M_{0}$.  Let $B=A_{2}\times_{A_{0}}A_{1}$ and 
$N=M_{2}\times_{M_{0}}M_{1}$.  Then $N$ is a flat $B$-module and 
$N\tensor A_{i}\cong M_{i}$.
\end{lem}
\begin{proof} The proof of this result given in \cite{schl} only treats 
a special case which does not 
suffice for our purposes and the reference given there for the general 
case is not publicly available.  Thus, we give a proof which works for 
Noetherian rings and indicate how to generalize it to arbitrary 
commutative rings.

    To see that $N$ is $B$-flat, we use the local criterion of flatness 
\cite[\S22]{matsumura}.  Since $A_{2}\to A_{0}$ is surjective (say with kernel 
$J$), we see that $B\to A_{1}$ is surjective with (nilpotent) kernel 
$I:=J\times_{A_{0}}0_{A_{1}}$.  It is easy to see that $N/IN\cong M_{1}$ 
as $A_{1}$-modules.  To show that $N$ is flat over $B$, it remains to 
show that the natural map $\phi:I\tensor_{B} N\to IN$ is an isomorphism.  
We may assume (after filtering $J$ and proceeding inductively) that 
$J$ is generated by a single element $t$ of square $0$.  (This step 
of the proof only works in the Noetherian case, but the usual 
``equational criterion'' for flatness \cite[7.6]{matsumura} will work in the general 
case.  We choose to analyze this case for the sake of simplicity, and 
because it suffices for our purposes.)
The 
statement that $I\tensor N\to IN$ is an isomorphism is then 
equivalent to the statement that if $n=m_{2}\times m_{1}$ satisfies 
$(t\times 0)n=0$ then $(t\times 0)\tensor n=0$.  But if $(t\times 
0)n=0$, then $tm_{2}=0$.  As $M_{2}$ is flat over $A_{2}$, we have 
$m_{2}=tm_{2}'$, so that $m_{2}\mapsto 0\in A_{0}$.  Thus, 
$m_{1}\mapsto 0\in A_{0}$, so $m_{1}=\sum k_{j}m^{(j)}_{1}$ for some 
$k_{j}\in\ker(A_{1}\to A_{0})$, and so $m_{2}\times m_{1}=(t\times 
k_{1})(m_{2}\times m^{(1)}_{1})+(0\times k_{2})(m_{2}\times 
m^{(2)}_{2})+\cdots$.  Plugging this in, we find that $(t\times 
0)\otimes n=0$ as required.
\end{proof}

\begin{proof}[Proof of \ref{P:abstract existence}] We recall Artin's 
conditions: let $F$ be the stack of twisted sheaves, $\widebar F$ the 
associated presheaf of isomorphism classes.  Given a morphism of 
rings $B\to A$ and an element $a\in F(A)$, we will denote $F_{a}(B)$ 
the fiber of $F(B)\to F(A)$ over $a$ (and similarly for $\widebar 
F$).  The first conditions which must be satisfied to apply Artin's 
theorem are the Schlessinger-Rim criteria (our versions are slightly 
more general then are necessary; see \cite{artin} for Artin's list): 
\begin{enumerate}
    \item[(S1a)] given a diagram $A'\to 
A\leftarrow B$ with $A'\to A$ surjective with nilpotent kernel, and 
given $a\in F(A)$, the canonical map $$\widebar 
F_{a}(A'\times_{A}B)\to\widebar F_{a}(A')\times\widebar F_{a}(B)$$
is surjective.  

    \item[(S1b)] If $B\to A$ is a surjection, $b\in F(B)$ with image 
    $a\in F(A)$, and $M$ is a finite $A$-module then the canonical map 
    $$\widebar F_{b}(B\oplus M)\to\widebar F_{a}(A\oplus M)$$
    is bijective.

    \item[(S2)] Given $a\in F(A)$, the $A$-module $\widebar 
    F_{a}(A\oplus M)$ is finite.  (The module structure comes about 
    via S1b.  See \cite{artin, rim, schl} for details.)

    \item[(Aut)] Given $a\in F(A)$, the module $Aut_{a}(A\oplus M)$ of 
    infinitesimal automorphisms of $a$ is a finite $A$-module.
\end{enumerate}
In our case, these are easy to check.  (S1a) follows from 
\ref{L:schlessinger} by an argument similar to \cite[3.1]{schl}.  
(S1b) follows from the cher \`a Cartan 
isomorphism \ref{C:general cher}.  (S2) comes from the coherence of derived pushforwards of 
coherent sheaves on proper morphisms.  (Aut) follows from 
\ref{L:finite coho}.

In addition to the ``local versality'' conditions, one must check 
effectivization and constructibility conditions.  In particular, one 
must check that the map $F(\widehat A)\to\invlim F(\widehat A/\mf 
m^{n})$ is a 1-isomorphism of groupoids for a local Noetherian $A$ 
over $S$.  This follows from \ref{L:existence theorem} above.  The constructibility conditions are the 
following: the deformation and obstruction theories are compatible 
with \'etale localizations and completion (\ref{P:loc and constr}(1) 
and (2)), and there is a dense open where they are compatible with base 
change to fibers (\ref{P:loc and constr}(3)).  One requires that similar 
conditions hold for the group of infinitesimal automorphisms; this is 
also subsumed in \ref{P:loc and constr}.  The last condition to check 
is that given a reduced finite type $S$-affine $\spec A_{0}\to S$ and an 
element $a_{0}\in F(A_{0})$, any automorphism which induces the 
identity in the fiber at a dense set of points $A_{0}\to k$ of finite 
type over $S$ is the identity morphism.  This is local on $\ms X$, so 
it reduces to the case where $\ms X=X$ is affine and $\ms F$ is an $S$-flat 
coherent sheaf on $X$.  This reduces to showing that a section 
$\sigma$ of $\ms F$ which vanishes in fibers over a dense set of points of 
$\spec A_{0}$ is the zero 
section.  By flatness, the locus of points in $\spec A_{0}$ over 
which $\sigma$ vanishes is closed under specialization.  On the 
other hand, one easily sees that the set is constructible.  (The only 
non-trivial point comes in checking that if $A_{0}$ is integral and 
$\sigma$ does not vanish on the generic fiber, then 
there is an open subset of $\spec A_{0}$ consisting of fibers where 
$\sigma$ does not vanish.  There is an open subset $U$ of $X$ consisting 
of points $x\in X$ such that $\sigma_{\kappa(x)}\neq 0$, as $\ms F$ 
is coherent.  But $X\to\spec A_{0}$ is of finite type and $A_{0}$ is 
Noetherian, so the image of $U$ contains an open subset of $\spec A_{0}$ 
by Chevalley's theorem.)  Thus, the set of fibers where $\sigma$ 
vanishes is closed, so if it contains a dense set it is all of $\spec 
A_{0}$, as required.
\end{proof}

\subsubsection{Semistability and boundedness}\label{S:semistability}
We wish in this section to define a reasonable stability 
condition for twisted sheaves on smooth projective varieties.  
Variations on this theme occur throughout the study of moduli of 
sheaves.  The basic goal is to produce a condition which cuts out a 
well-behaved substack of the stack of pure sheaves.  Historically, 
this has meant two things: from the differential-geometric angle, 
stability conditions are related to the existence of certain 
types of metrics on bundles; from the algebro-geometric direction, 
the choice of a stability condition is influenced by the use of 
Geometric Invariant Theory to construct the moduli of such sheaves as 
a quotient stack.  

Using the geometric Hilbert polynomial, we define 
a stability condition for twisted sheaves analogous to the classical definition 
for untwisted sheaves.  As usual, a coarsening of our relation will define 
$\mu$-semistability (also known as Mumford-Takemoto semistability and slope semistability).  (In 
characteristic zero, this condition is probably equivalent to the existence of 
certain metrics on the associated analytic stack [orbifold].  We do 
not pursue this matter here.)  On surfaces, we relate our 
construction to GIT via a Morita equivalence and fundamental work of Simpson on 
moduli of modules for a sheaf of algebras \cite{simpson}.  (More generally, we 
make this comparison when there exists a locally free twisted sheaf 
with sufficiently many vanishing Chern classes.)

\begin{defn}\label{D:semistability} An $\ms X$-twisted sheaf $\ms F$ of dimension $d$ is \emph{semistable\/} 
(respectively \emph{stable\/}) if for any subsheaf $\ms G\subset\ms F$ 
we have $\alpha_{d}(\ms F)P^{g}_{\ms G}\leq\alpha_{d}(\ms G)P^{g}_{\ms F}$ 
(respectively $\alpha_{d}(\ms F)P^{g}_{\ms G}<\alpha_{d}(\ms G)P^{g}_{\ms F}$).
\end{defn}
\begin{lem} A semistable coherent $\ms X$-twisted sheaf $\ms F$ is 
pure.
\end{lem}
\begin{proof} If $\ms G\subset\ms F$ is a torsion subsheaf, then 
$\dim\ms G<\dim\ms F$, which means that $P_{\ms G}\leq 0$ (as 
$\alpha_{d}(\ms F)\neq 0$).  Thus, $P^{g}_{\ms G}=0$, and therefore $\ms 
G=0$ by \ref{R:vanish if P does}.
\end{proof}
\begin{defn} The \emph{reduced geometric Hilbert polynomial\/} of a 
coherent $\ms X$-twisted sheaf $\ms F$ of dimension $d$ is $p^{g}_{\ms 
F}:=(1/\alpha_{d})P^{g}_{\ms F}$.
\end{defn}
By definition, an $\ms X$-twisted sheaf $\ms F$ is semistable if and only 
if it is pure and for any subsheaf $\ms G\subset\ms F$ we have 
$p^{g}_{\ms G}\leq p^{g}_{\ms F}$.

\begin{defn} The \emph{slope\/} of a coherent $\ms X$-twisted sheaf 
$\ms F$ of dimension $d$ is 
$$\mu(\ms F):=\frac{\deg\ms F}{\rk\ms F}.$$
\end{defn}

\begin{defn} A coherent $\ms X$-twisted sheaf $\ms F$ of dimension 
$d$ is \emph{$\mu$-semistable\/} (respectively \emph{$\mu$-stable\/}) 
if $\ms F$ is pure and for any subsheaf $\ms G\subset\ms F$ we have 
$\mu(\ms G)\leq\mu(\ms F)$ (respectively $\mu(\ms G)<\mu(\ms F)$).
\end{defn}
\begin{remark} If we define the \emph{modified slope\/} of an $\ms 
X$-twisted sheaf $\ms F$ of dimension $d$ as $\muhat(\ms F):=\alpha_{d-1}/\alpha_{d}$, then 
we get the same notion of slope semistability as above.  We will 
use both notions of slope interchangeably.\end{remark}

\begin{defn} Given a sheaf of algebras $\ms A$ on $X$, an $\ms A$-module $\ms 
F$ is \emph{Simpson (semi)stable\/} if the inequality of 
\ref{D:semistability} holds for subsheaves $\ms G$ which are $\ms 
A$-submodules.
\end{defn}

\begin{lem}\label{L:semistab same as simpson} Let $X$ be a smooth projective 
variety and $\ms X$ a 
$\m_{n}$-gerbe on $X$ with a locally free twisted sheaf $\ms V$ of 
rank $v$ such 
that $c_{i}(V)=0\in A_{\Q}(\ms X)$ for all $1\leq i<n$.  
Let $\ms A=\send_{\ms X}(\ms V)$.  Then semistability of $\ms 
X$-twisted sheaves is equivalent to Simpson-semistability of $\ms 
A$-modules via the fibered Morita equivalence $\ms W\mapsto\shom(\ms V,\ms W)$.
\end{lem}
\begin{proof} By the Riemann-Roch theorem,
\begin{equation*}\begin{split}
    \chi(\shom(\ms V,\ms W))&=\deg(\chern(\ms V^{\vee})\cdot\chern(\ms 
    W)\cdot\td_{\ms X})\\
    &=v\deg(\chern(\ms W)\cdot\td_{\ms X})+\deg(\chern(\ms 
V^{\vee})_{1}\cdot(\chern(\ms W)\cdot\td_{\ms 
X})_{n-1})+\cdots\\
&\qquad\qquad\qquad\qquad\qquad\qquad\qquad\qquad\cdots+\rk(\ms W)\deg(\chern(\ms V^{\vee})\cdot\td_{\ms 
X})
\end{split}\end{equation*}
The assumption about the Chern classes of $\ms V$ kills all 
of the terms but the first and the last.  We see that $\chi(\shom(\ms 
V,\ms W))=v\chi^{g}(\ms W)+\rk(\ms W)\cdot\text{constant},$ whence the result 
follows.
\end{proof}

\begin{prop} Let $X$ be a smooth projective variety and $\ms X$ a 
$\m_{n}$-gerbe on $X$.  The category of $\muhat$-semistable 
coherent $\ms X$-twisted sheaves with fixed geometric Hilbert 
polynomial is bounded.
\end{prop}
\begin{proof} By \ref{P:de jong gabber makes me quiver}, there is a locally free 
$\ms X$-twisted sheaf $\ms V$ with $\det\ms V={\ms O}$.  
Applying a result of Simpson \cite[3.3]{simpson}, one sees that if $\ms F$ comes from a set of coherent twisted sheaves with fixed geometric Hilbert polynomial $P$, than the $\ms A:=\send(\ms V)$-module $\shom(\ms V,\ms F)$ has fixed slope and $\mu_{\text{max}}$ bounded above by a constant depending only upon $P$.

To show boundedness, first consider the subset of reflexive sheaves.  
Given a reflexive sheaf $F$ on $X$, temporarily write 
$$P^{g}_{F}(m)=\sum_{i=0}^{\dim X}a_{i}(F)\binom{m+\dim X-i}{\dim 
X-i}.$$
By a result of Langer (proving a theorem of Maruyama in arbitrary 
characteristic), the set of coherent reflexive 
sheaves $F$ on $X$ with a fixed upper bound on $\mu_{\text{max}}(F)$, 
$a_{0}(F)=r$, $a_{1}(F)=a_{1}$, and $a_{2}(F)\geq a_{2}$ for fixed 
$a_{0},a_{1},a_{2}$ is bounded \cite[4.3]{langer}.  Thus, to apply this to 
our situation,
it remains to show that the ``codimension 2'' coefficient $a_{2}$ is 
bounded below.  Looking at the formula in 
\ref{L:semistab same as simpson} and using the formula for the Chern 
character, we see that if $\det\ms V={\ms O}$ then 
the correction to the codimension 2 term of the Hilbert polynomial of $\shom(\ms V,\ms 
W)$ coming from higher terms (= after the first term) is given by 
$-\kappa(c_{2}(\ms V)\cdot c_{1}({\ms O}(1))^{d-2})t^{d-2}$, where $\kappa$ is 
a coefficient which depends only on $d$ (the dimension of $X$) and 
the degree of ${\ms O}(1)$ on $X$.  In particular, after dualizing $\ms V$ 
if necessary, we may assume that this correction term is always 
non-negative.  Thus, we see that the Morita equivalence we apply will 
yield Hilbert polynomials with bounded below codimension 2 terms.  
Applying Langer's theorem [\textit{ibid\/}.], we are done for reflexive twisted 
sheaves.

Given a torsion free twisted sheaf, taking its reflexive hull 
preserves $\mu$-semistability, fixes 
$a_{0}$ and $a_{1}$, and increases $a_{2}$.  Thus, we have just shown 
that the set of reflexive hulls of the sheaves we are interested in 
is bounded.  In particular, only finitely many geometric Hilbert 
polynomials occur.  We can now apply \ref{P:quot proper} finitely many 
times to yield the desired result.
\end{proof}
\begin{cor}\label{C:boundedness} The category of semistable $\ms X$-twisted 
sheaves with fixed geometric Hilbert polynomial is bounded.
\end{cor}
\begin{proof} It is elementary that any semistable sheaf is 
$\muhat$-semistable.
\end{proof}

\begin{cor}\label{C:openness of stuff} Let $X\to S$ be a smooth projective morphism, $\ms X$ a 
$\m_{n}$-gerbe on $X$, and $\ms F$ 
an $S$-flat family of coherent $\ms X$-twisted sheaves.  The locus 
of $\muhat$-semistable (resp.\ semistable, resp.\ geometrically 
$\muhat$-stable, resp.\ geometrically stable) fibers of $\ms F$ is 
open in $S$.
\end{cor}
\begin{proof} It suffices to prove this when $S$ is affine, whence we 
may assume that $\ms X$ is a quotient stack (\ref{C:gabber de jong} 
again!) and that the theory 
developed above applies.  Now one can apply \cite[2.3.1]{h-l} 
verbatim, with the additional remark that their proof also works 
for $\muhat$-semistability (even though they do not state this 
explicitly).\end{proof}

\subsubsection{Applications of GIT}
Let $(X,{\ms O}(1))$ be a (polarized) smooth projective variety over an algebraically closed 
field $k$ and $\ms X\to X$ a $\m_{n}$-gerbe with $n\in k^{\times}$.  
According to \ref{C:openness of stuff}, there is an algebraic stack 
$\Tw^{ss}_{\ms X/k}(r,P)$ of 
semistable twisted sheaves of fixed rank $r$ and geometric Hilbert 
polynomial $P$ containing an open substack $\Tw^{s}_{\ms X/k}(r,P)$ 
of geometrically stable points (which contains a further open substack of 
geometrically $\mu$-stable points).  We will make use of the algebraic Picard stack $\ms Pic_{\ms X/k}$ parametrizing invertible sheaves on the stack $\ms X$.  The reader uncomfortable with this (in fact, the methods used in this paper show it is algebraic) may feel free to consider only sheaves of rank $n$ in 
this section (in which case the determinant will in fact be a section of $\ms Pic_{X/k}$).  Recall that there is a determinant 1-morphism
$$\Tw_{\ms X/k}(n)\to\ms Pic_{\ms X/k}.$$
Let $L$ be an invertible sheaf on $\ms X$, corresponding to a 1-morphism 
$\phi_{L}:k\to\ms Pic_{\ms X/k}$.

\begin{defn} With the above notation, the \emph{stack of semistable twisted sheaves of rank 
$r$, determinant $L$ and geometric Hilbert polynomial $P$\/} is 
$$\Tw^{ss}_{\ms X/k}(r,L,P):=\Tw^{ss}_{\ms X/k}(r,P)\times_{\ms 
Pic_{\ms X/k},\phi_{L}}k.$$
The open substack of geometrically stable sheaves will be denoted $\Tw^{s}_{\ms 
X/k}(r,L,P)$.
\end{defn}

The usual computation \cite[2.2.2]{l-mb} of the 1-fiber product of stacks shows 
that $\Tw^{ss}_{\ms X/k}(r,L,P)$ has as objects over $T\to\spec k$ 
pairs $(\ms V,\phi)$, where $\ms V$ is a flat family of torsion free 
semistable $\ms X$-twisted sheaves parametrized by $T$, 
$\phi:\det\ms V\simto L_{T}$ is an isomorphism, and for all points 
$t\to T$ one has $P^{g}_{\ms V_{t}}=P$.  
As usual, isomorphisms in the groupoid are given by 
isomorphisms of the sheaves $\ms V$ which respect the 
trivializations $\phi$.  Combining \ref{P:abstract existence} with 
\ref{C:openness of stuff} shows that $\Tw^{ss}$ and $\Tw^{s}$ are 
algebraic stacks, locally of finite presentation over $k$, hence (as 
$\ms Pic$ is algebraic) the 
same is true for $\Tw^{ss}_{\ms X/k}(r,L,P)$ and $\Tw^{s}_{\ms 
X/k}(r,L,P)$.

\begin{lem} The stack $\Tw^{ss}(r,P)$ (resp.\ $\Tw^{s}(r,L,P)$) is quasi-compact and universally 
closed over $k$.  The substack $\Tw^{s}(r,P)$ (resp.\ 
$\Tw^{s}(r,L,P)$) is quasi-compact and separated over $k$.
\end{lem}
\begin{proof} The numerical properties of the geometric Hilbert 
polynomial allow for a transcription of Langton's proof 
\cite[\S2.B]{h-l}.  The uncomfortable reader may use the Morita equivalence of 
\ref{L:semistab same as simpson} to reduce this to \cite[\S4]{simpson} (but 
only when 
there exists a locally free twisted sheaf with sufficiently many vanishing Chern 
classes, e.g., if $X$ is a surface).
\end{proof}

Suppose $X$ is a surface.  Given $L$, fixing the geometric Hilbert 
polynomial of $\ms V$ is the same as fixing $\deg c_{2}(\ms V)$ by 
the Riemann-Roch formula.  In 
this case, we will often write $\Tw(r,L,c)$ in place of $\Tw(r,L,P)$ 
in order to align ourselves with the classical literature on surfaces.  
When all of the adornments are clear from context (or irrelevant), we will omit them from 
the notation.

Historically, moduli of semistable sheaves (and more generally modules) 
were studied using the tools of Geometric Invariant Theory, as 
developed in Mumford's thesis \cite{git}.  The basic consequence of these methods 
is a proof that $\Tw^{ss}$ is corepresented by a projective scheme; in 
fact, one can say quite a bit more about the corepresenting object 
using the full theory.  The philosophy adopted here is that 
the stack is really a more fundamental object.  (It is galling that 
the semistability of a sheaf still lacks a convincing explanation in 
intrinsic terms without recourse to GIT.  However, as we remind the 
reader, stable sheaves \emph{do\/} have a convincing description in 
terms of unitary connections in characteristic 0.  In fact, 
these bundles arose independently of GIT and it was only discovered 
later that they solve a GIT problem \cite{raynaud}.)  We will apply some of the 
classical results in this section to deduce GIT-like properties of our 
own moduli problem.  When the underlying projective variety is a 
surface, techniques of Simpson will yield the result for all of 
$\Tw^{ss}$.  In general, even without GIT, one can show that 
$\Tw^{s}$ has a coarse moduli space.

\begin{lem}\label{L:make coarse space} Let $\ms X$ be an algebraic stack, and suppose $\ms I(\ms 
X)\to\ms X$ is fppf.  Then the big \'etale sheaf $\sh(\ms X)$ associated to $\ms X$ is an 
algebraic space and $\ms X\to\sh(\ms X)$ is a coarse moduli space.
\end{lem}
\begin{proof} This is essentially the 
content of item $2$ of the appendix to \cite{artin}.  A similar 
construction may also be found in \cite{abramovich-corti-vistoli}.  
\end{proof}

\begin{prop} $\Tw^{s}(r,L,P)\to\sh(\Tw^{s}(r,L,P)$ is a $\m_{r}$-gerbe on an 
algebraic space of finite type over $k$
\end{prop}
The class of this $\m_{r}$-gerbe 
in $\H^{2}(X,\G_{m})$ is the famous ``Brauer obstruction'' to the 
existence of a tautological twisted sheaf with determinant $L$ on $\sh(\Tw^{s}(r,L,P)\times\ms X$.  

\begin{defn} The algebraic space $\mTw^{s}_{\ms 
X/k}(r,L,c):=\sh\Tw^{s}_{\ms X/k}(r,L,c)$ is the \emph{moduli space 
of stable twisted sheaves\/}.
\end{defn}
Thus, that $\mTw^{s}$ is an algebraic space, as we have seen above, is 
quite easy to prove using the proper abstract foundations.  
The interesting challenge is to show the existence of an 
ample invertible sheaf on $\mTw^{s}$.  This really does seem like a 
difficult problem.  Of course, by \ref{C:boundedness} if we apply a 
Morita equivalence there are only finitely many possible Hilbert 
polynomials occurring, so we see that it suffices to prove 
quasi-projectivity under the assumption that both the geometric 
Hilbert polynomial and the Morita-Simpson-Hilbert polynomial is 
constant in fibers.  Showing 
abstractly that there is an ample invertible sheaf is an interesting 
problem.  

\begin{prop}\label{P:GIT applies} Let $\ms X\to X$ be a $\m_{n}$-gerbe on a smooth 
projective variety of dimension $d$.  Suppose there is a twisted sheaf $\ms V$ such that all 
Chern classes but (possibly) $c_{d}(\ms V)$ are zero in $A(\ms 
X)_{\Q}$.  Then $\Tw^{ss}_{\ms X/k}(r,L,P)$ is a GIT quotient stack 
with stable sublocus $\Tw^{s}_{\ms X/k}(r,L,P)$.
\end{prop}
\begin{proof} By \ref{L:semistab same as simpson}, this reduces to 
work of Simpson \cite[\S4]{simpson}.
\end{proof}

\begin{cor} Given the hypotheses of \ref{P:GIT applies}, there is a 
morphism to a projective scheme $$\Tw^{ss}_{\ms X/k}(r,L,P)\to \mTw^{ss}_{\ms X/k}(r,L,P)$$ 
corepresenting $\Tw^{ss}_{\ms X/k}(r,L,P)$ in the category of schemes 
and an open subscheme $U\subset\mTw^{ss}$ such that the 
restriction of the morphism 
$\Tw^{ss}\to\mTw^{ss}$ to $U$ yields an isomorphism 
$\Tw^{s}\to\mTw^{s}\simto U$.
\end{cor}

\begin{ques} In the absence of a $\ms V$ with sufficiently many vanishing Chern 
classes, is it still true that the coarse moduli space $\mTw^{s}$ is 
quasi-projective?  Attempting to prove this in various na\"\i ve ways 
always leads one back to GIT.  If the space is quasi-projective, can 
one find a projective scheme corepresenting $\Tw^{ss}$ by taking 
a projective closure of $\mTw^{s}$?
\end{ques}

\begin{rem} The work described in section \ref{S:essentially trivial} gives an 
indication of how one might go about proving that $\mTw^{ss}$ exists 
and has projective components \emph{as long as there is a locally 
free $\ms X$-twisted sheaf, say $\ms V$\/}.  Indeed, in this case one is easily led to 
conjecture that the space $\mTw^{ss}$ is related by Mumford-Thaddeus-type 
flips to the space of Simpson-semistable $\send(\ms V)$-modules, which 
is shown to be a GIT quotient in \cite{simpson}.  The mechanism should 
be very similar to the notion of ``twisted stability'' investigated by 
Matsuki-Wentworth and Yoshioka, as described below.  We intend to 
return to this question in the future.
\end{rem}

\subsubsection{Essentially trivial gerbes}\label{S:essentially trivial}
In this section we describe the situation for a $\m_{n}$-gerbe $\ms 
X$ on a projective variety $X$ which is the gerbe of $n$th roots of an 
invertible sheaf.  These correspond to the kernel of the natural map 
$\H^{2}(X,\m_{n})\to\H^{2}(X,\G_{m})$.  If one chooses the ``correct'' 
polarization of $X$, then the stack of semistable twisted 
sheaves is canonically isomorphic to the stack of semistable sheaves 
on the underlying variety $X$.  These spaces have essentially been 
studied by Ellingsrud-G\"ottsche, Thaddeus, Yoshioka, and Matsuki-Wentworth, in 
the guise of ``twisted stability.''  These authors did not write in 
terms of gerbes, but rather investigated what happens when instead of 
computing the Hilbert polynomial of a torsion free sheaf $F$ one 
computes the Hilbert polynomial of $F\tensor{\ms O}(\alpha)$, where 
$\alpha$ is some $\Q$-divisor (with stability now being called 
``$\alpha$-twisted stability'').  We refer the reader to their work 
(\cite{yoshioka} and the references therein) 
for a detailed description of the situation (in characteristic 0); we 
will only use a small bit of the theory in what follows.
At the end of the section we will spend a few moments considering what happens when 
the base field is not algebraically closed.  Let $X\to \spec k$ be a 
geometrically connected smooth projective variety over a field.

\begin{defn} A $\m_{n}$-gerbe $\ms X\to X$ is \emph{(geometrically) essentially 
trivial\/} if the class $[\ms X]$ has trivial image in 
$\H^{2}(X,\G_{m})$ (respectively, 
$\H^{2}(X\tensor_{k}\widebar k,\G_{m})$).
\end{defn}

\begin{lem}\label{L:essentially trivial criterion} A gerbe $\ms X$ is 
essentially trivial if and only if there exists an invertible $\ms X$-twisted 
sheaf $\ms L$.
\end{lem}
\begin{proof} If $\ms X$ is essentially trivial, then the associated $\G_m$-gerbe is isomorphic to $\B{\G_m}$.  We can thus pullback a twisted line from $\B{\G_{m,\spec k}}$.  Conversely, if there is an invertible $\ms X$-twisted sheaf $\ms L$, then $\send(\ms L)\cong\ms O_{X}|_{\ms X}$ is (the pullback to $\ms X$) of an Azumaya algebra with Brauer class $[\ms X]$, whence $[\ms X]=0\in\H^2(X,\G_m)$.
\end{proof}

Given $\ms L$ as in \ref{L:essentially trivial criterion}, it follows that $\ms L^{\tensor n}$ will be the pullback of 
an invertible sheaf on $X$.  This observation allows to classify essentially trivial gerbes using the Kummer sequence.

\begin{defn}\label{D:gerbe of nth roots}  Let $\ms M$ be an 
invertible sheaf on $X$.  The \emph{gerbe of $n$th roots\/} of $\ms 
M$, denoted $[\ms M]^{1/n}$, is the stack whose objects over $T$ are 
pairs $(\ms L,\phi)$, where $\ms L$ is an invertible sheaf on 
$X\times T$ and $\phi:\ms L^{\tensor n}\simto\ms M$ is an isomorphism.
\end{defn}
It is immediate that $[\ms M]^{1/n}$ is a $\m_{n}$-gerbe.

\begin{prop}\label{P:essentially trivial} The cohomology of the Kummer sequence 
$1\to\m_{n}\to\G_{m}\to\G_{m}\to 1$ yields an exact sequence
$$0\to\Pic(X)/n\Pic(X)\to\H^{2}(X,\m_{n})\to\Br(X)[n]\to 0.$$
Under this identification the cohomology class 
of any essentially trivial gerbe $\ms X$ in $\H^{2}(X,\m_{n})$ equals the image of the class of 
$\ms L^{\tensor n}$, where $\ms L$ is an $\ms X$-twisted invertible sheaf.
\end{prop}
\begin{proof} The construction of the (moderately) long exact sequence 
in non-abelian cohomology shows that given $\ms M\in\Pic(X)$, the 
coboundary $\delta(\ms M)\in\H^{2}(X,\m_{n})$ is just $[\ms M]^{1/n}$.
Up to isomorphism, this gerbe depends 
only on the residue of $\ms M$ modulo $n\Pic(X)$.  The sequence shows 
that any essentially trivial gerbe has the form $\delta(\ms M)$ for some 
$\ms M$.  On $\delta(\ms M)$, there is a universal $n$th root $\ms L$ 
with $\ms L^{\tensor n}\simto\ms M$.  If $\ms N$ is any other 
invertible twisted sheaf then $\ms L\tensor\ms N^{\vee}$ is an 
untwisted invertible sheaf, say $\ms M'$, and one has $\ms L^{\tensor 
n}\tensor(\ms N^{\tensor n})^{\vee}\cong(\ms M')^{\tensor n}$.  Thus, 
the $n$th tensor power of any invertible twisted sheaf lies in the 
same class as $\ms M$ modulo $n\Pic(X)$.
\end{proof}

In the following proposition, we use the notion of twisted stability and from 3.2 of \cite{matsuki-wentworth}.  
Given a class $\gamma\in\Pic(X)\tensor\Q$, the $\gamma$-twisted Hilbert polynomial of a coherent sheaf $\ms F$ on $X$ is $n\mapsto\chi(\ms F\tensor\gamma(n))$, computed formally using the Riemann-Roch formula.  (This is what Matsuki and Wentworth use to define $\gamma$-twisted (semi)stability.)

\begin{prop}\label{P:essentially trivial no problem} Suppose $\ms X$ is essentially trivial.  Then there 
exists $\gamma\in\Pic(X)\tensor\frac{1}{n}\Z$ such that there is an 
isomorphism of $\Tw^{ss}_{\ms X/k}(n,L,P)$ with the stack 
$\Sh^{\gamma-ss}_{X/k}(n,L(n\gamma),P)$ of 
$\gamma$-twisted semistable sheaves on $X$ of rank $n$, determinant $L(n\gamma)$, and 
twisted Hilbert polynomial $P$.
\end{prop}
\begin{proof} There is some $\ms M$ such that $\ms X=[\ms M]^{1/n}$.  
Let $M$ be a universal $n$th root of $\ms M$ on $\ms X$.  The 
functor $\ms V\mapsto\ms V\tensor M^{\vee}$ yields an equivalence of 
categories from $\ms X$-twisted sheaves to sheaves on $X$.  It is easy to see that 
semistability of $\ms V$ as a twisted sheaf translates into $(1/n)\ms M$-twisted 
stability of $\ms V\tensor M^{\vee}$. 
\end{proof}

\begin{cor}\label{C:essentially trivial} Suppose $\ms X$ is geometrically essentially trivial. 
There is an isomorphism $$\Tw^{\mu}_{\ms X/k}(n,L)\tensor\widebar 
k\simto\Sh^{\mu}_{X/k}(n,L')\tensor\widebar 
k,$$ where the superscript $\mu$ denotes the open substacks of 
$\mu$-stable sheaves and $\ms X=[L\tensor (L')^{\vee}]^{1/n}$.  If $X$ 
is a surface, there is an isomorphism $$\Tw^{\mu}_{\ms X/k}(n,L,P)\tensor\widebar 
k\simto\Sh^{\mu}_{X/k}(n,L',Q)\tensor\widebar 
k,$$ with $Q$ an appropriate polynomial.
\end{cor}
\begin{proof} The first part follows from the fact that slope 
stability is independent of a Matsuki-Wentworth twisting by a 
$\Q$-divisor and is left to the reader.  To see that the Hilbert 
polynomial is constant when $X$ is a surface, note that the functor 
$\ms V\mapsto\ms V\tensor M^{\vee}$ described in \ref{P:essentially trivial no problem} fixes the 
determinant and also preserves the discriminant.  It follows that it 
sends $\ms X$-twisted sheaves with a given second Chern class $c$ to ordinary sheaves 
with a fixed second Chern class $c'$.  Fixing the determinant and 
the $c_{2}$ then fixes the Hilbert polynomial.
\end{proof}

\section{Curves and surfaces}\label{S:curves and surfaces}
In this section, we develop the theory of $\Tw^{ss}$ when the 
underlying variety $X$ is a curve or a surface.  Over an algebraically 
closed field, there is a guiding meta-theorem: \emph{Anything which happens 
in the theory of $\Sh^{ss}$ happens in the theory of $\Tw^{ss}$\/}.  
For curves, this is not just a meta-theorem: as we will show in 
section \ref{S:twisted sheaves on curves}, $\Sh^{ss}$ and $\Tw^{ss}$ 
are isomorphic (with the proper adornments added to the symbols).  For surfaces, there is not 
a similar direct comparison, but the classical structure theory for 
$\Sh^{ss}$ carries over to $\Tw^{ss}$.  In particular, as the second Chern 
class grows, $\Tw^{ss}$ becomes irreducible.  One can further compute 
examples on e.g.\ K3 surfaces, but we have unfortunately not included these examples 
here.  (However, Yoshioka has worked out quite a bit for K3 surfaces over the complex numbers.  See \cite{yoshioka3}.)  Despite the excessively abstract 
foundations, we thus have a reasonable understanding of the geometry of 
these moduli spaces for low-dimensional varieties over algebraically 
closed fields.  There are many gems from the untwisted world waiting to 
be properly twisted which we have not been able to include here.  
They will hopefully appear in future work.

When the base field is allowed to be non-algebraically closed, things get more 
interesting, and the stacks $\Tw^{ss}$ carry arithmetic information 
which $\Sh^{ss}$ knows nothing about.  The straightforward geometry of the moduli spaces can 
now be brought to bear on arithmetic problems.  We will exploit this extra 
information in \cite{period-index-paper} when we study the Brauer group 
of a surface over an algebraically closed field, a 
finite field, and a local field.  The work here also appears to be just 
the beginning of a possibly fruitful line of investigation.

\subsection{Twisted sheaves on curves}\label{S:twisted sheaves on 
curves}
We illustrate the theory developed up to this point with the example 
of semistable twisted sheaves on curves.  This serves two purposes: 
first, twisted sheaves are easy to understand.  Second, we will 
use the results mentioned in this section when we study semistable 
twisted sheaves on surfaces.

By a \emph{curve\/} $C$ we will always mean a proper smooth geometrically connected curve over a field $k$.  

\begin{remark}
Much of what we say here can be generalized to singular curves, but even the classical theory of sheaves has not been 
very well worked out in the non-smooth case.  Furthermore, it would be slightly more complex to develop the theory of the geometric Hilbert polynomial in this context, but it can be done using the theory of localized Chern classes on a smooth embedding of the singular gerbe.  We spare the reader the technical details in this work.  The 
reason to consider more general curves is that in the 
relative case it is nice to be able to handle degenerate fibers.  For 
example, if a surface $X$ carries a generically nice pencil 
$\widetilde X\to\P^{1}$, it is likely 
(usually necessary) that there will be singular fibers in the pencil.  
We would still like to relate the space of semistable twisted sheaves 
on $\widetilde{X}$ to the relative space of twisted sheaves of $X$ viewed as a 
family of curves over $\P^{1}$.  For the applications of the theory here and in \cite{period-index-paper}, it will suffice to consider only smooth curves, however, so we omit the more general theory for the sake of brevity.
\end{remark}

\subsubsection{A curve over a point}
Let $C\to\spec k$ be a curve over an algebraically closed field, and let 
$\ms C\to C$ be a $\m_{n}$-gerbe over $C$ with $n\in k^{\times}$.

\begin{lem} If $X$ is a scheme of dimension at most $1$ over an algebraically closed field, then $\Br(X)=0$.
\end{lem}
\begin{proof} We sketch the proof.  One first reduces to the case 
where $X$ is reduced.  (E.g., consider $0\to\ms 
I\to{\ms O}_{X}\to{\ms O}_{X_{\text{red}}}\to 0$.  Then $1\to 
1+I\to{\ms O}_{X}^{\ast}\to{\ms O}_{X_{\text{red}}}^{\ast}\to 1$ is exact, and 
taking cohomology we find 
$0=\H^{2}(X,I)\to\H^{2}(X,\G_{m})\to\H^{2}(X_{\text{red}},\G_{m})\to\H^{3}(X,I)=0$ is exact.)  
It then suffices to show that the 
Brauer group of any irreducible component vanishes (see 
\ref{L:abstract gluing}ff for the type of reasoning used in this 
argument).  This follows from Tsen's theorem.
\end{proof}

In other words, by \ref{L:essentially trivial criterion} there exists an invertible $\ms C$-twisted sheaf, 
say $\ms L$.  Recall that for any torsion free coherent sheaf $\ms G$ on $C$, one 
has $\deg\ms G:=\chi(\ms G)-\rk\ms G\chi({\ms O}_{C})$.  

\begin{defn} Given a $\ms C$-twisted sheaf $\ms F$, the \emph{degree\/} of $\ms F$ 
is $$\deg(\ms F):=n\deg c_1(\ms F).$$
\end{defn}
Here the degree of $c_1(\ms F)$ is computed by applying proper pushforward in rational Chow theorem from $\ms C$ to a point.  It is easy to see that multiplying by $n$ is necessary in order for the pullback from $C$ to $\ms C$ to preserve degrees.

\begin{defn} The \emph{slope\/} of $\ms F$, 
denoted $\mu(\ms F)$, is $\deg\ms F/\rk\ms F$.  The twisted sheaf $\ms F$ is 
\emph{(semi-)stable\/} if for every twisted subsheaf $\ms G\subset\ms 
F$, we have $\mu(\ms G)(\leq)\mu(\ms F)$. 
\end{defn}
It is easy to see that tensoring with 
$\ms L^{\vee}$ creates a bijection between the semistable $\ms 
C$-twisted sheaves of rank $r$ and degree $d$ and the 
semistable sheaves on $C$ with rank $r$ and degree $d-r\deg\ms L$.  
Note that this last number must be an integer.  In fact, $\deg\ms 
L\in\frac{1}{n}\Z$, so $d\in\frac{1}{\gcd{r,n}}\Z$.

Given a $\m_n$-gerbe $\ms C\to C$, let $\widebar{\delta}(\ms C)\in(1/n)\Z/\Z$ be the fraction corresponding to the image of $[\ms C]$ under $\H^2(C,\m_n)\simto\Z/n\Z\simto(1/n)\Z/\Z$.

\begin{prop}\label{P:on curves nada} With the above notation, the stack of semistable $\ms C$-twisted sheaves 
of rank $r$ and degree $d$
is non-canonically isomorphic to the stack of semistable 
sheaves on $C$ of rank $r$ and degree $d-r\widebar{\delta}$.
\end{prop}
\begin{proof} It well-known 
that the stack of semistable sheaves on $C$ of rank $r$ and 
degree $d$ is non-canonically isomorphic to the stack of semistable 
sheaves of rank $r$ and degree $d+nr$ for any $n$.  Using the above notation, it is straightforward to check that any invertible $\ms C$-twisted sheaf $\ms L$ will have degree $q+\widebar{\delta}$ 
for some integer $q$.  The result follows by combining these statements.
\end{proof}

In particular, it is a GIT stack (hence corepresented by a projective 
variety).

The usual structure theory for moduli spaces of semistable sheaves 
on smooth curves developed by Seshadri, Ramanan, Ramanathan, 
Narasimhan, Mumford, Newstead, etc., now carries over to the twisted 
setting.  (See \cite[Appendix 5C]{git} for a relatively exhaustive 
list of references.)  We omit proofs for the sake of brevity.
\begin{cor}\label{C:looks like it should} If $C$ is smooth of genus $g\geq 2$, the moduli space of semistable $\ms C$-twisted sheaves of 
rank $r$ and any fixed determinant is unirational of dimension 
$(r^{2}-1)(g-1)$.  The stack of semistable $\ms C$-twisted sheaves of 
rank $r$ and fixed degree $d$ is irreducible, smooth, and unirational over $k$ of dimension 
$r^{2}(g-1)+1$ at stable points.  The stable locus is a gerbe over a smooth unirational quasi-projective variety.
\end{cor}
Note that, as usual, even though the stack is smooth, its 
corepresenting GIT quotient need not be smooth away from the stable 
locus.

\begin{prop} Suppose $d-r\widebar{\delta}\in\Z$ and $r$ are relatively 
prime.  The open immersion $\Tw^{s}_{\ms C/k}(r,d)\inj\Tw^{ss}_{\ms 
C/k}(r,d)$ is an isomorphism.  In this case, $\mTw^{ss}$ is a smooth 
rational projective variety isomorphic to $\sh(\Tw^{ss})$.  There is 
a tautological sheaf $\ms F$ on $\mTw^{ss}\times\ms C$, and $\Pic(\mTw^{ss})\cong\Z$.
\end{prop}

\subsubsection{The relative case}

When $C$ is allowed to move over a base (or descend over a 
non-algebraically closed base field) things get more 
interesting.  In this section, we let $\pi:C\to S$ denote a proper flat 
morphism of finite presentation whose geometric fibers are curves as 
above, and we let $\ms C\to C$ be a $\m_{n}$-gerbe with 
$n\in{\ms O}_{S}(S)^{\times}$.  We fix a rank $r$ and a rational number $d$, the degree.

Note that the degree map on the relative Picard scheme induces a 
morphism $$\phi:\H^{0}(S,\R^{2}\pi_{\ast}\mu_{n})\to\H^{0}(S,\Z/n\Z)$$ 
which is a relative version of the map considered in the proof of \ref{P:on curves nada}:  
The 
image over a connected component $S'\subset S$ is equal to $n$ times the 
constant value for the minimal degree of an invertible $\ms C$-twisted sheaf on a  
fiber.  If $S$ is connected, write $\widebar{\delta}:=(1/n)\phi([\ms C])$ 
as above (where $\phi([\ms C])$ is chosen to lie between 0 and 
$n-1$).

\begin{lem}\label{L:line bundle exists} Let $C$ be a geometrically connected smooth proper curve over a separably closed field $k$ of characteristic exponent $p$ and $\ms C\to C$ a $\m_n$-gerbe.  If $(n,p)=1$, then there is an invertible $\ms C$-twisted sheaf of degree $\widebar{\delta}(\ms C)$.
\end{lem}
\begin{proof} It suffices to show the existence of an invertible twisted sheaf $\ms L$, as it is then clear (as in section \ref{S:essentially trivial}) that $\deg\ms L$ has fractional part $\widebar{\delta}$, so we can tensor with an untwisted invertible sheaf to bring the degree down to $\widebar{\delta}$.  Furthermore, showing that $\ms L$ exists is equivalent to showing that the Brauer class of $[\ms C]$ is $0$.   By Tsen's theorem, this is true of $\ms C\tensor\widebar k$, hence it is true for some finite extension $L/k$.  Since $k$ is separably closed, $[L:k]$ is a power of $p$.  The inflation-restriction sequence in Galois cohomology over the function field of $C$ then shows that some power of $p$ kills $[\ms C]\in\H^2(C,\G_m)$.  Since $n[\ms C]=0$ (as it comes from a $\m_n$-gerbe), the result follows.
\end{proof}

\begin{prop}\label{P:twisted is form} Suppose $S$ is connected.  The stack of semistable $S$-flat $\ms C$-twisted sheaves 
of rank $r$ and degree $d$ is an \'etale form of the stack of 
semistable $S$-flat sheaves of rank $r$ and degree 
$d-r\widebar{\delta}$.
\end{prop}
\begin{proof} We simply need to note that one can \'etale-locally on 
the base find an invertible twisted sheaf $\ms L$ on $\ms C$ of degree 
$\widebar{\delta}$.  (The obstruction to the gluing of these local 
invertible sheaves is the image of $[\ms C]$ in 
$\H^{1}(S,\R^{1}\pi_{\ast}\G_{m})$.)  The comparison is made by tensoring with $\ms L^{\vee}$; this 
will not change the $S$-flatness of the sheaf because it will not 
change its local structure.  Applying \ref{L:line bundle exists} gives such an $\ms L$ in the fiber over a separable closure of the residue field of any point $s\in S$.  Applying the deformation theory of section \ref{S:deformations of twisted sheaves}, it is easy to see that $\ms L$ extends to an \'etale neighborhood of $s$.
\end{proof}

If one takes $S$ to be the spectrum of a field $k$, then one can 
describe the space $\mTw^{ss}$ in terms of Galois twists of 
$\mSh^{ss}$.  Similarly, if $X\to S$ is a surface fibered 
over a curve, 
one can use relative stacks of twisted sheaves of 
rank $1$ (``twisted Picard spaces'') to reconstruct Artin's 
isomorphism between the Brauer group of $X$ and the 
Tate-Shafarevich group of the Jacobian of the generic fiber $X_{K(S)}$.  We have excluded 
these topics from this paper, as our main concern here is with the 
geometry of the moduli spaces and not arithmetic.
They are discussed in detail in \cite{mythesis} and touched on in 
\cite{period-index-paper}.  

\subsubsection{Moving twisted sheaves on curves}
Given a divisor moving in a surface and a twisted sheaf on the 
divisor, we can push it along the moving curve.  This gives us a way 
of connecting two stable twisted sheaves on linearly equivalent 
smooth divisors in a family.  Throughout this section, $X$ is a 
smooth projective surface over an algebraically closed field $k$ and 
$\ms X\to X$ is a fixed $\m_{n}$-gerbe on $X$.

\begin{prop}\label{P:moving sheaves on curves} Let 
$C_{0}$ and $C_{1}$ be smooth curves in $X$ and let $\ms P_{i}$ be the 
pushforward to $\ms X$ of a stable locally free twisted sheaf on 
$C_{i}\times_{X}\ms X$. If $C_{0}$ 
is linearly equivalent to $C_{1}$ and $P^{g}_{\ms X,\ms 
P_{0}}=P^{g}_{\ms X,\ms P_{1}}$ then there is an irreducible $k$-variety $S$, 
two points $s_{0},s_{1}\in S(k)$, and an $S$-flat family of 
$\ms X$-twisted sheaves $\ms F$ on $\ms X\times S$ such that $\ms F_{s_{i}}\cong\ms P_{i}$.
\end{prop}
\begin{proof} The idea is to push $\ms P_{0}$ along an embedded 
deformation of $C_{0}$ into $C_{1}$ and then move the image through 
the moduli space of twisted twisted sheaves on $C_{1}$.  We can 
actually do both simultaneously (which is more likely to yield an 
\emph{irreducible\/} parameter space for the family).

Since $C_{0}$ and $C_{1}$ are linearly equivalent, there is a flat Cartier divisor $\mc C\subset 
X\times\P^{1}\to\P^{1}$ such that $\mc C_{0}=C_{0}$ and $\mc 
C_{1}=C_{1}$.  (E.g., one can take the total space of the pencil of 
sections of ${\ms O}(C_{0})$ generated by $C_{0}$ and $C_{1}$.)  
Passing to an open subset $U\subset\P^{1}$ if necessary, we may assume 
$\mc C\to U$ is smooth.  Consider the stack $\ms M:=\Tw^{ss}_{\ms X\times\mc 
C/U}(n,P)$.  It is a classical result that the stack $\sh^{ss}_{C/k}(n,d)$ 
is an irreducible GIT quotient stack \cite[Appendix 5C]{git}.  Thus, applying 
\ref{P:twisted is form} and using quasi-properness, we see that $\ms 
M$ is irreducible and smooth over $U$ (and thus smooth over $k$).

Let $M\to\ms M$ be a 
smooth cover.  Write $M_{i},i=1,\ldots,t$ for the connected 
components of $M$.  Then each $M_{i}$ is an open irreducible subspace of $M$, 
hence has open image in $\ms M$.  Since $\ms M$ is irreducible, there 
is some $i$ such that $M_{i}\to\ms M$ is surjective.  In 
other words, $\ms M$ has an irreducible smooth cover.  Choosing points 
$m_{0},m_{1}\in M(k)$ mapping to $\ms P_{0}$ and $\ms P_{1}$ 
respectively, we see that we can make a family of semistable sheaves 
on $\mc C\times_{U}M$ containing $\ms P_{0}$ and $\ms P_{1}$.  
Since $\mc C\subset X\times U$, we see that $\mc 
C\times_{U}M\subset X\times M$.  Pushing forward the family yields 
the result.
\end{proof}

\begin{cor}\label{C:connecting picard mamas} The conclusion of \ref{P:moving 
sheaves on curves} holds when $\ms P_{i}$ are invertible twisted 
sheaves, without explicit stability hypotheses.
\end{cor}
\begin{proof} This follows from the fact that any 
invertible sheaf is stable and the fact that 
$\Sh^{s}_{C/k}(1,d)=\Pic_{C/k}^{d}$ is smooth and irreducible.
\end{proof}

\subsubsection{Moduli of restrictions}
We use the above machinery to study what happens when restricting 
stable twisted sheaves on a surface $X$ to a very ample smooth curve 
$D$.  In 
particular, we show that there are no positive-dimensional complete families of locally 
free stable twisted sheaves on $X$ which all restrict to the same stable 
twisted sheaf (up to isomorphism) on $D$.  This will ultimately be used to show that 
asymptotically, the irreducible components of the stack of 
semistable twisted sheaves on $X$ contain both locally free and 
non-locally free points.

Throughout this section, $X$ is a smooth projective surface (with 
fixed very ample invertible sheaf ${\ms O}(1)$) over an 
algebraically closed field $k$ and $\ms X\to X$ is a $\m_{n}$-gerbe 
on $X$, with $n\in k^{\times}$.


\begin{notation} Whenever $\ms S$ is a stack of sheaves (on a curve or surface), we will let $\ms S_{lf}$ denote the open substack parametrizing locally free sheaves.
\end{notation}

Let $D\in|{\ms O}(1)|$ be a general member.
\begin{situation} Let $C$ be a smooth projective curve over $k$ and 
$\phi:C\to\Tw^{s}_{\ms X/k,lf}$ a 1-morphism to the locally free locus 
corresponding to $\ms F$ on $C\times X$.  Suppose that 
every object $\phi(c)$, $c\in C(k)$, restricts to a fixed stable locally 
free $\ms X_{D}$-twisted sheaf $\ms F_{0}$.  
\end{situation}

\begin{defn} A divisor $D\subset X$ is \emph{$\phi$-sticky\/} if there exists a simple locally free $\ms X_D$-twisted sheaf $\ms F_0$ such that for every $c\in C(k)$, the object $\phi(c)|_D$ is isomorphic to $\ms F_0$.
\end{defn}
We will suppress the $\phi$ from the notation when it is clear (essentially, for the rest of this section).

\begin{prop}\label{P:no curves inside}
$\phi$ is essentially constant (i.e., isotrivial).
\end{prop}
(Since $\phi$ lands in the stable locus, being isotrivial is 
equivalent to the map to the coarse moduli space $\mTw^{s}$ being 
constant.  Indeed, if $\phi$ is isotrivial, then there is a finite 
\'etale extension $C'\to C$ such that the induced map $C'\to\mTw^{s}$ 
is constant, whence the original map must be constant.  Conversely, 
if $C\to\mTw^{s}$ is constant, then $\phi$ lands in the fiber $\ms T$ of 
$\Tw^{s}\to\mTw^{s}$ over a point.  Since $\ms T$ is a $\m_{n}$-gerbe, 
the map $pt\to\ms T$ is finite \'etale, and pulling back by this map 
yields a finite \'etale cover $C'\to C$ such that the restriction of 
$\phi$ to $C'$ is constant.)

\begin{lem}\label{L:open set of yummies} There is an open subset of $|{\ms O}(1)|$ consisting of smooth sticky
divisors $D'$.
\end{lem}
\begin{proof} Write $P$ for $|{\ms O}(1)|$.  Let $I\subset X\times P$ be 
the incidence correspondence of ${\ms O}(1)$; the fiber of the second 
projection over a point $p\in|{\ms O}(1)|$ is the divisor corresponding to 
$p$.  The family $\ms F$ on $C\times X$ corresponding to $\phi$ pulls back to give a flat family 
$\F$ of twisted sheaves on $C\times I\to C\times P$.  (The 
sheaf $\F$ is flat by e.g.\ a Hilbert polynomial calculation after 
applying a Morita equivalence.)  The stickiness condition on $D$ says that the locus $\Psi$ of stable 
fibers contains all of $C\times\{[D]\}$.  By openness of stability and 
properness of $C$, we conclude that there is an open $U\subset P$ such 
that $C\times U\subset\Psi$.  We have a map $C\times 
U\to\Tw^{s}_{I_{U}/U}$ over $U$ such that the fiber $C\times\{[D]\}$ 
collapses to a point in the fiber $\mTw^{s}_{D/k}$.  By the usual 
rigidity lemma, it follows that an open subset of fibers get 
collapsed in the map to $\mTw^{s}_{I_{U}/U}$.  By Tsen's theorem, 
there is an open set of sticky divisors $D'$.  (The 
careful reader will note that the formation of the coarse moduli 
space $\mTw^{s}$ commutes with arbitrary base change in this case 
because $\Tw^{s}\to\mTw^{s}$ is a gerbe, so $\mTw^{s}$ is equal to the 
sheafification of $\Tw^{s}$ in the big \'etale topology, which is 
tautologically of formation compatible with 
base change.)
\end{proof}

\begin{lem}\label{L:structure of one yummy} Suppose $D$ is sticky.  The twisted sheaf $\ms 
F_{C\times D}$ has the form $\pr_{1}^{\ast}(\ms 
M)\times\pr_{2}^{\ast}(\ms F_{0})$, where $\ms M$ is an invertible 
sheaf of ${\ms O}_{C}$-modules and $\ms F_{0}$ is a stable twisted sheaf on $D$.
\end{lem}
\begin{proof}  Write $\ms D:=\ms X\times_{X}D$.  Since $D$ is sticky, the family $\ms F_{C\times D}$ gives rise to a diagram
    $$\xymatrix{C\ar[r]^{\widetilde\phi}\ar[dr]_{\phi} & \Tw^{s}_{\ms 
    D/k}\ar[d]^{\pi}\\
    & \mTw^{s}_{\ms D/k}}$$
such that $\phi$ is constant with value $[\ms F_{0}]$.  There is 
also a constant lift $\psi$ of $\phi$ given by the family 
$\pr_{2}^{\ast}\ms F_{0}$ on $C\times D$.  Since $\pi$ is a 
$\G_{m}$-gerbe, we see that $\psi$ and $\widetilde\phi$ are 
identified with two sections of a trivial $\G_{m}$-gerbe.  Using one 
of them to trivialize the gerbe, they differ by a map 
$C\to\B{\G_{m}}$, which gives the invertible sheaf $\ms M$.
\end{proof}

\begin{defn} A (possibly singular) divisor $E$ is called \emph{$\phi$-slippery\/} if 
there is an invertible sheaf $\ms M$ on $C$ and a fixed twisted sheaf $\ms 
F_{0}$ on $E$ such that $\ms F_{C\times E}\cong\pr_{1}^{\ast}\ms M\tensor\pr_{2}^{\ast}\ms 
F_{0}$.  
\end{defn}
We will similarly suppress the $\phi$ from the notation when it is clear from context.  We just showed that for a smooth divisor $D$, if $D$ is sticky then it is slippery, and that if $D$ is sticky for one smooth very ample divisor, 
then an open set in $|D|$ parametrizes sticky points.  Using these two facts, we now provide an 
inductive procedure for enlarging a slippery divisor $D$.  

\begin{lem}\label{L:abstract gluing} Let 
    $$\xymatrix{A\ar[r]\ar[d] & B\ar[d] \\
    C \ar[r] & D}$$
be a Cartesian diagram of surjections of sheaves of groups in a topos $T$.  The 
natural map $\B{A}\to\B{B}\times_{\B{D}}\B{C}$ is a 1-isomorphism of 
classifying stacks.
\end{lem}
\begin{proof} The natural map $\B{A}\to\B{B}\times_{\B{D}}\B{C}$ is 
given by sending a right $A$-torsor $F_{A}$ to the triple 
$(F_{A}\times^{A}B,F_{A}\times^{A}C,\phi)$, where 
$\phi:(F_{A}\times^{A}B)\times^{B}D\simto(F_{A}\times^{A}C)\times^{C}D$ 
is the natural isomorphism arising from the associativity of the 
contracted product (i.e., $(F_{A}\times^{A}B)\times^{B}D\simto 
F_{A}\times^{A}(B\times^{B}D)\simto F_{A}\times^{A}D$ and similarly 
for $C$).  There is a 1-morphism in the other direction arising as 
follows.  An object of $\B{B}\times_{\B{D}}\B{C}$ is given by a triple 
$(F_{B},F_{C},\psi)$, where $\psi:F_{B}\times^{B}D\simto 
F_{C}\times^{C}D$ is an isomorphism of right $D$-torsors.  Given such 
an object, one can produce a right $A$-torsor $F_A$ by forming the 
fiber square 
$$\xymatrix{& F_{A}\ar[dl]\ar[dr] &\\
F_{B}\ar[d]\ar[dr] & \Box & F_{C}\ar[dl]\\
F_{B,D}\ar[r]^{\psi} & F_{C,D} &}$$
where $F_{B,D}:=F_{B}\times^{B}D$, etc.  That $F_{A}$ is in fact an 
$A$-torsor follows from the surjectivity of $B\to D$.  We leave it as an exercise to  
check that these maps of stacks are 1-inverse to one another.
\end{proof}

We can use \ref{L:abstract gluing} to prove a (twisted) classical result about vector bundles on a 
union of curves meeting transversely.  Let $D$ and $D'$ be curves 
with transverse intersection $D\cap D'=\{q_{1},\ldots,q_{r}\}$.  Let 
$X$ be a $k$-scheme.  The transversality of the intersection of $D$ 
and $D'$ says that the diagram of surjections of sheaves of 
rings on $X\times(D\cup D')$
$$\xymatrix{{\ms O}_{X\times(D\cup D')}\ar[r]\ar[d] & {\ms O}_{X\times D'}\ar[d]\\
{\ms O}_{X\times D}\ar[r] & {\ms O}_{X\times(D\cap D')}}$$
is Cartesian, where all schemes are given their reduced structures.  
(More generally, given a ring ${\ms O}$ in a topos and two ideals $I$ and 
$I'$ such that $I\cap I'=0$, one has a corresponding diagram.  For 
non-CM schemes, there can be complex information at embedded 
intersection points.)  
Here we write (by abuse of notation) ${\ms O}_{D}$ for the pushforward of the 
structure sheaf of $D$ and similarly for $D'$ and $D\cap D'$.
It follows that given any $k$-scheme $X$ the diagram
$$\xymatrix{\GL_{n}{\ms O}_{X\times(D\cup D')}\ar[r]\ar[d] & 
\GL_{n}{\ms O}_{X\times D'}\ar[d]\\
\GL_{n}{\ms O}_{X\times D}\ar[r] & \GL_{n}{\ms O}_{X\times(D\cap D')}}$$
is a Cartesian diagram of surjections of sheaves of groups on 
$X\times(D\cup D')$.

Suppose $\ms V$ and $\ms V'$ are locally free sheaves of rank $n$ on 
$D$ and $D'$, respectively.  Our goal is to describe the space of locally free 
sheaves $\ms W$ on $D\cup D'$ which restrict to $\ms V$ on $D$ and 
$\ms V'$ on $D'$.  

Define a stack
$\Sigma$ on $k$-schemes as follows.  Given a 
$k$-scheme $X$, the fiber category $\Sigma_{X}$ is the groupoid of triples $(\ms 
W,\alpha,\beta)$ where $\ms W$ is locally free of rank $n$ on $X\times(D\cup 
D')$ and $\alpha:\ms W|_{X\times D}\simto\ms V_{X\times D}$ and 
$\beta:\ms W|_{X\times D'}\simto\ms V'_{X\times D'}$ are 
isomorphisms.  It is easy to see that in fact this groupoid is 
discrete, i.e., $\Sigma$ is the stack associated to a sheaf.

\begin{prop}\label{P:sheafified gluing} With the above notation, there 
is an isomorphism 
$$\Sigma\simto\isom_{D\cap D'}(\ms 
V|_{D\cap D'},\ms V'|_{D\cap D'}).$$
\end{prop}
\begin{proof} By \ref{L:abstract gluing} and transversality (as 
discussed above), we have a 1-isomorphism of stacks
$$\B{\GL_{n}{\ms O}_{D\cup 
D'}}\simto\B{\GL_{n}{\ms O}_{D}}\times_{\B{\GL_{n}{\ms O}_{D\cap D'}}}\B{\GL_{n}{\ms O}_{D'}}.$$
This shows that the functor sending $(\ms W,\alpha,\beta)$ to 
$\alpha\circ\beta^{-1}$ defines an isomorphism $\Sigma\to\isom$.  The 
details are left to the reader.
\end{proof}
\begin{cor}\label{C:untwisted gluing} Suppose $\ms V$ and $\ms V'$ are locally free simple sheaves 
of rank $n$.  The moduli space of locally free sheaves $\ms W$ of rank $n$ on 
$D\cup D'$ such that $\ms W|_{D}\cong\ms V$ and $\ms W|_{D'}\cong\ms 
V'$ is isomorphic to $\GL_{n}^{r}/\G_{m}$, where $\G_{m}$ is embedded 
along the diagonal.  Moreover, this scheme is affine.
\end{cor}
\begin{proof} That the scheme is affine follows from the fact that the 
quotient is the complement of a hypersurface (cut out by the product of 
the determinants) in a projective space.  Since $\ms V$ and $\ms V'$ 
are simple, it is easy to see that the moduli space $M$ parametrizing $\ms W$ 
restricting to $\ms V$ and $\ms V'$ exists as an algebraic space.  
Furthermore, there is a surjection $\Sigma\to M$ which is a 
$\G_{m}$-bundle, and in fact $M$ is identified with 
$\Sigma/\G_{m}$, where $\G_{m}$ acts in the natural way on the 
isomorphism $\beta$.  Applying \ref{P:sheafified gluing} completes the 
proof.
\end{proof}
\begin{cor}\label{C:twisted gluing} Suppose $\ms C\to D\cup D'$ is 
a $\m_{n}$-gerbe and $\ms V$ and $\ms V'$ are locally free simple 
twisted sheaves of rank $n$ on $\ms C_{D}$ and $\ms C_{D'}$.  The 
moduli space $M(\ms V,\ms V')$ of locally free twisted sheaves $\ms W$ of rank $n$ on 
$\ms C$ such that $\ms W|_{\ms C_{D}}\cong\ms V$ and $\ms W|_{\ms 
C_{D'}}\cong\ms V'$ is (non-canonically) isomorphic to $\GL_{n}^{r}/\G_{m}$.
\end{cor}
\begin{proof} This follows from \ref{C:untwisted gluing} after 
twisting down by a $\ms C$-twisted invertible sheaf.
\end{proof}

\begin{lem}\label{L:add a yummy} Suppose $D$ and $D'$ are (not necessarily smooth) slippery elements of 
$|{\ms O}(1)|$ which intersect transversely.  
Then $D\cup D'$ is slippery.
\end{lem}
\begin{proof}
In the decompositions $\ms F_{C\times D}\cong\pr_{1}^{\ast}\ms 
M\tensor\pr_{2}^{\ast}\ms F_{0}$ and $\ms F_{C\times D'}\cong\pr_{1}^{\ast}\ms 
M'\tensor\pr_{2}^{\ast}\ms F_{0}'$, we claim that $\ms M\cong\ms 
M'$.  Indeed, let $q\in D\cap D'$ be a point.  Restricting $\ms F$ 
to $C\times \{q\}$ and using the two decompositions, we find that 
$\ms M\tensor(\ms F_{0}\tensor\kappa(q))\cong\ms M'\tensor(\ms 
F_{0}'\tensor\kappa(q))$.  Both $\ms F_{0}\tensor\kappa(q)$ and $\ms 
F_{0}'\tensor\kappa(q)$ are non-zero finite-dimensional $\kappa(q)=k$-vector 
spaces of dimension $r$.  Thus, we conclude that both $\ms M\tensor\ms M'^{-1}$ and 
$\ms M'\tensor\ms M^{-1}$ have non-zero global sections, whence $\ms 
M\cong\ms M'$.  Choosing such an isomorphism and twisting down by 
$\pr_{1}^{\ast}\ms M$, there results a map from $C$ to the moduli 
space $M(\ms V,\ms V')$ of \ref{C:twisted gluing}, with $\ms V=\ms 
F_{0}$ and $\ms V'=\ms F_{0}'$.  Since $M(\ms V,\ms 
V')$ is affine and $C$ is proper, the map $C\to M(\ms V,\ms V')$ must 
be constant.  As moduli of simple sheaves are a $\G_{m}$-gerbe over moduli and $C$ is a 
curve over an algebraically closed field, Tsen's theorem shows that 
the family $\pr_{1}^{\ast}\ms M^{\vee}\tensor\ms F_{C\times (D\cup 
D')}$ is constant.  Thus, $D\cup D'$ is slippery.
\end{proof}

\begin{proof}[Proof of \ref{P:no curves inside}]   
Note that $\Tw^{s}$ is a $\G_{m}$-gerbe over its moduli 
space $\mTw^{s}$.  This means that any curve $C$ in $\mTw^s$ admits a 1-morphism 
$C\to\Tw^{s}(\ms X)$ lifting the inclusion $C\inj \mTw^s$.  Replacing $C$ 
by the normalization of the lift of its image in $T$, we may assume 
that the map $C\to\Tw^{s}$ is separably generated.  Thus, to show that 
it is essentially constant, it suffices to show that the map on 
tangent spaces is the zero map, i.e., that the first-order 
deformations of any point in $C$ induce the trivial deformation of 
the image point in moduli.  We will do this by showing that they 
induce the trivial deformation on a sufficiently ample divisor.  
It is easy to see that given a locally 
free twisted sheaf $\ms G$ on $X$, the space of first-order infinitesimal deformations of $\ms G$ 
which restrict to the trivial deformation on an effective divisor $D$ 
is principal homogeneous under the kernel of the restriction map 
$\H^{1}(X,\send(\ms G))\to\H^{1}(D,\send(\ms G_{D}))$; in the case 
where $\ms G$ and $\ms G_{D}$ are simple, this is precisely 
$\H^{1}(X,\send(\ms G)(-D))$.  Thus, if $D$ is sufficiently ample, 
the deformations of $\ms G$ inject into the deformations of $\ms G_{D}$.  By \ref{L:open set of yummies} and 
\ref{L:add a yummy}, we see that 1) $R(D)$ holds for some $D$ in 
${\ms O}(1)$, and 2) when $R(D)$ holds for some 
$D\in|{\ms O}(1)|$, there is an arbitrarily ample 
divisor $D^{(n)}=D_{1}\cup D_{2}\cup\cdots\cup D_{n}\in|{\ms O}(n)|$ such that 
$R'(D^{(n)})$ holds.  But $R'(D^{(n)})$ says precisely that the 
infinitesimal deformation of $\ms F_{D^{(n)}}$ induced by a 
tangent vector $t$ of $C$ is trivial.  As $D^{(n)}$ is arbitrarily ample, 
we see that the deformation of $\ms F$ induced by $t$ is also trivial.
\end{proof}

\subsection{Twisted sheaves on surfaces}
In this section, we discuss the moduli of twisted sheaves on surfaces.  In the 
process, we develop tools to reduce certain twisted statements to 
their classical counterparts.  This should be viewed as a preliminary survey of a 
theory which is certainly amenable to significant further 
development.  In particular, ongoing work of Langer 
(\cite{langer-castelnuovo}) should help 
clarify the classical situation in positive characteristic (and therefore, in 
our view, in characteristic 0 as well), and we believe that his methods 
will ultimately prove useful in the twisted case.

Throughout, we focus on moduli of twisted sheaves of rank $n$.  This 
is technically simpler, as then determinants naturally take values in the 
Picard group of $X$ itself.  This is also the case one is naturally 
led to consider when approaching the classification of (generalized) 
Azumaya algebras of degree $n$ in a Brauer class of order $n$, which 
is the most natural (and na\"\i ve) thing to do on a surface.  In 
general, if one wants to consider rank $r$ twisted sheaves on a 
$\m_{n}$-gerbe $\ms 
X$, then there is a $\m_{r}$-gerbe $\ms X_{r}$ carrying them all with the same 
Brauer class as $\ms X$.  Moreover, the natural map $\ms X_{r}\to\ms 
X_{n}$ serves to identify the stacks of \emph{semistable\/} sheaves via pullback.  
Thus, we lose nothing by assuming that $r=n$.

In order to orient the reader, we sketch the contents of this 
section: in \ref{S:discrim} we discuss the discriminant and estimates for 
the dimension of $\Tw^{ss}$ (at stable points).  In \ref{S:preparation} and 
\ref{S:restriction theorems}, we then discuss how (semi)stability behaves under 
restriction to curves in the surface, giving twisted forms of results 
of Langer \cite{langer}.  We prove a Bogomolov inequality for twisted 
sheaves and use it to re-prove a result of Artin and de Jong bounding 
the second Chern class of an Azumaya algebra.  
Finally, in \ref{S:asymptotic properties}, we generalize 
O'Grady's results on asymptotic smoothness and irreducibility to the 
space of twisted sheaves on an optimal gerbe (\ref{D:optimal}).  The 
restriction to optimal gerbes is made in order to have results which 
hold in all characteristics.  As we discuss below, it is likely that 
this restriction is unnecessary, using arguments of O'Grady in 
characteirstic $0$ and recent work of Langer in positive 
characteristic.  Unfortunately, we have not written out the proofs, so 
we must exclude those results from our treatment.

The reader will observe throughout this section evidence for our 
meta-theorem (``All phenomena which occur for moduli spaces of 
semistable sheaves on surfaces also occur for moduli spaces of 
semistable twisted sheaves'').  Unlike the case of curves, the 
evidence in this case is purely behavioral and not attributable to any 
direct comparison of the twisted and untwisted situations.

\subsubsection{Discriminants and dimension estimates}\label{S:discrim}

\begin{defn}\label{D:discriminant} Let $X$ be a smooth projective surface and $\ms X\to X$ a 
$\m_{n}$-gerbe.  Given a coherent $\ms X$-twisted sheaf $\ms F$ of rank $r$, the 
\emph{discriminant\/} of $\ms F$ is the quantity 
$$\Delta(\ms F):=\deg(2rc_{2}(\ms F)-(r-1)c_{1}(\ms 
F)^{2})\in\Z.$$
\end{defn}
\begin{proof}[Proof that $\Delta(\ms F)\in\Z$]  Since $X$ is smooth 
and projective, 
\ref{C:gabber de jong} shows that $\ms F$ has a finite global resolution 
by locally free twisted sheaves.  A formal calculation (in $K^{0}$) shows that 
$$\deg\chern^{\vee}(\ms F)\chern(\ms 
F)\td_{X}=\sum_{i=0}^{2}(-1)^{i}\dimext^{i}(\ms F,\ms F),$$ where 
$\chern^{\vee}(\ms F)_{i}=(-1)^{i}\chern(\ms F)_{i}$.  Another formal 
calculation shows that $$\deg\chern^{\vee}(\ms F)\chern(\ms F)=\rk(\ms 
F)^{2}-\Delta(\ms F).$$  We thus conclude that $\chi(\ms F,\ms 
F)=\Delta(\ms F)-(\rk(\ms F)^{2}-1)\chi({\ms O}_{X})$.  (Such ``formal 
calculations'' show at least that the Chern character and all related 
results -- Grothendieck-Hirzebruch-Riemann-Roch, discriminant 
calculations, etc. -- can be extended to $K^{0}$, hence to the 
homotopy category of strict perfect complexes.  When $X$ is 
projective, any perfect complex on $X$ admits a left resolution by a 
strict perfect complex.  This can be used to show that in fact such 
formal calculations apply to objects of $\D(X)_{\text{parf}}$, the 
derived category of perfect complexes.)
\end{proof}
When $\ms F$ is locally free, $\Delta(\ms F)=c_{2}(\send(\ms F))$.  
(More generally, using the remarks above, one has $\Delta(\ms F)=c_{2}(\rsend(\ms F))$.)  
The discriminant plays an important role in the behavior of the moduli space.

\begin{lem} The discriminant is locally constant in flat families: 
given an $S$-flat family of coherent twisted sheaves $\ms F$ on $X\times 
S$ with $S$ connected, the number $\Delta(\ms F_{s})$ is constant for 
all (geometric) points $s\in S$.
\end{lem}
\begin{proof} Implicit is the statement that $\Delta(\ms F)$ may be 
computed after making any base field extension, which is clear.  One 
easy way to see that $\Delta(\ms F_{s})$ is locally constant in our 
case is to (locally on $S$) resolve $\ms F$ by a complex of locally 
free twisted sheaves $\ms V^{\bullet}\to\ms F$, use the fact that 
$\Delta(\ms F)=c_{2}(\shom^{\bullet}(\ms V^{\bullet},\ms V^{\bullet}))$ 
and then use the fact that intersection products and geometric Hilbert 
polynomials are constant in a flat family (\ref{P:hilb poly 
constant}).  Another proof is based on the equality $\Delta(\ms 
F)=\chi(\ms F,\ms F)+(\rk(\ms F)^{2}-1)\chi({\ms O}_{X})$ and the 
semicontinuity theorems for higher Exts (whose methods are 
demonstrated somewhat in \ref{P:loc and constr}(3)).
\end{proof}

In fact, when the determinant is fixed, it is equivalent to specify 
$\Delta$, $P^{g}$, or $c_{2}$.  
Since we will usually fix a determinant in what follows, this means we 
can use any of these surrogates to divide the moduli problem into 
clusters of connected components.  

Recall that the deformation theory of $\Tw^{ss}_{\ms X/k}(n,P)$ at 
a point $[\ms F]$ is governed by the vector spaces $\ext^{1}(\ms F,\ms F)$ and 
$\ext^{2}(\ms F,\ms F)$, while the deformation theory with fixed 
determinant is determined by $\ext^{1}(\ms F,\ms F)_{0}$ and 
$\ext^{2}(\ms F,\ms F)_{0}$ (\ref{S:deformations of twisted sheaves}), where the subscript $0$ denotes 
traceless elements (\ref{S:equi}).  We can use this to 
estimate the dimension of $\Tw^{ss}_{\ms X/k}(n,L,P)$.  
We remind the reader of a well-known lemma, whose proof may be 
extracted from Schlessinger's thesis \cite{schl} and which is written 
up explicitly in 2A.11 of \cite{h-l}.

\begin{lem}\label{L:hull dimension} Let $k$ be a field and $F:Art_{k}\to\operatorname{Set}$ a 
functor with a hull $R$.  If the embedding dimension $\dim_{k}\mf 
m_{R}/\mf m_{R}^{2}=d$ and $F$ has an obstruction theory with values 
in an $r$-dimensional vector space $\mc O$, then $d\geq\dim R\geq d-r$.
\end{lem}

\begin{prop} Suppose $\ms F$ is a semistable $\ms X$-twisted sheaf of 
rank $n$, geometric Hilbert polynomial $P$, and determinant $L$.  Given an algebraic stack $\ms M$ containing $\ms F$ as a 
point, write $\dim_{\ms F}\ms M$ for the dimension of the miniversal 
deformation space of $\ms F$.
\begin{enumerate}
    \item[(i)] $\dimext^{1}(\ms F,\ms F)\geq\dim_{\ms F}\Tw^{ss}_{\ms 
    X/k}(n,P)\geq 
    \dimext^{1}(\ms F,\ms F)-\dimext^{2}(\ms F,\ms F);$

    \item[(ii)] $\dimext^{1}(\ms F,\ms F)_{0}\geq\dim_{\ms 
    F}\Tw^{ss}_{\ms X/k}(n,L,P)\geq 
    \dimext^{1}(\ms F,\ms F)_{0}-\dimext^{2}(\ms F,\ms F)_{0}.$
\end{enumerate}
In both cases, the moduli stack is a local complete intersection at 
$\ms F$ if the lower bound is achieved and formally 
smooth at $\ms F$ if and only if the upper bound is achieved.
\end{prop}
\begin{proof} This is an application of the results of 
\ref{S:deformations of twisted sheaves} and \ref{S:equi} 
along with \ref{L:hull dimension}.
\end{proof}

\begin{defn} Given a semistable twisted sheaf of rank $n$, 
geometric Hilbert polynomial $P$, and determinant $L$, the 
\emph{expected dimension\/} of $\Tw^{ss}(n,L,P)$ at $\ms F$ 
is the quantity 
$$\expdim_{\ms F}\Tw^{s}(n,L,P):=\dimext^{1}(\ms F,\ms F)_{0}-\dimext^{2}(\ms F,\ms F)_{0}.$$
\end{defn}
\begin{lem}\label{L:dimension estimate} The expected dimension at stable 
points is independent of the choice of $\ms 
F\in\Tw^{s}(n,L,P)$ and is equal to $\Delta(\ms 
F)-(n^{2}-1)\chi({\ms O}_{X})$.  The expected dimension jumps at properly 
semistable points.  There is a constant $\beta_{\infty}$ such that 
for all points $\ms F\in\Tw^{s}_{\ms X/k}(n,L,P)_{k}$, 
$$\expdim\Tw^{s}(n,L,P)\leq\dim_{\ms 
F}\Tw^{s}(n,L,P)\leq\expdim\Tw^{s}(n,L,P)+\beta_{\infty}.$$
\end{lem}
\begin{proof} The formula for the expected dimension follows from the identity
$$-\dimhom(\ms F,\ms F)_{0}+\dimext^{1}(\ms F,\ms F)_{0}-\dimext^{2}(\ms F,\ms 
F)_{0}=\chi({\ms O}_{X})-\sum_{i=0}^{2}(-1)^{i}\dimext^{i}(\ms F,\ms F)$$
and formal calculations.  One uses the fact that stable sheaves $\ms 
F$ are simple ($\End(\ms F)=k$), which immediately implies that $\hom(\ms 
F,\ms F)_{0}=0$, and the rest follows.  (The given identity also uses the trace map 
splitting \ref{L:trace splitting} and thus requires that the rank of $\ms F$ be 
relatively prime to the characteristic of $X$.)
Details of this type of calculation may be found in \cite[4.5, 6.1, 
8.3]{h-l}.  Since $\Delta(\ms F)$ is determined by the determinant and 
Hilbert polynomial, we see that this 
is independent of the stable twisted sheaf $\ms F$.

The jumping of the expected dimension at properly semistable points 
comes from the fact that $\hom(\ms F,\ms F)_{0}$ need not be zero.  
The identity above shows that $$\expdim_{\ms F}\Tw^{ss}_{\ms 
X/k}(n,L,P)-\hom(\ms F,\ms F)_{0}$$ is constant, so the expected 
dimension jumps whenever there are traceless endomorphisms (i.e., 
infinitesimal automorphisms acting trivially on the determinant).

The last inequality follows immediately from the fact that there is a constant 
$\beta_{\infty}$ such that for all semistable twisted sheaves of 
rank $r$ with fixed discriminant (and no restrictions on Chern 
classes if $\ch X=0$), $\dimext^{2}(\ms F,\ms F)_{0}\leq\beta_{\infty}$.
In characteristic 0, this follows easily (using the methods of 
section \ref{S:restriction theorems}) from the Le 
Potier-Simpson estimate and the fact that the endomorphism sheaf of a 
semistable sheaf is semistable \cite[4.5.7]{h-l}, a fact which does not 
hold in positive characteristic.  In general, this is slightly 
subtle (whence the restriction on the discriminant, which is not 
present in characteristic 0) and will be 
proven in \ref{L:universal constant} below.
\end{proof}

\subsubsection{Preparation for restriction theorems}\label{S:preparation}

In this section, we study the following question: given a 
$\m_{n}$-gerbe $\ms X\to X$, how can one construct a finite flat cover 
$Y\to\ms X$ with $Y$ smooth and such that a general member of 
$|{\ms O}_{X}(1)|$ has smooth preimage on $Y$?  Slight complications 
arise in positive characteristic, but this is nonetheless always 
possible.  In the end of the section, we recall a result of Artin and 
de Jong which can be used to ensure that $\deg Y/\ms X=\ind\ms X$.

Let $X$ be a smooth projective surface over an algebraically closed 
field $k$ and $P\to X$ a Brauer-Severi 
variety of relative dimension $n$.  Note that the Brauer class of $P$ 
is split by $P$, hence by any subscheme of $P$.  Choose a projective 
embedding of $P$.  Let $D\subset X$ be a smooth divisor.  We start 
with a lemma about generic hyperplane sections of a Brauer-Severi 
scheme, which is essentially a 
refinement of a special case of a lemma of Vistoli and 
Kresch \cite{vistoli-kresch}.

\begin{lem}\label{L:cut me down} Let $P\to X$ be a surjective map of smooth 
projective varieties with fibers of equidimension $n$ which is 
generically smooth over $D$.  Let $P\inj\P^{N}$ be a closed 
immersion.  A generic hyperplane section $P_{H}$ of $P$ has the 
following properties: $P_{H}$ is smooth and irreducible, $P_{H}\times_{X}D\subset 
P_{H}$ is an irreducible smooth divisor, $P_{H}\to X$ is surjective and 
generically smooth over $D$ with fibers of 
equidimension $n-1$.
\end{lem}
\begin{proof} Let $\Xi$ be the projective space parametrizing hyperplane sections 
of $P$.  The smoothness of the hyperplane section of $P$ and its 
intersection with the pre-image of $D$ defines an open subset 
$U\subset\Xi$.  Let 
$d\in D(k)$ be a smooth point with smooth fiber $P_{d}\subset P$.  The 
condition that a hyperplane $H\in\Xi$ intersect $P_{d}$ in a smooth 
variety of dimension $n-1$ defines an open subset $V$ in $\Xi$. Let 
$W=U\cap V$.  We claim that the hyperplane sections parametrized by 
$W$ have the properties of the lemma.  Indeed, if $H\in W$, then 
$P_{H}$ and $P_{H}\times_{X}D$ are smooth and irreducible since 
$W\subset U$.  Furthermore, the fiber of $P_{H}\to X$ over $d$ is 
smooth of dimension $n-1$ (and hence also irreducible, incidentally) 
since $H\in V$.  We claim that this forces $P_{H}\to X$ be be 
surjective, generically smooth, with equidimensional fibers.  Indeed, 
we have $\dim P-\dim X=n$, hence $\dim P_{H}-\dim X=n-1$.  If $\im 
P_{H}=I$, then the usual inequalities \cite[\S15]{matsumura} show that 
$\dim P_{H}-\dim I\leq n-1$ (as $\dim (P_{H})_{d}=n-1$).  Thus, $I=X$ 
and $P_{H}\to X$ is surjective.  Applying the identity once more 
shows that any closed fiber has dimension at least $n-1$ at any closed 
point.  Thus, every closed fiber is equidimensional of dimension $n-1$.
\end{proof}

\begin{lem}\label{L:removing etale good} Let $f:C\to \spec K$ be a normal curve over a field.  If 
$S\subset C$ is a closed subscheme which is finite \'etale over $K$ 
and $f$ is smooth along $C\setminus S$, then $f$ is smooth.
\end{lem}
\begin{proof} The scheme $C$ is Noetherian and reduced.  Thus, to show 
that the sheaf $\Omega^{1}_{C/K}$ is locally free of rank $1$, it 
suffices to show that for every point $P\in C$, the $\kappa(P)$-vector space 
$\Omega^{1}_{C/K}\tensor_{C}\kappa(P)$ is 1-dimensional.  For points 
$P\in C\setminus S$, this holds by assumption.  On the other hand, 
given a point $Q\in S$, there is a canonical sequence
$$\mf m_{Q}/\mf 
m_{Q}^{2}\to\Omega^{1}_{C/K}\tensor\kappa(Q)\to\Omega^{1}_{\kappa(Q)/K}\to 0.$$
Since $Q$ is a Weil divisor on a normal separated scheme, it is a 
Cartier divisor and therefore the left-most term is $1$-dimensional 
over $\kappa(Q)$.  Since $\kappa(Q)$ is separable over $K$, the 
right-most term vanishes.
\end{proof}

\begin{lem}\label{L:nodal pencil good} Suppose $Y$ is a smooth surface over an algebraically 
closed field $k$ and $D,D'\in|{\ms O}(1)|$ are 
very ample divisors such that $D$ is at worst nodal and $D$ and $D'$ intersect 
transversely.  Then the general member of the pencil spanned by $D$ 
and $D'$ is smooth.
\end{lem}
\begin{proof} Write $\widetilde Y\to\P^{1}$ for the total space of 
the pencil.   The non-smooth locus of $\widetilde Y\to\P^{1}$ has 
the property that it is unramified over $\P^{1}$ at $[D]$.  Indeed, 
the fiber over $[D]$ is a nodal curve, so this follows from the 
standard construction of the scheme structure on the non-smooth locus 
using Fitting ideals \cite[2.21]{dejong-alterations}.  (Really, this just comes down to showing 
that the relative differentials of a node are supported precisely on 
the node with length 1.)  
Thus, 
all components of the non-smooth locus which intersect the generic fiber must be generically \'etale over $\P^{1}$.  This 
implies that any non-smooth points of the generic fiber have separable 
residue fields over $k(\P^{1})$.  The result follows by 
\ref{L:removing etale good}.

An alternative (well-known) argument (rather than exploit the scheme structure of 
the non-smooth locus) comes from the miniversal deformation space of a 
node.  Completing $\widetilde Y$ with respect to the uniformizing 
parameter of $[D]$ at one of the nodes over $[D]$ 
yields an effective formal deformation of the node over $k\[t\]$ with 
the property that the total space is regular.  On the other hand, the 
versal deformation of a node is isomorphic to 
$k\[\xi,X_{0},X_{1}\]/(X_{0}X_{1}-\xi)$ parametrized by $k\[\xi\]$.  
(In other words, given any family of curves $\ms C\to S$ with a node 
$c$ in a closed fiber $\ms C_{s}$, there is a map 
$k\[\xi\]\to\widehat{{\ms O}}_{S,s}$ such that $\widehat{{\ms O}}_{\ms 
C,c}=\widehat{{\ms O}}_{S,s}\ctensor_{k\[\xi\]}k\[\xi,X_{0},X_{1}\]/(X_{0}X_{1}-\xi)$.)
Thus, there is some map $k\[\xi\]\to k\[t\]$ giving rise to 
$\widehat{\widetilde Y}$, and the condition of regularity forces 
$\xi$ to map to $ut$ where $u$ is a unit of $k\[t\]$.  This shows 
that the generic fiber is smooth in the generizations of the node.  
(Indeed, the compatibility properties of $\Omega^{1}$ allow us to 
assume that the base is $k\[t\]$.  Now the map from $\widetilde Y_{\text{node}}$ to its completion 
is regular as $\widetilde Y$ is excellent.  Thus, the map from the 
generic fiber of $\widetilde Y_{\text{node}}$ to the 
generic fiber of the completion is regular.  But given a regular map $A\to B$ 
of Noetherian rings over a field, it follows 
that $A$ is geometrically regular over the field if and only if $B$ is 
geometrically regular over the field.  This applies to our situation 
to show that $\widetilde Y_{\text{node}}$ is smooth over $k\(t\)$.)
\end{proof}

\begin{prop}\label{P:finite happy cover} There exists a smooth 
subvariety $Y\subset P$ which is finite flat generically \'etale
over $X$ such that for every $n$, the pullback of a general member 
of $|{\ms O}_{X}(n)|$ to $Y$ is smooth.
\end{prop}
\begin{proof} By \ref{L:cut me down} and induction, we may carry this out 
for $n=1$.  (Indeed, once a single smooth member pulls back to a 
smooth divisor, it will hold for a general smooth member.  This 
follows from a consideration of the pullback of the incidence 
correspondence for ${\ms O}(1)$ on $X$ to $Y$ and the standard results 
about generization of smoothness in a flat family.)
Let $f:Y\to X$ be the restriction of the projection $P\to X$.  We 
will show that once the result holds for $n=1$, it holds for all $n$.  
Indeed, once it holds for $n=1$, there is a dense open in 
$|{\ms O}_{X}(1)|$ of smooth members whose preimages in $Y$ are smooth.  Given 
$n$, we may choose 
$n$ such general members which intersect transversely away from the 
branch curve of $f$.  Call such a resulting nodal divisor $D_{n}$.  Choose 
$D'_{n}\in|{\ms O}_{X}(n)|$ which is at worst nodal and intersects $D_{n}$ 
transversely away from the branch curve.  Then the pencil generated 
by $D_{n}\times_{X}Y$ and $D_{n}'\times_{X}Y$ satisfies the conditions 
of \ref{L:nodal pencil good}, hence has smooth general member.  (So 
does the pencil generated by $D_{n}$ and $D'_{n}$ on $X$.)
\end{proof}

\begin{para}\label{R:smaller degree} It is likely that the cover produced by \ref{P:finite 
happy cover} is not ideal, in the sense that the degree of the map 
$Y\to X$ is far too large.  (We can know this ``abstractly'' because 
the proof of \ref{P:finite happy cover} is so easy.)  In fact, if $P\to X$ is a Brauer-Severi 
variety of relative dimension $d-1$ representing a Brauer class of 
index $d$, the lowest degree 
for the map $Y\to X$ arising in \ref{P:finite happy cover} will be $d^{d-1}$.  On the 
other hand, we know by results of Artin and de Jong 
\cite[\S8.1]{artin-dejong} that there will be 
a finite flat surjection $Y'\to X$ from a smooth surface to $X$ of 
degree $d$.  Moreover, using methods similar to those already shown 
above, one can actually find such a cover such that a general member 
of $|{\ms O}(1)|$ has smooth preimage on $Y'$.  The interested reader can 
find a few more details in \cite{period-index-paper}.  Given the Azumaya algebra $\ms A$ 
on $X$ of degree $d$ representing $P\to X$, the idea of Artin and de 
Jong's construction is to let $Y'$ be determined by the 
``characteristic polynomial'' of a general section of $\ms A\tensor 
L$, where $L$ is a sufficiently ample invertible sheaf on $X$.   In 
other words, thinking of $\ms A$ as a form of $\M_{d}({\ms O}_{X})$, one 
can see that for any invertible sheaf $L$, the reduced norm yields an 
algebraic map $\ms A\tensor L\to\operatorname{Sym}^{\bullet}L$ with 
image in the polynomial sections of degree $d$.  The locus of zeros 
of such a polynomial function on $L^{\vee}$ gives a finite cover of 
$X$ of degree $d$, which in this case will factor through the gerbe $\ms X$.  Taking 
general $L$ and a general section yields a smooth such cover.  If, in 
addition to the arguments of Artin and de Jong, one pays attention to 
the generic branching behavior of such a cover (which may require 
making $L$ more ample), one gets the following statement.
\end{para}

\begin{prop}\label{C:artin de jong cover} Given a smooth surface $X$ and a $\m_{n}$-gerbe $\ms X$, 
if there is a locally free $\ms X$-twisted sheaf of rank $d$ then 
there is a finite flat surjection of smooth surfaces $Y\to X$ of 
degree $d$ such that 
\begin{enumerate}
    \item there exists an invertible $\ms X\times_{X}Y$-twisted 
    sheaf, and

    \item for every very ample invertible sheaf ${\ms O}(1)$ on $X$, a general 
member has smooth preimage in $Y$. 
\end{enumerate}
\end{prop}

We will use this in the sequel to make better numerical estimates.

\begin{remark} The method of \ref{L:cut me down} and \ref{P:finite 
happy cover} seems likely to generalize to higher dimensional 
varieties $X$.  The only difficulty in the argument is in ensuring 
that general members of ${\ms O}(n)$ have smooth preimages once it is true 
for $n=1$.  For the applications envisioned, it is in fact sufficient 
that such divisors have normal preimages, which may be easier to 
arrange.  In either case, it seems likely that a similar (more 
subtle) analysis of the behavior of a pencil with a fiber consisting 
of a divisor with sufficiently transverse crossings will yield a
geometrically normal generic fiber, which is enough for 
applications.  In other words, there would result a finite flat cover 
$Y\to X$ by a smooth variety such that the general member of ${\ms O}(n)$ 
has normal integral preimage.

On the other hand, the method of Artin and de Jong seems 
harder to generalize directly, because their construction can produce 
singularities in the cover in codimension 3.  Nevertheless, if one is willing to allow $Y$ to be normal, it 
is conceivable that a refinement of their method could yield a finite 
flat covering with a better degree and all of the properties necessary 
to carry out analogues of our proofs below.  This of course has the 
advantage of yielding a better numerical answer, hence more effective 
bounds, but at the present time it is not clear if having a non-smooth 
cover $Y$ is compatible with the methods used here.  We leave this investigation to future work(ers).
\end{remark}

\subsubsection{Restriction theorems and the Bogomolov inequality}\label{S:restriction theorems}

Classically, Mehta and Ramanathan proved that the restriction of a slope-semistable sheaf to 
a general sufficiently ample divisor is again slope-semistable.  An effective version (which specifies 
what ``sufficiently'' means) was first proven in characteristic 0 by 
Bogomolov; a recent paper of Langer \cite{langer} gives a much more general 
statement, valid in all characteristics.  Using Langer's results, we will give twisted 
versions of these theorems in this section.  One of the 
(future) uses of these theorems is to construct the Uhlenbeck 
compactification of the space of twisted sheaves (and then, hopefully, 
the space of $\PGL_{n}$-bundles).  We also use the work of Langer to 
provide a twisted Bogomolov inequality, recovering earlier work of 
Artin and de Jong \cite[\S7.2]{artin-dejong} in the context of Azumaya algebras.  
Throughout, $X$ is a smooth projective surface over an algebraically 
closed field $k$.

\begin{para} We first study restriction theorems.  Fix a 
$\m_{n}$-gerbe $\ms X\to X$.
\end{para}
\begin{lem}\label{L:take cover} Let $f:Y\to X$ be a finite separable morphism of smooth 
surfaces.  A torsion free coherent twisted sheaf $\ms F$ on $X$ is 
$\mu$-semistable if and only if $f^{\ast}\ms F$ is $\mu$-semistable.
\end{lem}
\begin{proof} This may be found in \cite[3.2.2]{h-l}.
\end{proof}

\begin{lem} Let $f:Y\to X$ be a finite flat map of smooth surfaces of 
degree $d$, ${\ms O}_{X}(1)$ a very ample invertible sheaf on $X$, $\ms 
X\to X$ a $\m_{n}$-gerbe, $n\in k^{\times}$.  Write $\ms Y=\ms 
X\times_{X}Y$.  The diagram
$$\xymatrix{A^{2}(\ms X)_{\Q}\ar[r]^{f^{\ast}}\ar[d]_{\deg} & A^{2}(\ms 
Y)_{\Q}\ar[d]^{\deg}\\
\Q\ar[r]^{d} & \Q}$$
commutes.  In particular, given a torsion free $\ms X$-twisted 
sheaf $\ms F$, one has $$\mu_{f^{\ast}{\ms O}_{X}(1)}(f^{\ast}\ms 
F)=d\mu_{{\ms O}_{X}(1)}(\ms F)$$ and $\Delta(f^{\ast}\ms F)=d\Delta(\ms F)$.
\end{lem}
\begin{proof} It suffices to show that the similar diagram with $X$ 
and $Y$ in place of $\ms X$ and $\ms Y$ commutes.  That can be seen 
easily on the level of 0-cycles.
\end{proof}

By \ref{C:artin de jong cover}, we may fix a finite map $f:Y\to X$ of 
smooth surfaces of degree $d=\ind(\ms X)$ with the property that a 
general member of \emph{any\/} very ample linear system on $X$ has smooth 
preimage in $Y$, and such that there is an invertible twisted sheaf 
$\ms L$ on $Y$.  (The salient feature of such a cover is that the ramification curve is generically unramified over the branch curve.)  Fix a very ample linear system ${\ms O}_{X}(1)$ on $X$, with 
associated divisor class $H$.  
Following Langer \cite{langer}, we choose a nef divisor $A$ 
on $Y$ such that $T_{Y}(A)$ is globally generated, and we set 
$$\beta_{r}=\left(\frac{r(r-1)}{p-1}AH\right)^{2},$$
where we assume that $\ch X=p$.  This depends upon $A$, and it is 
slightly unfortunate that this fact is not recorded in the notation.  (When $\ch X=0$, set $\beta_{r}=0$.)
Our method has the perverse consequence that 
effective restriction theorems are easier to prove than generic 
restriction theorems.

\begin{prop}[Twisted Langer]\label{P:twisted langer} Let $\ms E$ be a torsion free $\ms X$-twisted 
sheaf of rank $r$.  Let $D\in |kH|$ be a smooth divisor such that 
$\ms E_{D}$ is torsion free and $D\times_{X}Y$ is smooth.
\begin{enumerate}
    \item[(i)] If $\ms E$ is $\mu$-stable and 
    $$k>\frac{r-1}{r}\ind(\ms X)\Delta(\ms E)+\frac{1}{\ind(\ms 
    X)\deg_{H}(X)(r-1)}+\frac{(r-1)\beta_{r}}{\ind(\ms 
    X)\deg_{H}(X)}$$
    then $\ms E_{D}$ is $\mu$-stable.

    \item[(ii)] If $\ms E$ is $\mu$-semistable and all of the 
    Jordan-H\"older factors of $\ms E$ have torsion free restrictions 
    to $D$, and the inequality of (i) holds, then $\ms E_{D}$ is 
    $\mu$-semistable.
\end{enumerate}
\end{prop}
\begin{proof} After twisting by $\ms L^{\vee}$, the pullback of $\ms E$ to $Y$ is naturally identified 
with a torsion free coherent untwisted sheaf $\ms F$, satisfying the 
stability conditions of (i) or (ii).  Furthermore, $\Delta(\ms 
F)=\ind(\ms X)\Delta(\ms E)$ and $\deg_{f^{\ast}H}(Y)=\ind(\ms 
X)\deg_{H}(X)$.  The inequalities reduce to those of Langer's 
effective restriction theorems \cite[5.2 and 5.4]{langer}, whence $\ms 
F_{D}$ is (i) $\mu$-stable or (ii) $\mu$-semistable.  Applying 
\ref{L:take cover}, we see that $\ms E_{D}$ is (i) $\mu$-stable or 
(ii) $\mu$-semistable, as required.
\end{proof}

\begin{remark} It is irritating to have to pay attention to 
$D\times_{X}Y$, as this makes the result quite a bit less effective.  
One might be tempted to see \ref{P:twisted langer} (as we have proven 
it) as an ``effective generic restriction theorem,'' as the integer 
$k$ is effectively bounded, whereas by \ref{C:artin de jong cover} we 
know that a general member of $|kH|$ will have smooth preimage in 
$Y$.  It would be interesting to find a more effective version which does away with the (abstract) selection of a general member of the linear system in favor of a criterion depending upon the geometry of a member.
\end{remark}

\begin{cor}[Twisted Mehta-Ramanathan] If $\ms F$ is a torsion free $\mu$-semistable $\ms 
X$-twisted sheaf then the restriction of $\ms F$ to a general 
sufficiently ample curve $C\subset X$ is $\mu$-semistable.
\end{cor}
\begin{proof} This is immediate from \ref{P:twisted langer} and the 
properties of preimages of divisors ensured by \ref{C:artin de jong cover} 
(or \ref{P:finite happy cover}, which will just change the estimates 
in \ref{P:twisted langer}).
\end{proof}
As promised in \ref{L:dimension estimate}, we prove the existence of the universal constant 
$\beta_{\infty}$ such that $\dimext^{2}(\ms F,\ms 
F)_{0}\leq\beta_{\infty}$ for all $\mu$-semistable $\ms F$ with rank 
$n$ and fixed discriminant $\Delta$.  The notation 
grates slightly with the notation $\beta_{r}$ of this section, but we 
have chosen to retain the notation of both Huybrechts and Lehn 
($\beta_{\infty}$) and Langer ($\beta_{r}$).  In future sections, we 
will not return to the restriction theorems, so $\beta_{r}$ will 
vanish, which makes this annoyance temporary.
\begin{lem}\label{L:universal constant} There exists a constant 
$\beta_{\infty}$ depending only on $X, \ms X, Y, n, H$ and $\Delta$ such that 
for any $\mu$-semistable twisted sheaf $\ms F$ of 
rank $n$ and discriminant $\Delta$, one has $\dimext^{2}(\ms F,\ms 
F)_{0}\leq\beta_{\infty}$.
\end{lem}
\begin{proof} It suffices to prove this after pulling back to $Y$.  
(Indeed, by the obvious twisted Serre duality, one can see that 
$\dimext^{2}(\ms F,\ms F)_{0}=\dimhom(\ms F,\ms 
F\tensor\omega_{X})_{0}$.  Furthermore, 
$f^{\ast}\omega)_{X}\inj\omega_{Y}$, so 
$$\hom(\ms F,\ms F\tensor\omega_{X})_{0}\leq\hom_{Y}(\ms F_{Y},\ms 
F_{Y}\tensor f^{\ast}\omega_{X})_{0}\leq\hom_{Y}(\ms F_{Y},\ms 
F_{Y}\tensor\omega_{Y})_{0}$$
and we may apply Serre duality again on $Y$.)
Thus, we may assume that $\ms F$ is a semistable untwisted sheaf.  We 
can then suppress $Y$ from the notation; the dependence of $\beta_{\infty}$ on $Y$ only comes in the 
form of a $\beta_{r}$ in the formula.  Pushing the formulas given here 
back down to $X$ will result in multiplying each $\deg_{H}(X)$ and 
each $\Delta(\ms F)$ by $\ind(\ms X)$.

In general, we have $\Delta(\send(\ms F))\leq 2n^{2}\Delta(\ms F)$.  
Indeed, $\ms F$ injects into its reflexive hull $\ms F^{\vee\vee}$, 
yielding an injection $\send(\ms F)\inj\send(\ms F^{\vee\vee})$.  It is not hard to see that 
\begin{equation}
    \ell(\send(\ms F^{\vee\vee})/\send(\ms F))\leq n\ell(\ms F^{\vee\vee}/\ms 
    F).      
    \label{eq:length estimate}
\end{equation}
On the other hand \cite[3.4.1]{h-l}, we have 
\begin{equation}
    \Delta(\ms F)=\Delta(\ms F^{\vee\vee})+2n\ell(\ms F^{\vee\vee}/\ms F)    
    \label{eq:delta comparison}
\end{equation}
and similarly for $\send(\ms F)$.  Combining (\ref{eq:length estimate}) 
with (\ref{eq:delta comparison}) for $\ms F$ and for $\send(\ms F)$ 
shows that $\Delta(\send(\ms F))\leq 2n^{2}\Delta(\ms F)$.  Now
a theorem of Langer \cite[5.1]{langer} 
combined with the inequality $\Delta(\send(\ms F))\leq 2r^{2}\Delta(\ms 
F)$ and the fact that $\mu(\send(\ms F))=0$ shows that 
$$\mu_{\text{max}}(\send(\ms F))\leq 2n\deg_{H}(X)\Delta(\ms 
F)+\beta_{r}.$$  Another theorem of Langer \cite[3.3]{langer-moduli} says 
(in the case of surfaces) 
that for any torsion free sheaf of rank $n$ on $X$,
$$h^{0}(X,E)\leq 
n\deg_{H}(X)\binom{\frac{\mu_{\text{max}}(E)}{\deg_{H}(X)}+f(n)+2}{2},$$
where $f(n)=-1+\sum_{i=1}^{n}1/i$.  Combining this with the estimate 
for $\mu_{\text{max}}(\send(\ms F))$ yields a bound for $\hom(\ms 
F,\ms F)$.  Similarly, we get a bound for $\hom(\ms F,\ms 
F\tensor\omega_{X})$ which differs from the first by a constant depending 
only on $X$.  By Serre duality, $$\dimext^{2}(\ms F,\ms F)_{0}=\dimhom(\ms F,\ms 
F\tensor\omega_{X})-\operatorname{h}^{0}(\omega_{X}),$$ 
so we are done.
\end{proof}

\begin{remark} Note that bounding the discriminant does not suffice 
to bound the Hilbert polynomial when the determinant is not fixed.  
Thus, \ref{L:universal constant} is non-trivial.  Of course, when 
working with a fixed determinant and therefore a bounded set of 
sheaves, some constant $\beta_{\infty}$ will exist by virtue of the 
boundedness and the usual semicontinuity theorems for Ext sheaves.  
In characteristic 0 (or for strongly semistable sheaves in general, 
which we will briefly describe below), the dependence upon the discriminant disappears; 
it is not clear to me if this should still be true in positive 
characteristic.
\end{remark}

\begin{para} We can also use the work of Langer and the coverings of 
\ref{C:artin de jong cover} to produce a version of the Bogomolov 
inequality for twisted sheaves.  After defining a notion of Frobenius 
pullback and strict semistability for twisted sheaves, we can use these 
methods to recover a Bogomolov-like inequality first proven by Artin 
and de Jong in the context of Azumaya algebras.  This inequality will 
be important at one point during the study of asymptotic properties 
of the moduli spaces.

We begin by defining a Frobenius map which is appropriate for our 
situation.  First, note that the (absolute) Frobenius can be defined for 
stacks of characteristic $p$.  If $\ms S\to S$ is such a stack (with 
$\ch(S)=\{p\}$), which we may assume split 
as a fibered category, then the Frobenius 
1-morphism $F_{\ms S}:\ms S\to\ms S$ sends a 1-morphism $T\to\ms S$ to the 
composition $\xymatrix{T\ar[r]^{F_{T}} & T\ar[r] & \ms S}$ (and fixes all 
morphisms in fiber categories).
\end{para}
\begin{lem} If $\ms X\to X$ is any stack and $\chi:\ms I(\ms 
X)\to\G_{m}$ is any character, then the Frobenius map $F_{\ms X}$ 
pulls back $\chi$-twisted sheaves to $p$-fold $\chi$-twisted sheaves.
In particular, if $\ms X\to X$ is a $\m_{n}$-gerbe, then  
the Frobenius map $F_{\ms X}:\ms X\to\ms X$ pulls back $\ms X$-twisted 
sheaves to $p$-fold twisted sheaves.
\end{lem}
\begin{proof} Note that the map on the site of $\ms X$ induced by the 
Frobenius is the identity.  In particular, there is a natural 
isomorphism $$F_{\ms X}^{\ast}(\ms I(\ms X))\simto\ms I(\ms X)$$
(as this is true for any sheaf).
It is not hard to see that the composition
$$\ms I\simto F^{\ast}\ms I\simto\ms 
I$$
is equal to the identity, where the left-hand map 
in the composition is the natural map \ref{L:inertia is functorial}.  
Under this identity, given any sheaf $\ms F$ on $\ms X$, the action 
$\ms F\times\ms I\to\ms F$ pulls back under $F_{\ms X}$ to be the 
\emph{same\/} action $\ms F\times\ms I\to\ms F$.  
On the other hand, given any ${\ms O}_{\ms X}$-module $\ms M$, the 
${\ms O}$-structure on $F^{\ast}\ms M$ is given by $\ms 
M\tensor_{{\ms O},F_{{\ms O}}}{\ms O}$, with the map $F_{{\ms O}}:{\ms O}\to{\ms O}$ given by 
sending a section $s$ to $s^{p}$.  Thus, if a section of $\ms I$ acts 
by $\chi$ on $\ms M$, when pulled back it acts by $\chi^{p}$.

The second statement of the lemma is just a restatement of the first one for readers 
who cleverly skipped section \ref{S:prelim}!
\end{proof}

\begin{defn} Let $\ell$ be the order of $p$ in $(\Z/n\Z)^{\times}$.  
The power $F^{\ell}_{\ms X}$ is called the \emph{twisted Frobenius\/} 
of $\ms X$, denoted $F_{\ms X,\tau}$.  The resulting map
$$\xymatrix{\ms X\ar[rr]^{F^{\ell}_{\ms X}}\ar[dr] & & \ms X\times_{X,F^{\ell}_X} X\ar[dl]\\ & X &}$$
is an isomorphism of $\m_{n}$-gerbes which pulls back twisted sheaves 
to twisted sheaves.
\end{defn}
\begin{defn} An $\ms X$-twisted sheaf $\ms F$ is \emph{strictly 
($\mu$-)
semistable\/} if $(F^{q}_{\ms X,\tau})^{\ast}\ms F$ is ($\mu$-) 
semistable for all $q\geq 0$.
\end{defn}
As with untwisted sheaves, it is the strictly $\mu$-semistable 
twisted sheaves which have the best properties.

\begin{prop}[Twisted Langer-Bogomolov Inequality]\label{P:langer bogomolov} Let $\ms E$ be a torsion 
free $\ms X$-twisted sheaf, $Y\to X$ a cover as in \ref{C:artin de 
jong cover} and $\beta_{r}$ as in \ref{P:twisted langer}.
\begin{enumerate}
    \item[(i)] If $\ms E$ is $\mu$-semistable then $\ind(\ms 
    X)^{2}\Delta(\ms E)+\beta_{r}\geq 0.$

    \item[(ii)] If $\ms E$ is strongly $\mu$-semistable then 
    $\Delta(\ms E)\geq 0$
    
    \item[(iii)] If $\rk(\ms E)=\ind(\ms X)$ then $\Delta(\ms E)\geq 0$.
\end{enumerate}
\end{prop}
\begin{proof} Parts (i) and (ii) follow immediately from Langer's 
version of the Bogomolov inequality \cite[3.2]{langer} (which is our statement if 
$\ind(\ms X)=1$) and \ref{L:take cover}.  Part (iii) follows from the 
fact that if $\rk(\ms E)=\ind(\ms X)$, then $\ms E$ has no proper 
torsion free submodules of strictly smaller rank, so $\ms E$ is 
$\mu$-stable.  Thus, since the 
rank of $\ms E$ is unchanged by Frobenius pullback, $\ms E$ is 
strongly $\mu$-stable and we may apply (ii).
\end{proof}
\begin{rem} The fact that $n$ is prime to the characteristic figures 
essentially into part (iii).  We see from (i) that in general there 
is still a lower bound for the second Chern class of any Azumaya 
algebra of class $[\ms X]$, depending only upon $\ms X$ (and 
possibly the choice of covering $Y\to X$).
\end{rem}

\begin{cor}[Artin, de Jong {\cite[7.2.1]{artin-dejong}}]\label{C:artin de jong c2} Let $X$ be a smooth projective surface 
with function field $K$, and let $A$ be an Azumaya algebra over $X$ 
such that $A_{K}$ is a division ring of degree prime to the 
characteristic.  Then $c_{2}(A)\geq 0$.
\end{cor}
\begin{proof} Let $\deg A=d$.  There is a $\m_{d}$-gerbe $\ms X\to X$ and a 
locally free $\ms X$-twisted sheaf $\ms V$ of rank $d$ and trivial 
determinant such that $\send(\ms V)=A$.  It is easy to see that 
$c_{2}(A)=2rc_{2}(\ms V)$, so we are done by \ref{P:langer bogomolov}(iii).
\end{proof}

\begin{rem} Artin and de Jong's original proof of \ref{C:artin de jong c2} is not 
very difficult, but in their approach positive characteristic and 
characteristic 0 are treated in completely different ways.   
Our method ``explains'' what is going on in a characteristic free manner.
They must also bound the second Chern class from below by a different 
method before showing it is 0,  
while both things happen at once in our approach (which also applies 
to more general Azumaya algebras with possibly non-division generic 
points).  Finally, our proof 
gives a reason for the failure of \ref{C:artin de jong c2} when the 
characteristic divides the degree, namely the failure of 
strong stability of $\ms V$.  We feel that this is 
another demonstration of the usefulness of working with twisted 
sheaves (and thus thinking sheaf-theoretically).
\end{rem}

\subsubsection{Asymptotic properties for optimal classes}\label{S:asymptotic properties}
In this section we study the behavior of $\Tw^{ss}_{\ms 
X/k}(n,L,c_{2})$ as $\Delta\to\infty$.  We will always work with 
spaces of twisted sheaves with a fixed determinant.  Due to 
inadequacies in the classical theory of semistable sheaves on 
surfaces 
in positive characteristic (currently being ameliorated by Langer), we 
only prove these theorems in the optimal case in all 
characteristics.  For the arithmetic applications of 
\cite{period-index-paper}, this is the only case that is needed.

The approach is essentially that of O'Grady, described beautifully by 
Huybrechts and Lehn in \cite[Chapter 9]{h-l}.  The biggest difference between the 
approach here and their approach is \ref{P:coming together}, which is 
an alternative ending step in the proof of asymptotic irreducibility.  
Other than this, the rest of the proof is essentially identical to the 
classical proof.  In the optimal case, certain better numerical 
estimates can be made, which we present here.  Otherwise, we quote 
the book of \cite{h-l} for certain proofs.  While they were written in an 
untwisted context, they carry over verbatim (as indicated) to the 
twisted (arbitrary characteristic) context.  I believe (but have not 
carefully checked) that in the non-optimal characteristic 0 case, one can carry 
out a similar transcription of the classical proofs.  However, I have 
avoided dealing with $e$-stability and related numerical estimates in 
this work, so the reader should take this belief with a grain of 
salt.  It is likely that the current characteristic-free work of Langer 
(\cite{langer-castelnuovo}) will prove amenable to a twisted transcription.

\emph{Throughout this section, $\ms X\to X$ is an optimal $\m_{n}$-gerbe 
with $n$ prime to the characteristic of the base field $k$\/}.  
Thus, any rank $n$ torsion free twisted sheaf will be $\mu$-stable.  
We will continue to use the notation $\Tw^{ss}$, even though in this 
case there are equalities 
$\Tw(n,L,P)=\Tw^{s}(n,L,P)=\Tw^{ss}(n,L,P)$.  Furthermore, all of 
these stacks are Deligne-Mumford and are gerbes over their moduli 
spaces.  We are therefore free to conflate 
their closed substacks and closed subspaces of their coarse moduli 
spaces; in particular, the dimension theory does not change.

We write $\Tw$ for 
$\Tw_{\ms X/k}$, etc.  We will also use the notation $\Tw(n,L,c)$, 
where $c=c_{2}$, rather than $\Tw(n,L,P)$, where $P$ is the geometric 
Hilbert polynomial.  (By the Riemann-Roch theorem, these are 
equivalent sets of data.)  Finally, as we will always work with fixed 
rank and determinant, we will write $\Tw^{ss}(\Delta)$ for 
$\Tw^{ss}(n,L,c)$, where $\Delta$ is the discriminant.

\begin{para} We first outline the asymptotic properties and their 
proofs.  The statements will be proven in \ref{Para:proofs} below.
\end{para}
\begin{defn} The closed subspace in $\Tw^{ss}(n,L,c)$ parametrizing 
non-locally free twisted sheaves is the \emph{boundary\/}, denoted 
$\partial\Tw^{ss}(n,L,c)$. 
\end{defn}
For any map $T\to\Tw^{ss}(n,L,c)$ 
corresponding to a family of twisted sheaves on $T\times X$, the 
preimage of $\partial\Tw^{ss}$ in $T$ is a closed subspace 
$\partial T$, which we will also call the boundary of $T$.

\begin{defn} A ($\mu$-stable) point $\ms F\in\Tw^{ss}$ is \emph{good\/} if $\ms F$ is 
locally free and $\dimext^{2}(\ms F,\ms F)_{0}=0$.
\end{defn}
(We include the $\mu$-stability so that the reader is aware of the 
general definition.)  In general, we will write $\beta(\ms 
F)=\dimext^{2}(\ms F,\ms F)_{0}$ and $\beta(Z)=\max\{\beta(\ms F)|\ms 
F\in Z\}$ for a substack $Z\subset\Tw^{ss}(\Delta)$.  The good locus 
is the vanishing set for $\beta$.

\begin{lem}\label{L:good open} There is an open substack of good points 
$\Tw^{ss}_{g}(\Delta)\subset\Tw^{ss}(\Delta)$ which is smooth over $k$ 
with smooth moduli space.
\end{lem}
\begin{proof} The openness follows from the semicontinuity properties 
of higher Exts (see \cite{banica} and \ref{P:loc and constr}(3) for 
an example of the method involved in the twisted case).  Smoothness of the stack is well known and 
comes from \ref{L:hull dimension} (which shows that the universal 
deformation space of a point is formally smooth).  Smoothness of the 
moduli space follows from the fact that 
$\Tw^{s}(\Delta)\to\mTw^{s}(\Delta)$ is a $\m_{n}$-gerbe.
\end{proof}

The asymptotic properties of $\Tw^{ss}(\Delta)$ come from an analysis 
of the substacks $\partial\Tw^{ss}(\Delta)$ and 
$\Tw^{ss}_{g}(\Delta)$.  We can first show that sufficiently large 
irreducible closed substack of $\Tw^{ss}(\Delta)$ must intersect 
$\partial\Tw^{ss}(\Delta)$.

\begin{prop}\label{P:finding boundary points} There are constants 
$A_{1}$, $C_{1}$, and $C_{2}$ such that if $\Delta\geq A_{1}$ and if $Z$ is an irreducible 
closed substack of $\Tw^{ss}(\Delta)$ such that  
$$\dim Z>\left(1-\frac{1}{n+2}\right)\Delta +C_{1}\sqrt\Delta+C_{2}$$ 
then $\partial Z\neq\emptyset$.
\end{prop}

Using \ref{P:finding boundary points}, we will then show that as $\Delta$ grows, so does 
the codimension of the complement of $\Tw^{ss}_{g}(\Delta)$.  More 
precisely, we have the following.  Let 
$W=\Tw^{ss}(\Delta)\setminus\Tw^{ss}_{g}(\Delta)$ (as a reduced 
closed substack).

\begin{prop}\label{P:bounding w} There is a constant $C_{3}\geq C_{2}$ 
and a constant $A_{2}\geq A_{1}$ such that for all 
$\Delta\geq A_{2}$,
$$\dim W\leq\left(1-\frac{1}{2n}\right)\Delta+C_{1}\sqrt\Delta+C_{3}.$$
\end{prop}

Thus, the stack will asymptotically become smooth in codimension $1$ and everywhere 
l.c.i.\ of the expected dimension, hence normal.  

\begin{prop}\label{P:generic smoothness} Suppose $\Delta$ satisfies 
   \begin{enumerate}
       \item[(1)] $\Delta>A_{1}$
   
       \item[(2)] 
       $\Delta-(n^{2}-1)\chi({\ms O}_{X})\geq\left(1-\frac{1}{2n}\right)\Delta+C_{1}\sqrt\Delta+C_{3}+2.$
   \end{enumerate}
Then every irreducible component of $\Tw^{ss}(\Delta)$ intersects 
$\Tw^{ss}_{g}(\Delta)$.  In particular, it is generically smooth of 
the expected dimension.  Furthermore, $\Tw^{s}(\Delta)$ is normal and a local 
complete intersection.    
\end{prop}
\begin{proof} The two properties and the fact that 
$\expdim\Tw^{ss}(\Delta)=\Delta-(n^{2}-1)\chi({\ms O}_{X})$ (at any point, 
hence on any irreducible component) shows that the locus of good 
points $\Tw^{ss}_{g}(\Delta)$ is dense in every component of 
$\Tw^{ss}(\Delta)$.  When $\dimext^{2}(\ms F,\ms F)_{0}=0$, one then 
has $$\dim\Tw^{ss}_{g}(\Delta)=\dimext^{1}(\ms F,\ms F)_{0}=\expdim\Tw^{ss}_{g}(\Delta),$$ so 
the stack $\Tw^{ss}(\Delta)$ is generically smooth of the expected in 
every irreducible component, hence at every point.  This implies by 
\ref{L:hull dimension} that $\Tw^{s}(\Delta)$ is a local complete 
intersection.  Furthermore, by condition (2) and \ref{P:bounding 
w}, $\Tw^{ss}(\Delta)$ is regular in codimension 1.  By Serre's 
theorem, $\Tw^{ss}(\Delta)$ is normal.
\end{proof}

Another use of \ref{P:finding boundary points} is in proving that 
$\Tw^{ss}(\Delta)$ is irreducible for sufficiently large $\Delta$.  Suppose $\ms 
F\in\Tw^{ss}(\Delta)$ is good.  This implies that $\ms F$ lies on a 
unique irreducible component of $\Tw^{ss}(\Delta)$.  Any subsheaf 
$\ms F'\subset \ms F$ of finite colength $\ell$ (i.e., such that the 
quotient $\ms F/\ms F'$ has finite length $\ell$) must also be good.  
Indeed, by twisted Serre duality (which is derived from the usual Grothendieck duality for the complex $\rshom(\ms F',\ms F')$ on $X$) and compatibility with trace, 
$\dimext^{2}(\ms F',\ms F')_{0}=\dimhom(\ms F',\ms 
F'\tensor\omega_{X})_{0}$, and similarly for $\ms F$.  Furthermore, 
taking the reflexive hull gives a natural injection $\hom(\ms F',\ms 
F'\tensor\omega_{X})_{0}\inj\hom(\ms F,\ms F\tensor\omega_{X})_{0}$.

\begin{lem} $\Delta(\ms F')=\Delta(\ms F)+2n\ell$.
\end{lem}
\begin{proof} This reduces to showing that $c_{2}(\ms F/\ms 
F')=\ell$, which itself reduces to showing that a twisted sheaf $\ms 
Q$ of length 1 has $c_{2}(\ms Q)=1$.  This follows from the 
twisted Hirzebruch-Riemann-Roch theorem \ref{P:rr} applied to the inclusion 
of $\operatorname{Supp}(\ms Q)$ in $\ms X$, along with a trivial 
calculation when $\ms X$ is a $\m_{n}$-gerbe over a geometric point.
\end{proof}

Thus, $\ms F'$ lies on a unique irreducible component of 
$\Tw^{ss}(\Delta+2n\ell)$.    It is trivial that 
every locally free twisted sheaf $\ms F$ contains a colength 1 subsheaf $\ms 
F_{1}$.  Let $\Lambda_{\Delta}$ denote the set of 
irreducible components of $\Tw^{ss}(\Delta)$.

\begin{lem}\label{L:map of components} Suppose $\Delta$ satisfies the 
conditions of \ref{P:generic smoothness}.  The map sending a good 
twisted sheaf $\ms F$ to $\ms F_{1}$ yields a well-defined map 
$\phi:\Lambda_{\Delta}\to\Lambda_{\Delta+2n}$.
\end{lem}
\begin{proof} It follows from \ref{L:irred quot} that the irreducible component containing $\ms 
F_{1}$ is independent of the choice of $\ms F_{1}$.
\end{proof}

The idea behind the proof of irreducibility of $\Tw^{ss}(\Delta)$ for 
large $\Delta$ is to show that $\phi$ is eventually surjective, and 
that any two points are eventually brought together under an iterate 
of $\phi$.

\begin{prop}\label{P:refined generic smoothness} There is a 
constant $A_{3}$ such that for all $\Delta\geq A_{3}$, the following 
hold.
\begin{enumerate}
    \item[(1)] Every irreducible component of $\Tw^{ss}(\Delta)$ 
    contains a locally free good twisted sheaf.

    \item[(2)] Every irreducible component of $\Tw^{ss}(\Delta)$ 
    contains a point $\ms F$ such that both $\ms F$ and $\ms 
    F^{\vee\vee}$ are good and $\ell(\ms F^{\vee\vee}/\ms F)=1$.
\end{enumerate}
\end{prop}

\begin{thm}\label{T:irred} There is a constant $A_{4}$ so that for 
all $\Delta\geq A_{4}$, the stack $\Tw^{ss}(\Delta)$ is irreducible.
\end{thm}
\begin{proof} By \ref{P:refined generic smoothness}(2), for 
$\Delta\geq A_{3}$ the map 
$\phi:\Lambda_{\Delta-2n}\to\Lambda_{\Delta}$ is surjective.  We wish to 
show that this implies that $\Lambda_{\Delta}$ is eventually a 
singleton.  In the twisted case, there is a slight wrinkle, as $c_{2}$ 
need not be an integer.  Thus, not all discriminants are congruent 
modulo $2n$.  
However, we do know that $\Delta$ is always 
an integer.  Consider the sequences of surjections
$$\xymatrix{\Lambda_{\Delta}\ar[r] & \Lambda_{\Delta+2n}\ar[r] & 
\Lambda_{\Delta+4n}\ar[r] & \cdots\\
\Lambda_{\Delta+1}\ar[r] & \Lambda_{\Delta+1+2n}\ar[r] & 
\Lambda_{\Delta+1+4n}\ar[r] & \cdots\\
\vdots & & & \\
\Lambda_{\Delta+2n-1}\ar[r] & \Lambda_{\Delta+2n-1+2n}\ar[r] & 
\Lambda_{\Delta+2n-1+4n}\ar[r] & \cdots.\\}$$
For any sufficiently large discriminant $\Delta'$, one of the sequences above will 
contain $\Lambda_{\Delta'}$.  If we show that any two components in 
the first set of the sequence eventually map to the same point, then 
we see that each sequence is eventually singletons, and hence that 
any $\Lambda_{\Delta'}$ is eventually a singleton (for large enough 
$\Delta'$).  

We claim that it is enough to show that given locally 
free $\ms V$ and $\ms W$ of rank $n$ with the same determinant and 
discriminant, there are finite colength subsheaves $\ms V'\subset\ms 
V$ and $\ms W'\subset\ms W$ and an irreducible flat family containing 
both $\ms V'$ and $\ms W'$.  This is not obviously the same as making 
colength 1 subsheaves of locally free good sheaves in each stage.  To 
see that these are the same, note that the irreducibility of the 
twisted Quot scheme shows that we may assume that the supports of $\ms 
V/\ms V'$ and $\ms W/\ms W'$ are finite sets of distinct reduced 
points.  Now suppose given a family of twisted sheaves $\ms 
F$ on $X\times S$.  The $S$-scheme of quotients of $\ms F$ of length 
$\ell$ with supports distinct reduced points disjoint from the 
singular locus of $\ms F$ in each fiber is easily seen to be 
irreducible when $S$ is irreducible (see e.g., the proof of 
\ref{L:irred quot}).  Thus, if $S$ is irreducible, so 
is this scheme of quotients.  So as we let a point move in it, it will 
end up in the same irreducible component of $\Tw^{ss}(\Delta+2n)$.  
Since at each stage \ref{P:refined generic smoothness} implies that 
each successive quotient may be irreducibly connected to a locally 
free sheaf, we see that $\ms V\mapsto \ms V'$ is the $\ell(\ms V/\ms 
V')$th iterate of $\phi$.

We will prove the existence of $\ms V'$ and $\ms W'$ below in 
\ref{Para:coming together}.
\end{proof}

\begin{para}\label{Para:proofs} We now prove everything!  First comes \ref{P:finding boundary points}.
\end{para}

\begin{lem}\label{L:unstable restriction} Let $C\in|{\ms O}(N)|$ be a smooth 
member (for any $N$) and let $\ms C=\ms X\times_{X}C$.  Let $Z\subset\Tw^{ss}(\Delta)$ 
be a closed irreducible substack with $\partial Z=\emptyset$.  If $\dim 
Z>\dim\Tw^{ss}_{\ms C/k}(n,\ms Q_{C})$ then there is a point of $Z$ 
parametrizing an $\ms X$-twisted sheaf $\ms F$ whose restriction to $C$ is unstable.
\end{lem}
\begin{proof} By \ref{P:no curves inside}, we see that if it is 
defined the restriction map $Z\to\Tw^{ss}_{\ms C/k}(n,\ms Q_{C})$ is 
finite.  Thus, if every restriction of a point of $Z$ to $C$ is 
stable, we see that $\dim Z\leq\dim\Tw^{ss}_{\ms C/k}(n,\ms Q_{C})$.
\end{proof}

\begin{prop}\label{P:unstable restriction makes boundary} Let 
$Z\subset\Tw^{ss}(\Delta)$ be a closed irreducible substack.  
Let $C\in|{\ms O}(N)|$ be smooth.  Suppose $Z$ contains a point $[\ms F]$ 
such that $\ms F_{C}$ is unstable.  If 
$$\dim Z>\expdim\Tw^{ss}(\Delta)+\beta_{\infty}+\frac{n^{2}}{4}-\frac{n-1}{2}C(C-K)$$
then $\partial Z\neq\emptyset$.
\end{prop}
\begin{proof} This may be copied almost verbatim from \cite[9.5.4]{h-l}, but omit 
the part about $e$-stability.
\end{proof}

\begin{proof}[Proof of \ref{P:finding boundary points}] This is an 
application of \ref{L:unstable restriction} and \ref{P:unstable restriction 
makes boundary}.  Indeed, these show that if $Z$ is an irreducible 
component such that
$$\dim Z>\dim\Tw^{ss}_{\ms 
C/k}(n,L_{C})=\frac{n^{2}-1}{2}(N^{2}H^{2}+NKH)$$
and
$$\dim 
Z>\Delta-(n^{2}-1)\chi({\ms O}_{X})+\beta_{\infty}+\frac{n^{2}}{4}-\frac{n-1}{2}C(C-K)$$
then $\partial Z\neq\emptyset$.  We seek a function of $\Delta$ which 
is greater than both right-hand sides for large $\Delta$ (and some 
choice of $N$) but which is smaller than 
$\Delta-(n^{2}-1)\chi({\ms O}_{X})$ by an amount which grows without bound 
as $\Delta$ increases.  (The second condition becomes necessary when 
trying to make the codimension of $W$ high.)  For the purposes of the 
present work, we do not make any attempt to be especially effective; 
this will make things easier.  Letting $N\sim c\sqrt\Delta$ and examining 
the resulting inequalities for that value of $N$ leads one to choose 
$c$ with 
$$c^{2}<\frac{2}{(n+1)(n-1)H^{2}}$$
to ensure that the ``leading term'' (coefficient of $\Delta$) of the top 
line is larger than that of the bottom line and less than $\Delta$.  
As we let $\Delta$ grow, this will eventually produce positive 
integers for $N$, and working through the arithmetic shows that there 
will be a function $f(\Delta)=C_{1}\sqrt\Delta+C_{2}$ such that for 
$N\sim C_{1}\sqrt{\Delta}$, the inequalities are satisfied and 
$f(\Delta)<\Delta-(n^{2}-1)\chi({\ms O}_{X})$.  Then any $Z$ 
with $\dim Z> f(Z)$ will satisfy both \ref{L:unstable restriction} 
and \ref{P:unstable restriction makes boundary} and have dimension 
strictly smaller than the expected dimension of $\Tw^{ss}(\Delta)$.  
For a similar argument, see \cite[pp.\ 209-210]{h-l}.
\end{proof}

\begin{para} Next come \ref{P:bounding w} and \ref{P:refined generic 
smoothness}.  We begin with some 
preparatory lemmas.
\end{para}

\begin{lem}\label{L:boundary codimension} If 
$\partial\Tw^{ss}(\Delta)\neq\emptyset$ then 
$\operatorname{codim}(\partial\Tw^{ss}(\Delta),\Tw^{ss}(\Delta))\leq 
n-1$.
\end{lem}
\begin{proof} The statement is local on the stack.  Locally on $\Tw^{ss}$, one may 
choose a locally free resolution of the universal object on 
$\Tw^{ss}(\Delta)\times X$ by two sheaves $\phi:L_{1}\to L_{0}\to\ms 
F_{\text{univ}}$ (as surfaces have homological dimension 2).  The 
result follows from studying the locus where the rank of $\phi$ drops, 
which is known from standard theorems about determinant schemes.  See 
\cite[9.2.2]{h-l} for more details.  Note that while the reference given 
for determinantal loci is written over $\C$, the estimates are 
independent of the characteristic.
\end{proof}

We need one more lemma, which is well known.  

\begin{lem}\label{L:double dual} If $\ms F$ is an $S$-flat family of torsion free twisted 
sheaves then the function $s\mapsto\ell(\ms F_{s}^{\vee\vee}/\ms 
F_{s})$ is upper semicontinuous.  If $S$ is reduced and the function 
is constant than the formation of the reflexive hull commutes with 
base change and $\ms F^{\vee\vee}$ is locally free.
\end{lem}
\begin{proof} See e.g.\ \cite[9.6.1]{h-l}.  One uses the fact that a surface has homological 
dimension 2 and that there are locally free resolutions (which is true 
in the twisted setting as well).
\end{proof}

\begin{defn} The \emph{double-dual stratification\/} of 
$\Tw^{ss}(\Delta)$ is given by subsets $$\Tw^{ss}(\Delta)_{\nu}=\{\ms 
F|\ell(\ms F^{\vee\vee}/\ms F)\geq\nu\}.$$
These are closed subsets by \ref{L:double dual}.  For any family of 
torsion free twisted sheaves over $S$, there is an induced 
stratification $S_{\nu}$ by pullback along the classifying map 
$S\to\Tw^{ss}(\Delta)$.
\end{defn}

The most important fact about this stratification is that formation 
of the double dual induces a map 
$$(\partial\Tw^{ss}(\Delta)_{\nu}\setminus\partial\Tw^{ss}(\Delta)_{\nu+1})_{\text{red}}\to\Tw^{ss}_{lf}(\Delta-2n\ell).$$
The fiber over a (locally free) point $\ms F$ is just 
(set-theoretically, at least) $\quot(\ms F,\ell)$.   Let 
$Z\subset\Tw^{ss}(\Delta)$ be a closed irreducible subspace with 
$\partial Z\neq\emptyset$ and $\beta(Z)>0$.  Following section 9.6 of \cite{h-l}, we define a sequence of triples
$$Y_{i}\subset Z_{i}\subset\Tw^{ss}(\Delta_{i})$$
as follows: $\Delta_{0}=\Delta$, $Z_{0}=Z$, and $Y_{i}\subset\partial 
Z_{i}$ is an irreducible component of the maximal open stratum of the 
double-dual stratification of $\partial Z$.  If $\ell$ is the 
constant colength on this stratum, then, as we just remarked, there is 
an induced map $Y_{i}\to\Tw^{ss}(\Delta_{i}-2n\ell)$.  Set 
$\Delta_{i+1}=\Delta_{i}-2n\ell$ and $Z_{i+1}$ equal to the closure of 
the image of $Y_{i}$.  There is some index $m$ such that $\partial 
Z_{m}=\emptyset$ by the twisted Langer-Bogomolov inequality 
$\Delta\geq 0$ \ref{P:langer bogomolov}(iii) (which applies since 
$\ind(\ms X)=n$).

Using \ref{L:irred quot} and \ref{L:boundary codimension}, one 
finds $\dim Z_{i}\geq\dim Y_{i-1}-\ell_{i}(n+1)$ and $\dim 
Y_{i-1}\geq\dim Z_{i-1}-(n-1)$, whence $\dim Z_{i}\geq \dim 
Z_{i-1}-(2n-1)\ell_{i}-1$.  A careful analysis of when equality can 
hold between $\dim Z_{i}$ and $\dim Z_{i-1}-(2n-1)\ell_{i}-1$ (which 
may be found in \cite[pp.\ 211-212]{h-l}) yields 
an inequality
$$\dim Z_{m}-\left(1-\frac{1}{2n}\right)\Delta_{m}\geq \dim 
Z-\left(1-\frac{1}{2n}\right)\Delta-\beta_{\infty}.$$
It is now clear what is going to happen: if $\dim Z$ is too large, 
then $\dim Z_{m}$ is too large, i.e., satisfies \ref{P:finding 
boundary points}, contradicting the fact that $\partial 
Z_{m}=\emptyset$.  The numerical details may be found in \cite[p.\ 
212-213]{h-l}, where it is shown that $$C_{3}:=\max\{C_{2}+\beta_{\infty}, 
A_{1}/2n +2\beta_{\infty}-(n^{2}-1)\chi({\ms O}_{X})\}$$ works in the 
statement of \ref{P:bounding w}.

Finally, the proof of \ref{P:refined generic smoothness} may be 
copied verbatim from \cite[p.\ 213]{h-l}.  

\begin{para}\label{Para:coming together}  As promised in 
\ref{T:irred}, we show that given two (good) locally free twisted sheaves 
$\ms V$ and $\ms W$ with the same rank, determinant, and discriminant, 
there are finite colength subsheaves $\ms V'\subset\ms V$ and $\ms 
W'\subset\ms W$ which belong to a common irreducible family of (good) 
twisted sheaves.
\end{para}

\begin{lem} A general map $\ms V\to\ms W(N)$ is injective with 
cokernel supported on a divisor where it has rank 1 in every fiber.
\end{lem}
\begin{proof} This is a Bertini type theorem.  Over any field, the space of $n\times 
n$-matrices which have rank at most $n-1$ is a divisor in $\M_{n}(k)$ 
with singular locus of codimension $3$ (in the divisor) given by 
matrices of rank at most $n-2$.  Thus, the cone of matrices of rank at 
most $n-2$ has codimension 4 in each fiber, and a standard  
argument shows that on a surface a generic section (for $N$ large 
enough that $\shom(\ms V,\ms W(N))$ is globally generated) will avoid 
this locus.  As the rank drops on a divisor, we are done.
\end{proof}
\begin{cor} A general map $\ms V\to\ms W(N)$ is injective with 
cokernel an invertible twisted sheaf supported on a \emph{smooth\/} curve in 
$|{\ms O}(nN)|$.
\end{cor}
\begin{proof} This involves a similar Bertini argument with the second jet 
bundle of a matrix algebra.  At a point $p$ with local coordinates 
$x$ and $y$, an element of the fiber of this bundle is a matrix 
$M_{0}+xM_{1}+yM_{2}$.  Taking the determinant yields a function 
$f_{0}+xf_{1}+yf_{2}$ (as $x^{2}=y^{2}=xy=0$ in the jet bundle).  In 
order for the determinant to vanish to order at least 2 at the point, 
all three functions $f_{i}$ must vanish.  This defines a ``forbidden 
cone'' of codimension 3 in every fiber (see \cite[8.1.1.6]{artin-dejong} for a 
verification that these conditions are independent), which is greater than the 
dimension of $X$.  The usual argument shows that once the jet bundle 
is globally generated, a general section will miss the forbidden cone 
in each fiber.
\end{proof}

\begin{prop}\label{P:coming together} Let $\ms V$ and $\ms W$ be two locally 
free twisted 
sheaves of rank $n$ with the same determinant and discriminant.  Then there exist torsion 
free twisted sheaves and finite colength inclusions
$\ms V'\subset\ms V$ and $\ms W'\subset\ms W$ (of the same colength) and 
an irreducible flat 
family of twisted sheaves containing $\ms V'$ and $\ms W'$.  If $\ms 
V$ and $\ms W$ are both (good) ($\mu$-)(semi)stable, then there exists an irreducible family 
consisting of (good) ($\mu$-)(semi)stable sheaves.
\end{prop}
\begin{proof} For $N$ sufficiently large, there are extensions
    $$0\to\ms V(-N)\to\ms V\to\ms P\to 0$$
and
    $$0\to\ms V(-N)\to\ms W\to\ms Q\to 0,$$
where $\ms P$ and $\ms Q$ are invertible twisted sheaves on smooth 
curves in the linear system $|{\ms O}(nN)|$.  Furthermore, $\ms P$ and $\ms 
Q$ have the same geometric Hilbert polynomial.  By \ref{C:connecting 
picard mamas}, there is an irreducible variety $S$ (which we may assume 
is affine) and an $S$-flat family 
of twisted sheaves $\ms D$ on $X\times S$ supported on an $S$-flat Cartier divisor 
which interpolates between $\ms P$ (the fiber over $s_{0}\in S(k)$) and $\ms 
Q$ (the fiber over $s_{1}\in S(k)$).  The idea is to make $\ms V'$ 
and $\ms W'$ by taking the inverse image of finite colength 
subsheaves $\ms P(-m)$ and $\ms Q(-m)$ for $m$ sufficiently large 
that we can connect torsion free extensions over this family.

We will use the usual Grauert semicontinuity results for $\ext$ 
spaces to make a connected family interpolating between finite colength 
subsheaves of $\ms P$ and $\ms Q$.  We can do this explicitly quite 
easily as follows.  
Let $\ms F^{\bullet}\to\ms D$ be a finite resolution by a complex of locally 
free twisted sheaves.  (In fact, it will have length at most 2!)  
Twisting $\ms F^{\bullet}$ by a very negative power of ${\ms O}(1)$, we see that 
the perfect complex $\mf C=\hom^{\bullet}(\ms F^{\bullet}(-m),\ms V(-N))$ on $S$ 
universally computes relative $\ext$ spaces.  In other words, for any 
$T\to S$, $\H^{i}(\mf C\tensor_{S}T)\cong\ext^{i}_{X\times_{S}T}(\ms 
D(-m){X\times_{S} T},\ms V(-N)_{X\times_{S} T})$.  Moreover, for large 
enough $m$, it is the case that the function $$s\in S\mapsto\dim\H^{0}(\mf 
C\tensor\kappa(s))$$ is constant and the function $$s\in 
S\mapsto\dim\H^{2}(\mf C\tensor\kappa(s))$$
is the 0 function (by Serre duality).  Standard methods (see 
\cite[III.12]{hartshorne} for example) now show that $\H^{1}(\mf C)$ is a locally 
free sheaf and that for all $f:T\to S$ the natural map 
$f^{\ast}\H^{1}(\mf C)\to\H^{1}(\mf C\tensor_{S}T)$ is an 
isomorphism.  Let $\mathbf V\to S$ be the vector bundle whose sections 
are $\H^{1}(\mf C)$.  Then we have shown that $\mathbf V$ represents 
the functor $T\to S\mapsto\ext^{1}(\ms D(-m),\ms V(-N))$ (for 
sufficiently large $m$).  (In fact, we could have easily shown that 
$\H^{2}(\mf C\tensor T)$ is universally 0 to begin with, by a trivial 
homological dimension calculation, but the method here generalizes 
slightly to higher dimensional ambient varieties.) 
As such, there is a universal extension
$$0\to\ms V(-N)_{\mathbf V\times X}\to\ms E\to\ms D(-m)_{\mathbf 
V\times X}\to 0.$$  Furthermore, once the existence of a vector bundle 
representing $\ext^{1}(\ms D(-m),\ms V(-N))$ is true for $m$, it will be 
true for all $m'>m$.  Thus, to get $\mathbf V$ to have nice 
properties, we can keep enlarging $m$.  

We claim that for sufficiently large $m$, given any $s\in S$ there is a 
non-empty open subset $U_{s}\subset\mathbf V_{s}$ parametrizing 
torsion free extensions.  It is enough to prove that there is a 
single torsion free extension by the openness of purity in families.  
Furthermore, the existence of such a point is stable under increases 
of $m$: if $\ms E$ is torsion free element of $\ext^{1}(\ms D(-m),\ms 
V(-N))$, then the preimage of $\ms D_{s}(-m-m_{0})$ in $\ms E$ gives a torsion 
free element in $\ext^{1}(\ms D_{s}(-m-m_{0}),\ms V(-N))$.
Let $s$ be a point of $S$, so that we are considering extensions 
$\ext^{1}_{X}(\ms D_{s}(-m),\ms V(-N))$.  We are reduced to proving 
that if $m$ is large enough, there is a point of this space 
representing a torsion free twisted sheaf.  Let $\ms E$ be any 
extension with torsion subsheaf $T(\ms E)$.  Since $\ms V(-N)$ is torsion 
free, the intersection $\ms V(-N)\cap T(\ms E)=0$, so $T(\ms E)\inj\ms 
D_{s}(-m)$.  Now consider the situation generically.  Over the local 
ring at the generic point of $\supp\ms D$ there is certainly a torsion free 
extension, so over the complement $W:=X\setminus D$ of some sufficiently ample hyperplane 
section $D\in|{\ms O}(m_{0})|$ there is a torsion free extension
$$0\to\ms V(-N)_{W}\to\ms E_{W}\to\ms D_{s}(-m)_{W}\to 0.$$
Since $W$ is affine, 
this may be realized as a map $\phi_{W}:\ms F^{-1}(-m)_{W}\to\ms V(-N)_{W}$ 
whose composition $\ms F^{-2}(-m)_{W}\to\ms F^{-1}(-m)_{W}\to\ms 
V(-N)_{W}$ is $0$.  Twisting by a high power $c$ of $D$ we find an 
extension $\phi:\ms F^{-1}(-m-cm_{0})\to\ms V(-N)$ of $\phi_{W}$.  It 
is immediate that $\phi$ satisfies the cocycle condition, hence gives 
rise to an extension $\ms E$ which restricts to $\ms E_{W}$ on $W$.  
Since $W$ contains the generic point of $C$ and $\ms D_{s}(-m-cm_{0})$ 
is torsion free, the inclusion $T(\ms E)\inj\ms D_{s}(-m-cm_{0})$ 
implies that $T(\ms E)=0$.  

Therefore, by the openness of the torsion free locus 
and Noetherian induction, we may choose a large $m$ so that the 
torsion free locus of every fiber of $\mathbf V\to S$ is open and 
dense.  This implies that the locus $\mathbf U\subset\mathbf V$ parametrizing 
torsion free sheaves is irreducible.  

Now consider the original points $s_{0}$ and $s_{1}$ over which lie 
$\ms P$ and $\ms Q$.  Choosing a section of ${\ms O}(m)$ yields an 
injection $\ms D(-m)\inj\ms D$, and taking preimages of $\ms P(-m)$ 
and $\ms Q(-m)$, we find finite 
colength subsheaves $\ms V'\subset\ms V$ and $\ms W'\subset\ms W$ 
parametrized by points of $\mathbf U$, hence lying in an irreducible
family of torsion free twisted sheaves.  If $\ms V$ and $\ms W$ are 
(good) ($\mu$-)(semi)stable, then the same is true of $\ms V'$ and $\ms W'$, 
and we are done by the openness of these loci in families and 
irreducibility.
\end{proof}

\begin{rem} This is the key step to proving that the stack 
of semistable twisted sheaves is asymptotically irreducible for 
non-optimal classes as well.  Our proof is sufficiently general to 
work in the general (non-optimal) case.  However, some of the other foundations 
(notably a study of $e$-stability) cannot be carried out in positive 
characteristic yet.  The general characteristic 0 case is likely to work 
precisely as it does in the classical case, but we have not yet 
checked the details.
\end{rem}

\begin{para} We can give a relative version of all of the 
constructions here.  Stack-theoretically, this extension is trivial.  
The GIT construction of Simpson also gives a good global projective 
corepresenting scheme (although in positive characteristic it is no longer 
clear whether or not this is universal on the base).  In the case of 
an optimal class, all points will be stable, so $\Tw^{s}_{X/S}$ is a 
gerbe over a projective scheme, which shows that in this case the 
formation of the coarse moduli space is universal on $S$.  (Universal 
means ``compatible with all base change.''  It is always true that the 
GIT quotient is \emph{uniform\/}, which means that it is compatible 
with flat base change.)
\end{para}

\begin{prop}\label{P:dedekind family} Let $\ms X\to X\to S$ be a 
$\m_{n}$-gerbe on a smooth proper morphism of finite presentation of 
algebraic spaces with geometrically connected fibers of dimension 2, and assume 
that $n$ is invertible on $S$.  Suppose $\ms X$ has optimal geometric fibers.  
The stack $\Tw^{ss}_{X/S}(\Delta)\to S$ is a proper flat local 
complete intersection morphism for large $\Delta$.
\end{prop}
\begin{proof} (This result can also be finagled when the fibers are 
either geometrically optimal or geometrically essentially trivial, 
and is likely to hold completely generally in characteristic 0 by a 
simple extension of our methods.  As above, the general positive 
characteristic case is still in progress.)  This follows from 
\ref{P:generic smoothness} which shows that $\Tw^{ss}_{X/S}(\Delta)$ is 
l.c.i.\ over $S$ (by \ref{L:hull dimension}), combined with the local criterion of flatness.
\end{proof}

\end{document}